\documentclass{amsart}
\usepackage{amssymb}
\usepackage{amsmath}
\usepackage{latexsym}
\usepackage{eepic}
\usepackage{epsfig}
\usepackage{graphicx}
\usepackage{pb-diagram,lamsarrow,pb-lams,amscd}
\usepackage{eufrak}
\usepackage{mathrsfs}
\usepackage[dutch,english]{babel}

\DeclareMathAlphabet{\mathpzc}{OT1}{pzc}{m}{it}
\newtheorem{theorem}{Theorem}[section]
\newtheorem{theorem-definition}[theorem]{Theorem-Definition}
\newtheorem{lemma-definition}[theorem]{Lemma-Definition}
\newtheorem{definition-prop}[theorem]{Proposition-Definition}
\newtheorem{corollary}[theorem]{Corollary}
\newtheorem{prop}[theorem]{Proposition}
\newtheorem{lemma}[theorem]{Lemma}
\newtheorem{cor}[theorem]{Corollary}
\newtheorem{definition}[theorem]{Definition}

\newtheorem{example}[theorem]{Example}

\newenvironment{remark}{\vspace{4pt}\noindent\textbf{Remark.}}{\qed\vspace{4pt}}

\newcommand{\LL}{\ensuremath{\mathbb{L}}}

\newcommand{\N}{\ensuremath{\mathbb{N}}}
\newcommand{\Z}{\ensuremath{\mathbb{Z}}}
\newcommand{\Q}{\ensuremath{\mathbb{Q}}}

\newcommand{\C}{\ensuremath{\mathbb{C}}}

\newcommand{\A}{\ensuremath{\mathbb{A}}}

\newcommand{\X}{\ensuremath{\mathscr{X}}}

\newcommand{\mY}{\ensuremath{\mathfrak{Y}}}

\renewcommand{\C}{\ensuremath{\mathbb{C}}}

\renewcommand{\A}{\ensuremath{\mathbb{A}}}
\renewcommand{\X}{\ensuremath{\mathfrak{X}}}

\renewcommand{\mY}{\ensuremath{\mathfrak{Y}}}
\newcommand{\mZ}{\ensuremath{\mathfrak{Z}}}
\newcommand{\mU}{\ensuremath{\mathfrak{U}}}
\newcommand{\mV}{\ensuremath{\mathfrak{V}}}
\newcommand{\mE}{\ensuremath{\mathfrak{E}}}

\newcommand{\Spec}{\ensuremath{\mathrm{Spec}\,}}
\newcommand{\Spf}{\ensuremath{\mathrm{Spf}\,}}

\numberwithin{equation}{section}

\begin{document}
\title[A trace formula for rigid varieties]{A trace formula for rigid varieties,
and motivic Weil generating series for formal schemes}
\author[Johannes Nicaise]{Johannes Nicaise}
\address{Universit\'e Lille 1\\
Laboratoire Painlev\'e, CNRS - UMR 8524\\ Cit\'e Scientifique\\59655 Villeneuve d'Ascq C\'edex \\
France} \email{johannes.nicaise@math.univ-lille1.fr}
\thanks{The research for this article was partially supported by ANR-06-BLAN-0183}

 \maketitle
 \begin{abstract}
 We establish a trace formula for rigid varieties $X$ over a complete discretely valued field, which relates the set of unramified points on $X$ to the Galois
action on its \'etale cohomology. We develop a theory of motivic
integration for formal schemes of pseudo-finite type over a
complete discrete valuation ring $R$, and we introduce the Weil
generating series of a regular formal $R$-scheme $\mathfrak{X}$ of
pseudo-finite type, via the construction of a Gelfand-Leray form
on its generic fiber. Our trace formula yields a cohomological
interpretation of this Weil generating series.

 When $\mathfrak{X}$ is the formal completion of a morphism $f$ from a smooth irreducible
 variety to the affine line, then its Weil generating series coincides (modulo normalization)
 with the motivic zeta function of $f$. When $\mathfrak{X}$ is the formal completion of $f$ at
 a closed point $x$ of the special fiber $f^{-1}(0)$, we obtain the local motivic zeta function of
 $f$ at $x$. In the latter case, the generic fiber of $\mathfrak{X}$ is the so-called
 analytic Milnor fiber of $f$ at $x$;
 we show that it completely determines the formal germ of $f$ at $x$.
\end{abstract}
\section{Introduction}
Let $R$ be a complete discrete valuation ring,
and denote by $k$ and $K$ its residue field, resp. its field of
fractions. We assume that $k$ is perfect. We fix a uniformizing
parameter $\pi$ in $R$ and a separable closure $K^s$ of $K$, and
we denote by $K^t$ and $K^{sh}$ the tame closure, resp. strict
henselization of $K$ in $K^s$.

 The main goal of the present
paper is to establish a broad generalization of the trace formula
in \cite[5.4]{NiSe}, to the class of special formal $R$-schemes
(i.e. separated Noetherian formal $R$-schemes $\X$ such that
$\X/J$ is of finite type over $R$ for each ideal of definition
$J$; these are also called formal $R$-schemes of pseudo-finite
type). For any special formal $R$-scheme $\X$, Berthelot
constructed in \cite[0.2.6]{bert} its generic fiber $\X_\eta$,
which is a rigid variety over $K$, not quasi-compact in general.
If $\X_\eta$ is smooth over $K$, our trace formula relates the set
of unramified points on $\X_\eta$ to the Galois action on the
\'etale cohomology of $\X_\eta$, under a suitable tameness
condition.

The set of unramified points on $\X_\eta$ is infinite, in general,
but it can be measured by means of the motivic Serre invariant,
first introduced in \cite{motrigid} and further refined in
\cite{NiSe2}. A priori, this motivic Serre invariant is only
defined if $\X_\eta$ is quasi-compact. However, using dilatations,
we show that there exists an open quasi-compact rigid subvariety
$X$ of $\X_\eta$ such that $X(K^{sh})=\X_\eta(K^{sh})$. Moreover,
the motivic Serre invariant of $X$ depends only on $\X_\eta$, and
can be used to define the motivic Serre invariant $S(\X_\eta)$ of
$\X_\eta$. If we denote by $\X_0$ the underlying $k$-variety of
$\X$ (i.e. the closed subscheme of $\X$ defined by the largest
ideal of definition), then this motivic Serre invariant takes
values in a certain quotient of the Grothendieck ring of varieties
over $\X_0$. Our trace formula states that the $\ell$-adic Euler
characteristic of $S(\X_\eta)$ coincides with the trace of the
action of any topological generator of the tame geometric
monodromy group $G(K^t/K^{sh})$, on the graded $\ell$-adic
cohomology space $\oplus_{i\geq
0}H^i(\X_\eta\widehat{\times}_K\widehat{K^t},\Q_\ell)$. As an
auxiliary result, we compute Berkovich' $\ell$-adic nearby cycles
\cite{berk-vanish2} associated to a special formal $R$-scheme
whose special fiber is a strict normal crossings divisor.

Next, we generalize the theory of motivic integration of
differential forms on formal schemes of finite type, to the class
of special formal $R$-schemes, and we extend the constructions and
result in \cite{NiSe} to this setting. In particular, if $k$ has
characteristic zero, we associate to any generically smooth
special formal $R$-scheme its motivic volume, which is an element
of the localized Grothendieck ring of $\X_0$-varieties
$\mathcal{M}_{\X_0}$. Moreover, if $\X$ is a regular special
formal $R$-scheme of pure relative dimension $m$, we associate to
any continuous differential form $\omega$ in
$\Omega^{m+1}_{\X/k}(\X)$ its Gelfand-Leray form $\omega/d\pi$,
which is a section of $\Omega^m_{\X_\eta/K}(\X_\eta)$. This
construction allows us to define for any regular special formal
$R$-scheme $\X$ its Weil generating series, a formal power series
over the localized Grothendieck ring $\mathcal{M}_{\X_0}$ whose
degree $d$ coefficient measures the set of unramified points on
$\X_\eta\times_K K(d)$. Here $K(d)$ denotes the totally ramified
extension $K((\sqrt[d]{\pi}))$ of $K$; the Weil generating series
depends on the choice of $\pi$ if $k$ is not algebraically closed.

If $X$ is a smooth irreducible $k$-variety, endowed with a
dominant morphism $f:X\rightarrow \mathrm{Spec}\,k[\pi]$, then we
denote by $\widehat{X}$ its $\pi$-adic completion. It is a regular
special formal $R$-scheme. We show that, modulo normalization, its
Weil generating series coincides with Denef and Loeser's motivic
zeta function associated to $f$, and its motivic volume coincides
with the motivic nearby cycles \cite{DL3}.

Finally, we study the analytic Milnor fiber $\mathcal{F}_x$ of $f$
at a closed point $x$ of the special fiber $X_s=f^{-1}(0)$. This
object was introduced in \cite{NiSe-Milnor}, and its points and
\'etale cohomology were studied in \cite{NiSe}. In particular, if
$k=\C$, the \'etale cohomology of $\mathcal{F}_x$ coincides with
the singular cohomology of the topological Milnor fiber of $f$ at
$x$, and the Galois action corresponds to the monodromy. We will
show that $\mathcal{F}_x$ completely determines the formal germ of
$f$ at $x$: it determines the completed local ring
$\widehat{\mathcal{O}}_{X,x}$ with its $R$-algebra structure
induced by the morphism $f$. The formal spectrum
$\X_x:=\mathrm{Spf}\,\widehat{\mathcal{O}}_{X,x}$ is a regular
special formal $R$-scheme, and its generic fiber is precisely
$\mathcal{F}_x$. The Weil generating series of $\X_x$ coincides
(modulo normalization) with the local motivic zeta function of $f$
at $x$, and its motivic volume with Denef and Loeser's motivic
Milnor fiber.

To conclude the introduction, we give a survey of the structure of
the paper.
In Section \ref{sec-special}, we
study the basic properties and constructions for special formal
$R$-schemes: generic fibers, formal blow-ups, dilatations, and
resolution of singularities. We also encounter an important
technical complication w.r.t. the finite type-case: if $\X$ is a
special formal $R$-scheme and $\omega$ is a differential form on
$\X_\eta$, there does not necessarily exist an integer $a>0$ such
that $\pi^a\omega$ is defined on $\X$. We call forms which
(locally on $\X$) have this property, $\X$-bounded differential
forms; this notion is important for what follows.

In Section \ref{sec-nearby}, we compute the $\ell$-adic tame
nearby cycles on a regular special formal $R$-scheme $\X$ whose
special fiber is a tame strict normal crossings divisor
(Proposition \ref{nearby}). We prove that any such formal scheme
$\X$ admits an algebraizable \'etale cover (Proposition
\ref{rev}), and then we use Grothendieck's description of the
nearby cycles in the algebraic case \cite[Exp. I]{sga7a}.

We generalize the theory of motivic integration to the class of
special formal $R$-schemes in Section \ref{sec-mot}. In fact,
using dilatations, we construct appropriate models which are
topologically of finite type over $R$, and for which the theory of
motivic integration was developed in \cite{sebag1} and
\cite{motrigid}. Of course, we have to show that the result does
not depend on the chosen model. In particular, we associate a
motivic Serre invariant to any generically smooth special formal
$R$-scheme $\X$ (Definition \ref{def-serre}) and we show that it
can be computed on a N\'eron smoothening (Corollary
\ref{serreneron}).

In Section \ref{sec-comput}, we construct weak N\'eron models for
tame ramifications of regular special formal $R$-schemes whose
special fiber is a strict normal crossings divisor (Theorem
\ref{neronsmooth}) and we obtain a formula for the motivic Serre
invariants of these ramifications. We define the order of a
bounded gauge form $\omega$ on the generic fiber $\X_\eta$ of a
smooth special formal $R$-scheme $\X$ along the connected
components of $\X_0$, and we deduce an expression for the motivic
integral of $\omega$ on $\X$ (Proposition \ref{comput-smooth}).

The trace formula is stated and proven in Section \ref{trace}
(Theorem \ref{trace}). It uses the computation of the motivic
Serre invariants of tame ramifications in Section
\ref{sec-comput}, the computation of the tame nearby cycles in
Section \ref{sec-nearby}, and Laumon's result on equality of
$\ell$-adic Euler characteristics with and without proper support
\cite{laumon-euler}.

In Section \ref{sec-volume}, we consider regular special formal
$R$-schemes $\X$ whose special fibers $\X_s$ are strict normal
crossings divisors, where $char(k)=0$. We define the order of a
bounded gauge form $\omega$ on the generic fiber $\X_\eta$ along
the components of $\X_s$, and we use this notion to compute the
motivic integral of $\omega$ on all the totally ramified
extensions of $\X_\eta$ (Theorem \ref{compvolume-d}). These values
are used to define the volume Poincar\'e series for any
generically smooth special formal $R$-scheme $\X$ and any bounded
gauge form $\omega$ on $\X_\eta$, and we obtain an explicit
expression for this series in terms of a resolution of
singularities (Corollary \ref{serrat}). In particular, the limit
of this series is well-defined, and does not depend on $\omega$,
and we use it to define the motivic volume of $\X$ (Definition
\ref{motvolume}). In Section \ref{sec-gelfand}, we introduce the
Gelfand-Leray form $\omega/d\pi$ associated to a top differential
form $\omega$ on $\X$ over $k$ (Definition \ref{def-gelfand}),
using the fact that the wedge product with $d\pi$ defines an
isomorphism between $\Omega^m_{\X_\eta/K}$ and
$(\Omega^{m+1}_{\X/k})_{rig}$ if $\X$ is a generically smooth
special formal $R$-scheme of pure relative dimension $m$
(Proposition \ref{gelfand}). If $\X$ is regular and $\omega$ is a
gauge form, then $\omega/d\pi$ is a bounded gauge form on
$\X_\eta$ and the volume Poincar\'e series of $(\X,\omega/d\pi)$
depends only on $\X$, and not on $\omega$; we get an explicit
expression in terms of any resolution of singularities
(Proposition \ref{mellin}). We call this series the Weil
generating series associated to $\X$.

Now let $X$ be a smooth, irreducible variety over $k$, endowed
with a dominant morphism $X\rightarrow \mathrm{Spec}\,k[t]$, and
denote by $\widehat{X}$ its $t$-adic completion. We show in
Section \ref{Milnor} that the analytic Milnor fiber of $f$ at a
closed point $x$ of the special fiber $X_s=f^{-1}(0)$ completely
determines the formal germ of $f$ at $x$ (Proposition
\ref{milnor}). Finally, in Section \ref{sectioncompar}, we prove
that the Weil generating series of $\widehat{X}$ and
$\mathcal{F}_x$ coincide, modulo normalization, with the motivic
zeta function of $f$, resp. the local motivic zeta function of $f$
at $x$ (Theorem \ref{comparzeta} and Corollary \ref{comparlocal}),
and we show that the motivic volumes of $\widehat{X}$ and
$\mathcal{F}_x$ correspond to Denef and Loeser's motivic nearby
cycles and motivic Milnor fiber (Theorem \ref{cycles}). This
refines the comparison results in \cite{NiSe}.

I am grateful to Christian Kappen for pointing out a mistake in an
earlier version of this article.
\subsection*{Notation and conventions}
Throughout this article, $R$ denotes a complete discrete valuation
ring, with residue field $k$ and quotient field $K$, and we fix a
uniformizing parameter $\pi$. Some of the constructions require
that $k$ is perfect or that $k$ has characteristic zero; this will
be indicated at the beginning of the section. For any field $F$,
we denote by $F^s$ a separable closure. We denote by $K^{sh}$ the
strict henselization of $K$, by $R^{sh}$ the normalization of $R$
in $K^{sh}$, and by $K^t$ the tame closure of $K$ in $K^s$. We
denote by $\widehat{K^t}$ and $\widehat{K^s}$ the completions of
$K^t$ and $K^s$. We denote by $p$ the characteristic exponent of
$k$, and we fix a prime $\ell$ invertible in $k$. We say that $R'$
is a finite extension of $R$ if $R'$ is the normalization of $R$
in a finite field extension $K'$ of $K$.

For any integer $d>0$ prime to the characteristic exponent of $k$,
we put $K(d):=K[T]/(T^d-\pi)$. This is a totally ramified
extension of degree $d$ of $K$.
 We denote by $R(d)$ the normalization
$R[T]/(T^d-\pi)$ of $R$ in $K(d)$. For any formal $R$-scheme $\X$
and any rigid $K$-variety $X$, we put $\X(d):=\X\times_R R(d)$ and
$X(d):=X\times_K K(d)$. Moreover, we put
$\overline{X}:=X\widehat{\times}_K \widehat{K^t}$.

If $S$ is a scheme, we denote by $S_{red}$ the underlying reduced
scheme. A $S$-variety is a separated reduced $S$-scheme of finite
type.

If $F$ is a field of characteristic exponent $p$, $\ell$ is a
prime invertible in $F$, and $S$ is a variety over $F$, then we
say an \'etale covering $T\rightarrow S$ is \textit{tame} if, for
each connected component $S_i$ of $S$, the degree of the \'etale
covering $T\times_S S_i\rightarrow S_i$ is prime to the
characteristic exponent $p$ of $F$. We call a $\Q^s_\ell$-adic
sheaf $\mathcal{F}$ on $S$ \textit{tamely lisse} if its
restriction to each connected component $S_i$ of $S$ corresponds
to a finite dimensional continuous representation of the
prime-to-$p$ quotient $\pi_1(S_i,s)^{p'}$, where $s$ is a
geometric point of $S_i$. We call a $\Q^s_\ell$-adic sheaf
$\mathcal{F}$ on $S$ \textit{tamely constructible} if there exists
a finite stratification of $S$ into locally closed subsets $S_i$,
such that the restriction of $\mathcal{F}$ to each $S_i$ is tamely
lisse. If $M$ is a torsion ring with torsion orders prime to $p$,
then tamely lisse and tamely constructible sheaves of $M$-modules
on $S$ are defined in the same way.

If $\X$ is a Noetherian adic formal scheme and $Z$ is a closed
subscheme (defined by an open coherent ideal sheaf on $\X$) we
denote by $\widehat{\X/Z}$ the formal completion of $\X$ along
$Z$. If $\mathcal{N}$ is a coherent $\mathcal{O}_{\X}$-module, we
denote by $\widehat{\mathcal{N}/Z}$ the induced coherent
$\mathcal{O}_{\widehat{\X/Z}}$-module. We embed the category of
Noetherian schemes into the category of Noetherian adic formal
schemes by endowing their structure sheaves with the discrete
(i.e. $(0)$-adic) topology. If $\X$ is a Noetherian adic formal
scheme and $\mathcal{J}$ a coherent ideal sheaf on $\X$, we'll
write $V(\mathcal{J})$ for the closed formal subscheme of $\X$
defined by $\mathcal{J}$.

For the theory of $stft$ formal $R$-schemes ($stft$=separated and
topologically of finite type) and the definition of the
Grothendieck ring of varieties $K_0(Var_Z)$ over a separated
scheme $Z$ of finite type over $k$, we refer to \cite{NiSe}. Let
us only recall that $\LL$ denotes the class of the affine line
$\A^1_Z$ in $K_0(Var_Z)$, and that $\mathcal{M}_Z$ denotes the
localized Grothendieck ring $K_0(Var_Z)[\LL^{-1}]$. The
topological Euler characteristic
$$\chi_{top}(X):=\sum_{i\geq 0}(-1)^i\mathrm{dim}\,H^i(X\times_k k^s,\Q_\ell)$$
induces a group morphism $\chi_{top}:\mathcal{M}_Z\rightarrow \Z$.
The definition of the completed localized Grothendieck ring
$\widehat{\mathcal{M}}_Z$ is recalled in \cite[\S 4.1]{NiSe2}.

If $V=\oplus_{i\in\Z}V_i$ is a graded vector space over a field
$F$, such that $V_i=0$ for all but a finite number of $i\in\Z$ and
such that $V_i$ is finite dimensional over $F$ for all $i$, and if
$M$ is a graded endomorphism of $V$, then we define its trace and
zeta function by \begin{eqnarray*}
Tr(M\,|\,V)&:=&\sum_{i\in\Z}(-1)^i Tr(M\,|\,V_i)\in F
\\ \zeta(M\,|\,V;T)&:=&\prod_{i\in\Z}\left(
\mathrm{det}(Id-TM\,|\,V_i)\right)^{(-1)^{i+1}}\in F(T)
\end{eqnarray*}

Finally, if $X$ is a rigid $K$-variety, we put
$H(\overline{X}):=\oplus_{i\geq 0}H^i(\overline{X},\Q_\ell)$ where
the cohomology on the right is Berkovich' $\ell$-adic \'etale
cohomology \cite{Berk-etale}.
\section{Special formal schemes}\label{sec-special}
We recall the following definition: if $A$ is an adic topological
ring with ideal of definition $J$, then the algebra of convergent
power series over $A$ in the variables $(x_1,\ldots,x_n)$ is given
by
$$A\{x_1,\ldots,x_n\}:=\lim_{\stackrel{\longleftarrow}{n\geq
1}}(A/J^n)[x_1,\ldots,x_n]$$
\begin{definition}
Let $\X$ be a Noetherian adic formal scheme, and let $\mathcal{J}$
be its largest ideal of definition. The closed subscheme of $\X$
defined by $\mathcal{J}$ is  denoted by $\X_0$, and is called the
reduction of $\X$. It is a reduced Noetherian scheme.
\end{definition}

This construction defines a functor $(\cdot)_0$ from the category
of Noetherian adic formal schemes to the category of reduced
Noetherian schemes. Note that the natural closed immersion
$\X_0\rightarrow \X$ is a homeomorphism.

\begin{definition}[Special formal schemes \cite{berk-vanish2}, \S 1]
A topological $R$-algebra $A$ is called special, if $A$ is a
Noetherian adic ring and, for some ideal of definition $J$, the
$R$-algebra $A/J$ is finitely generated.

 A special formal $R$-scheme is a separated Noetherian adic formal scheme
$\X$ endowed with a structural morphism $\X\rightarrow
\mathrm{Spf}\,R$, such that $\X$ is a finite union of open formal
subschemes which are formal spectra of special $R$-algebras. In
particular, $\X_0$ is a separated scheme of finite type over $k$.
%

We denote by $\X_s$ the special fiber $\X\times_R k$ of $\X$. It
is a formal scheme over $\mathrm{Spec}\,k$. If $\X$ is $stft$ over
$R$, then $\X_s$ is a separated $k$-scheme of finite type, and
$\X_0=(\X_s)_{red}$.
\end{definition}
Note that our terminology is slightly different from the one in
\cite[\S\,1]{berk-vanish2}, since we impose the additional
quasi-compactness condition.

Special formal schemes are called formal schemes \textit{of
pseudo-finite type} over $R$ in \cite{formal1}. We adopt their
definitions of \'etale, adic \'etale, smooth, and adic smooth
morphisms \cite[2.6]{formal1}. If $\X$ is a special formal
$R$-scheme, we denote by $Sm(\X)$ the open formal subscheme where
the structural morphism $\X\rightarrow \mathrm{Spf}\,R$ is smooth.

Berkovich shows in \cite[1.2]{berk-vanish2} that a topological
$R$-algebra $A$ is special, iff $A$ is topologically
$R$-isomorphic to a quotient of the special $R$-algebra
$$R\{T_1,\ldots,T_m\}[[S_1,\ldots,S_n]]=R[[S_1,\ldots,S_n]]\{T_1,\ldots,T_m\}$$
It follows from \cite{Va1,Va2} that special $R$-algebras are
excellent, as is observed in \cite[p.476]{conrad}.

Any $stft$ formal $R$-scheme is special. Note that a special
formal $R$-scheme is $stft$ over $R$ iff it is $R$-adic. If $\X$
is a special formal $R$-scheme, and $Z$ is a closed subscheme of
$\X_0$, then the formal completion $\widehat{\X/Z}$ of $\X$ along
$Z$ is special.

We say that a special formal $R$-scheme $\X$ is algebraizable, if
$\X$ is isomorphic to the formal completion of a separated
$R$-scheme $X$ of finite type along a closed subscheme $Z$ of its
special fiber $X_s$. In this case, we say that $X/Z$ is an
algebraic model for $\X$. If $\X$ is $stft$ over $R$ and $Z=X_s$
(i.e. $\X$ is isomorphic to the $\pi$-adic completion
$\widehat{X}$ of $X$), we simply say that $X$ is an algebraic
model for $\X$. Finally, if $\widehat{\mathcal{I}}$ is the
$\pi$-adic completion of a coherent $\mathcal{O}_{X}$-module
$\mathcal{I}$, we say that $(X,\mathcal{I})$ is an algebraic model
for $(\X,\widehat{\mathcal{I}})$.

If $\X$ is a special formal $R$-scheme, $R'$ is a finite extension
of $R$, and $\psi$ is a section in $\X(R')$, then we denote by
$\psi(0)$ the image of $\mathrm{Spf}\,R'$ in $\X_0$.
\subsection{The generic fiber of a special formal
scheme}\label{generic} Berthelot explains in \cite[0.2.6]{bert}
how to associate a generic fiber $\X_\eta$ to a special formal
$R$-scheme $\X$ (see also \cite[\S 7]{dj-formal}). This generic
fiber $\X_\eta$ is a separated rigid variety over $K$, not
quasi-compact in general, and is endowed with a canonical morphism
of ringed sites $sp:\X_\eta\rightarrow \X$ (the specialization
map). This construction yields a functor $(.)_\eta$ from the
category of special formal $R$-schemes, to the category of
separated rigid $K$-varieties. We say that $\X$ is generically
smooth if $\X_\eta$ is smooth over $K$.

If $Z$ is a closed subscheme of $\X_0$, then $sp^{-1}(Z)$ is an
open rigid subvariety of $\X_\eta$, canonically isomorphic to
$(\widehat{\X/Z})_\eta$ by \cite[0.2.7]{bert}. We call it the tube
of $Z$ in $\X$, and denote it by $]Z[$.

We recall the construction of $\X_\eta$ in the case where
$\X=\mathrm{Spf}\,A$ is affine, with $A$ a special $R$-algebra,
following \cite[7.1]{dj-formal}. The notation introduced here will
be used throughout the article.

Let $J$ be the largest ideal of definition of $A$. For each
integer $n>0$, we denote by $A[J^n/\pi]$ the subalgebra of
$A\otimes_R K$ generated by $A$ and the elements $j/\pi$ with
$j\in J^n$. We denote by $B_n$ the $J$-adic completion of
$A[J^n/\pi]$ (this is also the $\pi$-adic completion), and we put
$C_n=B_n\otimes_R K$. Then $C_n$ is an affinoid algebra, the
natural map $C_{n+1}\rightarrow C_n$ induces an open embedding of
affinoid spaces $\mathrm{Sp}\,C_n\rightarrow
\mathrm{Sp}\,C_{n+1}$, and by construction,
$\X_\eta=\cup_{n>0}\mathrm{Sp}\,C_n$.

For each $n>0$, there is a natural ring morphism $A\otimes_R
K\rightarrow C_n$ which is flat by \cite[7.1.2]{dj-formal}. These
morphisms induce a natural ring morphism
$$i:A\otimes_R K\rightarrow
\mathcal{O}_{\X_\eta}(\X_\eta)=\cap_{n>0}C_n$$
%

\begin{definition}[\cite{Kiehl},2.3]
A rigid variety $X$ over $K$ is called a quasi-Stein space if
there exists an admissible covering of $X$ by affinoid opens
$X_1\subset X_2\subset \ldots$ such that
$\mathcal{O}_{X}(X_{n+1})\rightarrow \mathcal{O}_{X}(X_{n})$ has
dense image for all $n\geq 1$.
\end{definition}

%
A crucial feature of a
quasi-Stein space $X$ is that $H^i(X,\mathcal{F})$ vanishes for
$i>0$ if $\mathcal{F}$ is a coherent sheaf on $X$, i.e. the global
section functor is exact on coherent modules on $X$. This is
Kiehl's ``Theorem B'' for rigid quasi-Stein spaces
\cite[2.4]{Kiehl}.

\begin{prop}\label{stein}
If $\X=\mathrm{Spf}\,A$ is an affine special formal $R$-scheme,
then $\X_\eta$ is a quasi-Stein space.
\end{prop}
\begin{proof}
Let $J$ be the largest ideal of definition in $A$.
Put
$X_n=\mathrm{Sp}\,C_n$; then
 $$X_1\subset X_2\subset \ldots$$ is an affinoid cover of $\X_\eta$.
Fix an integer $n>0$, and let $\{g_1,\ldots,g_s\}$ be a set of
generators of the ideal $J^n$ in $A$. Then by construction,
$X_n$ consists exactly of the points $x$ in $\X_\eta$ such that
$|(g_j/\pi)(x)|\leq 1$ for $j=1,\ldots,s$, by the isomorphism
\cite[7.1.2]{dj-formal} $$C_{n}\cong
C_{n+1}\{T_1,\ldots,T_s\}/(g_1-\pi T_1,\ldots,g_s-\pi T_s)$$ As
Kiehl observes right after Definition 2.3 in \cite{Kiehl}, this
implies that $\X_\eta$ is quasi-Stein.
  \end{proof}

Let $A$ be a special $R$-algebra, and $\X=\mathrm{Spf}\,A$.
Whenever $M$ is a finite $A$-module, we can define the induced
coherent sheaf $M_{rig}$ on $\X_\eta$ by
$$M_{rig}|_{\mathrm{Sp}\,(C_n)}:=(M\otimes_A C_n)^{\sim}$$
Here $(M\otimes_A C_n)^{\sim}$ denotes the coherent
$\mathcal{O}_{\mathrm{Sp}\,C_n}$-module associated to the finite
$C_n$-module $M\otimes_A C_n$.

If $\X=\mathrm{Spf}\,A$ is topologically of finite type over $R$,
then $\X_\eta$ is simply $\mathrm{Sp}\,(A\otimes_R K)$, and
$M_{rig}$ corresponds to the $(A\otimes_R K)$-module $M\otimes_R
K$.

\begin{lemma}\label{exact}
If $A$ is a special $R$-algebra and $\X=\mathrm{Spf}\,A$, the
functor
$$(.)_{rig}:(Coh_\X)\rightarrow (Coh_{\X_{\eta}}):M\mapsto
M_{rig}$$ from the category $(Coh_\X)$ of coherent
$\mathcal{O}_{\X}$-modules to the category $(Coh_{\X_\eta})$ of
coherent $\mathcal{O}_{\X_\eta}$-modules, is exact.
\end{lemma}
\begin{proof}
This follows from the fact that the natural ring morphism
$A\rightarrow C_n$ is flat for each $n>0$, by
\cite[7.1.2]{dj-formal}.   \end{proof}
\begin{prop}\label{rigmodule}
For any special formal $R$-scheme $\X$, there exists a unique
 functor
$$(.)_{rig}:(Coh_{\X})\rightarrow (Coh_{\X_{\eta}}):M\mapsto
M_{rig}$$ such that
$$M_{rig}|_{\mU_\eta}=(M|_{\mU})_{rig}$$
for any open affine formal subscheme $\mU$ of $\X$.

The functor $(.)_{rig}$ is exact. For any morphism of special
formal $R$-schemes $h:\mY\rightarrow \X$, and any coherent
$\mathcal{O}_{\X}$-module $M$, there is a canonical isomorphism
$(h^*M)_{rig}\cong (h_\eta)^*M_{rig}$. Moreover, if $h$ is a
finite adic morphism, and $N$ is a coherent
$\mathcal{O}_{\mY}$-module, then there is a canonical isomorphism
$(h_*N)_{rig}\cong (h_\eta)_*(N_{rig})$.
\end{prop}
\begin{proof}
Exactness follows immediately from Lemma \ref{exact}. It is clear
that $(.)_{rig}$ commutes with pull-back, so let $h:\mY\rightarrow
\X$ be a finite adic morphism of special formal $R$-schemes, and
let $N$ be a coherent $\mathcal{O}_{\mY}$-module. We may suppose
that $\X=\mathrm{Spf}\,A$ is affine; then $\mY=\mathrm{Spf}\,D$
with $D$ finite and adic over $A$, and $h_*N$ is simply $N$ viewed
as a $A$-module. By \cite[7.2.2]{dj-formal}, the inverse image of
$\mathrm{Sp}\,C_n\subset \X_\eta$ in $\mY_\eta$ is the affinoid
space $\mathrm{Sp}\,(D\otimes_{A}C_n)$, so both
$(h_*N)_{rig}|_{\mathrm{Sp}\,C_n}$ and
$(h_\eta)_*(N_{rig})|_{\mathrm{Sp}\,C_n}$ are associated to the
coherent $C_n$-module $N\otimes_{A}C_n$.   \end{proof}

\begin{example}\label{notadic}
The assumption that $h$ is finite and adic is crucial in the last
part of Proposition \ref{rigmodule}. Consider, for example, the
special formal $R$-scheme $\X=\mathrm{Spf}\,R[[x]]$, and denote by
$h:\X\rightarrow \mathrm{Spf}\,R$ the structural morphism. Choose
a series $(a_n)$ in $K$ such that $|a_n|\to \infty$ as
$n\to\infty$, but with $|a_n|\leq \log\,n$. Then the power series
$f=\sum_{n\geq 0}a_n x^n$ in $K[[x]]$ defines an element of
$\mathcal{O}_{\X_\eta}(\X_\eta)=(f_\eta)_*(\mathcal{O}_{\X})_{rig}$
since it converges on every closed disc $D(0,\rho)$ with $\rho<1$,
but it does not belong to $R[[x]]\otimes_R
K=(h_*\mathcal{O}_{\X})_{rig}$.   \end{example}

\begin{lemma}\label{inject}
If $\X$ is a special formal $R$-scheme and $M$ is a coherent
$\mathcal{O}_{\X}$-module, then the functor $(.)_{rig}$ induces a
natural map of $K$-modules
$$i:M(\X)\otimes_R K\rightarrow M_{rig}(\X_\eta) $$ and this map is
injective. If $M=\mathcal{O}_{\X}$ then $i$ is a map of
$K$-algebras.
\end{lemma}
\begin{proof}
The map $i$ is constructed in the obvious way: if
$\X=\mathrm{Spf}\,A$ is affine and $m$ is an element of $M(\X)$,
then the restriction of $i(m)$ to $\mathrm{Sp}\,C_n$ is simply the
element $m\otimes 1$ of $M_{rig}(\mathrm{Sp}\,C_n)=M\otimes_A
C_n$. The general construction is obtained by gluing.

To prove that $i$ is injective, we may suppose that
$\X=\mathrm{Spf}\,A$ is affine; we'll simply write $M$ instead of
$M(\X)$. Let $m$ be an element of $M$ and suppose that $i(m)=0$.
Suppose that $m$ is non-zero in $M\otimes_R K$, and let
$\mathfrak{M}$ be a maximal ideal in $A\otimes_R K$ such that $m$
is non-zero in the stalk $M_{\mathfrak{M}}$. By
\cite[7.1.9]{dj-formal}, $\mathfrak{M}$ corresponds canonically to
a point $x$ of $\X_\eta$ and there is a natural local homomorphism
$(A\otimes_R K)_{\mathfrak{M}}\rightarrow \mathcal{O}_{\X_\eta,x}$
which induces an isomorphism on the completions, so $i(m)=0$
implies that $m$ vanishes in the $\mathfrak{M}$-adic completion of
$M_{\mathfrak{M}}$. This implies at its turn that $m$ vanishes in
$M_{\mathfrak{M}}$ since $M_{\mathfrak{M}}$ is separated for the
$\mathfrak{M}$-adic topology \cite[7.3.5]{ega1}; this contradicts
our assumption.   \end{proof}
\begin{cor}\label{tor}
If $\X$ is a special formal $R$-scheme and $M$ is a coherent
$\mathcal{O}_{\X}$-module, then $M_{rig}=0$ iff $M$ is annihilated
by a power of $\pi$.
\end{cor}
\begin{proof}
We may assume that $\X$ is affine, say $\X=\mathrm{Spf}\,A$. By
Lemma \ref{inject}, $M(\X)\otimes_R K=0$, so $M$ is annihilated by
a power of $\pi$.   \end{proof}

\begin{lemma}\label{bounded-iso}
Let $\X=\mathrm{Spf}\,A$ be an affine special formal $R$-scheme,
and let $f:M\rightarrow N$ be a morphism of coherent
$\mathcal{O}_{\X}$-modules such that the induced morphism of
coherent $\mathcal{O}_{\X_\eta}$-modules
$f_{rig}:M_{rig}\rightarrow N_{rig}$ is an isomorphism. Then the
natural map
$$f:M(\X)\otimes_R K\rightarrow N(\X)\otimes_R K$$ is an
isomorphism, and fits in the commutative diagram
$$\begin{CD}
M(\X)\otimes_R K@>>> N(\X)\otimes_R K
\\@VVV @VVV
\\M_{rig}(\X_\eta)@>>> N_{rig}(\X_\eta)
\end{CD}$$
where the vertical arrows are injections and the horizontal arrows
are isomorphisms.
\end{lemma}
\begin{proof}
We extend the morphism $f$ to an exact sequence of coherent
$\mathcal{O}_{\X}$-modules
$$\begin{CD}0@>>> \mathrm{ker}(f)@>>> M@>f>> N@>>>
\mathrm{coker}(f)@>>> 0\end{CD}$$ Since $(.)_{rig}$ is an exact
functor by Lemma \ref{exact}, and $(f)_{rig}$ is an isomorphism by
assumption, $\mathrm{ker}(f)_{rig}$ and $\mathrm{coker}(f)_{rig}$
vanish, and hence, $\mathrm{ker}(f)$ and $\mathrm{coker}(f)$ are
$\pi$-torsion modules, by Corollary \ref{tor}. Since $\X$ is
affine, the above exact sequence gives rise to an exact sequence
of $A$-modules
$$\begin{CD}0@>>> \mathrm{ker}(f)(\X)@>>> M(\X)@>f>> N(\X)@>>>
\mathrm{coker}(f)(\X)@>>> 0\end{CD}$$ and by tensoring with $K$,
we obtain the required isomorphism. The remainder of the statement
follows from Lemma \ref{inject}.   \end{proof}

By \cite[7.1.12]{dj-formal}, there is a canonical isomorphism of
$\mathcal{O}_{\X_\eta}$-modules
$$\Omega^i_{\X_\eta/K}\cong \left(\Omega^i_{\X/R}\right)_{rig}$$
for any special formal $R$-scheme $\X$ and each $i\geq 0$.

\subsection{Bounded differential forms}
\begin{definition}\label{mbounded}
%
Let $\X$ be a special formal $R$-scheme. For any $i\geq 0$, we
call an element $\omega$ of $\Omega^i_{\X_\eta/K}(\X_\eta)$ an
$\X$-bounded $i$-form on $\X_\eta$, if there exists a finite cover
of $\X$ by affine open formal subschemes $\{\mU^{(j)}\}_{j\in I}$
such that for each $j\in I$, $\omega|_{\mU^{(j)}_\eta}$ belongs to
the image of the natural map
$$\Omega^i_{\X/R}(\mU^{(j)})\otimes_R K\rightarrow
\Omega^i_{\X_\eta/K}(\mU^{(j)}_\eta)$$
\end{definition}
By Lemma \ref{inject}, this definition is equivalent to saying
that $\omega$ belongs to the image of the natural map
$$(\Omega^i_{\X/R}\otimes_R K)(\X)\rightarrow
\Omega^i_{\X_\eta/K}(\X_\eta)$$ where $\Omega^i_{\X/R}\otimes_R K$
is a tensor product of sheaves on $\X$.

%


If $\X$ is $stft$ over $R$ then any differential form on
$\X_{\eta}$ is $\X$-bounded, by quasi-compactness of $\X_\eta$.
This is false in general: see Example \ref{notadic} for an example
of an unbounded $0$-form.

\begin{lemma}\label{bounded-affine}
If $\X=\mathrm{Spf}\,A$ is an affine special formal $R$-scheme,
and $i\geq 0$ is an integer, then an element $\omega$ of
$\Omega^i_{\X_\eta/K}(\X_\eta)$ is $\X$-bounded iff it belongs to
the image of the natural map
$$\Omega^i_{\X/R}(\X)\otimes_R K\rightarrow \Omega^i_{\X_\eta/K}(\X_\eta)$$
\end{lemma}
\begin{proof}
Since $\X$ is affine, $$ (\Omega^i_{\X/R}\otimes_R
K)(\X)=\Omega^i_{\X/R}(\X)\otimes_R K$$
  \end{proof}
\begin{cor}
Let $\X$ be a special formal $R$-scheme, and let $i\geq 0$ be an
integer. If $\omega$ is a $\X$-bounded $i$-form on $\X_\eta$, then
for any finite cover of $\X$ by affine open formal subschemes
$\{\mU^{(j)}\}_{j\in I}$,  and for each $j\in I$,
$\omega|_{\mU^{(j)}_\eta}$ belongs to the image of the natural map
$$\Omega^i_{\X/R}(\mU^{(j)})\otimes_R K\rightarrow
\Omega^i_{\X_\eta/K}(\mU^{(j)}_\eta)$$
\end{cor}
\begin{lemma}\label{bounded-function}
Let $\X$ be a special formal $R$-scheme, such that $\X_\eta$ is
reduced. An element $f$ of $\mathcal{O}_{\X_\eta}(\X_\eta)$ is
$\X$-bounded iff it is bounded, i.e. iff there exists an integer
$M$ such that $|f(x)|\leq M$ for each point $x$ of $\X_\eta$.
\end{lemma}
\begin{proof}
Since an element $f$ of $\mathcal{O}_{\X}(\X)$ satisfies
$|f(x)|\leq 1$ for each point $x$ of $\X_\eta$ by
\cite[7.1.8.2]{dj-formal}, it is clear that an $\X$-bounded
analytic function on $\X_\eta$ is bounded.

Assume, conversely, that $f$ is a bounded analytic function on
$\X_\eta$. To show that $f$ is $\X$-bounded, we may suppose that
$\X=\mathrm{Spf}\,A$ is affine and flat. Since the natural map
$A\otimes_R K\rightarrow \mathcal{O}_{\X_\eta}(\X_\eta)$ is
injective, $A\otimes_R K$ is reduced; since $A$ is $R$-flat, $A$
is reduced. If $A$ is integrally closed in $A\otimes_R K$, then
the image of the natural map $A\otimes_R K\rightarrow
\mathcal{O}_{\X_\eta}(\X_\eta)$ coincides with the set of bounded
functions on $\X_\eta$, by \cite[7.4.1-2]{dj-formal}. So it
suffices to note that the natural map $A\otimes_R K\rightarrow
B\otimes_R K$ is bijective, where $B$ is the normalization of $A$
in $A\otimes_R K$ ($B$ is a special $R$-algebra since it is finite
over $A$, by excellence of $A$; see \cite{conrad}).
%
  \end{proof}
\subsection{Admissible blow-ups and dilatations}
Let $\X$ be a Noetherian adic formal scheme, let $\mathcal{J}$ be
an ideal of definition, and let $\mathcal{I}$ be any coherent
ideal sheaf on $\X$. Following the $tft$-case in \cite[\S
2]{formrigI}, we state the following definition.

\begin{definition}[Formal blow-up]
The formal blow-up of $\X$ with center $\mathcal{I}$ is the
morphism of formal schemes
$$\X':=\lim_{\stackrel{\longrightarrow}{n\geq
1}}\mathrm{Proj}\left(\oplus_{d=0}^{\infty}\mathcal{I}^d\otimes_{\mathcal{O}_{\X}}(\mathcal{O}_{\X}/\mathcal{J}^n)\right)
\rightarrow \X$$
\end{definition}

\begin{prop}\label{blow-up1}
Let $\X$ be a Noetherian adic formal scheme with ideal of
definition $\mathcal{J}$, let $\mathcal{I}$ be a coherent ideal
sheaf on $\X$, and consider the formal blow-up $h:\X'\rightarrow
\X$ of $\X$ at $\mathcal{I}$.
\begin{enumerate}
\item If $\mathfrak{U}=\mathrm{Spf}\,A$ is an affine open formal
subscheme of $\X$, 
then the restriction of $h$ over $\mathfrak{U}$ is the
$\mathcal{J}(\mU)$-adic completion of the scheme-theoretic blow-up
of $\mathrm{Spec}\,A$ at the ideal $\mathcal{I}(\mU)$ of $A$.
\item The blow-up morphism $\X'\rightarrow \X$ is adic and
topologically of finite type. In particular, $\X'$ is a Noetherian
adic formal scheme. \item (Universal property) The ideal
$\mathcal{I}\mathcal{O}_{\X'}$ is invertible on $\X'$, and each
morphism of adic formal schemes $g:\mathfrak{Y}\rightarrow
\mathfrak{X}$ such that $\mathcal{I}\mathcal{O}_{\mY}$ is
invertible, factors uniquely through a morphism of formal schemes
$\mathfrak{Y}\rightarrow \X'$. \item The formal blow-up commutes
with flat base change: if $f:\mathfrak{Y}\rightarrow \X$ is a flat
morphism of Noetherian adic formal schemes, then
$$\X'\times_{\X}\mathfrak{Y}\rightarrow \mathfrak{Y}$$ is the formal blow-up of $\mathfrak{Y}$ at $\mathcal{I}\mathcal{O}_{\mY}$.
 \item If $\mathcal{K}$ is an open coherent ideal sheaf on $\X$, defining a closed subscheme $Z$ of $\X$,
 then the formal blow-up of $\widehat{\X/Z}$ at
$\widehat{\mathcal{I}/Z}$ is the formal completion along $Z$ of
the formal blow-up of $\X$ at $\mathcal{I}$.
\end{enumerate}
\end{prop}
\begin{proof}
Point (1) follows immediately from the definition, and (2) follows
from (1).
 In (3) and (4) we may assume that $\X$ and $\mY$ are affine; then the result follows from the
 corresponding properties for schemes, using (1). Point (5) is a special case of (4).
  \end{proof}
\begin{corollary}\label{blow-up2}
Let $\X$ be a special formal $R$-scheme with ideal of definition
$\mathcal{J}$, let $\mathcal{I}$ be a coherent ideal sheaf of
$\X$, and consider the admissible blow-up $h:\X'\rightarrow \X$ of
$\X$ at $\mathcal{I}$.
\begin{enumerate}
\item The blow-up $\X'$ is a special formal $R$-scheme. \item If
$\X$ is flat over $R$, then $\X'$ is flat over $R$.
\end{enumerate}
\end{corollary}
\begin{proof}
Point (1) follows immediately from Proposition \ref{blow-up1}(2).
 To prove (2), we may assume that
$\X=\mathrm{Spf}\,A$ is affine; then flatness of $\X'$ follows
from Proposition \ref{blow-up1}(1), the fact that the
scheme-theoretic blow-up of $\mathrm{Spec}\,A$ at
$\mathcal{I}(\X)$ is flat over $R$, and flatness of the completion
morphism.
  \end{proof}

Let $\X$ be a special formal $R$-scheme, and let $\mathcal{J}$ be
an ideal of definition. Let $\mathcal{I}$ be a coherent ideal
sheaf on $\X$, open w.r.t. the $\pi$-adic topology (i.e.
$\mathcal{I}$ contains a power of $\pi$). We will call such an
ideal sheaf $\pi$-\textit{open}. We do not demand that
$\mathcal{I}$ is open w.r.t. the $\mathcal{J}$-adic topology on
$\X$. If $\mathcal{I}$ is $\pi$-open, we call the blow-up
$\X'\rightarrow \X$ with center $\mathcal{I}$
\textit{admissible}\footnote{Contrary to the terminology used in
\cite{formrigI} for the $stft$ case, our definition of admissible
blow-up does not assume any flatness conditions on $\X$.}.

We can give an explicit description of admissible blow-ups in the
affine case, generalizing \cite[2.2]{formrigI}.
\begin{lemma}\label{blow-explicit}
Let $A$ be a special $R$-algebra, with ideal of definition $J$,
and let $I=(f_1,\ldots,f_q)$ be a $\pi$-open ideal in $A$. Put
$\X=\mathrm{Spf}\,A$. Let $h:\X'\rightarrow \X$ be the admissible
blow-up of $\X$ at $I$. The scheme-theoretic blow-up of
$\mathrm{Spec}\,A$ at $I$ is covered by open charts
$\mathrm{Spec}\,A_i$, $i=1,\ldots,p$, where
\begin{eqnarray*}
A_i'&=&A[\frac{\xi_1}{\xi_i},\ldots,\frac{\xi_p}{\xi_i}]/(f_i\frac{\xi_j}{\xi_i}-f_j)_{j=1,\ldots,p}
\\A_i&=&A'_i/(f_i-\mathrm{torsion})
\end{eqnarray*}
(here the $\xi_i/\xi_j$ serve as variables, except when $i=j$). We
write $\widehat{A_i}$ and $\widehat{A_i'}$ for the $J$-adic
completions of $A_i$, resp. $A_i'$. Then
\begin{eqnarray*}
\widehat{A_i'}&=&A\{\frac{\xi_1}{\xi_i},\ldots,\frac{\xi_p}{\xi_i}\}/(f_i\frac{\xi_j}{\xi_i}-f_j)_{j=1,\ldots,p}
\\\widehat{A_i}&=&\widehat{A_i'}/(f_i-\mathrm{torsion})
\end{eqnarray*}
and $\mathrm{Spf}\,\widehat{A_i'}$ is the open formal subscheme of
$\X'$ where $f_i$ generates $I\mathcal{O}_{\X'}$. In particular,
$\{\mathrm{Spf}\,\widehat{A_i'}\}_{i=1,\ldots,p}$ is an open cover
of $\X'$.

If, moreover, $A$ is flat over $R$, then
$$\widehat{A_i}=\widehat{A_i'}/(f_i-\mathrm{torsion})=\widehat{A_i'}/(\pi-\mathrm{torsion})$$
\end{lemma}
\begin{proof}
The proof is similar to the $stft$-case \cite[2.2]{formrigI}.
First, we show that
$\widehat{A_i}=\widehat{A_i'}/(f_i-\mathrm{torsion})$. Since $A_i$
is a finite $A'_i$-module,
$\widehat{A_i}=A_i\otimes_{A'_i}\widehat{A'_i}$. Since $A'_i$ is
Noetherian, $\widehat{A'_i}$ is flat over $A'_i$, so
$$A'_i/(f_i-\mathrm{torsion})\otimes_{A'_i}\widehat{A'_i}=\widehat{A'_i}/(f_i-\mathrm{torsion})$$

Now, we show that $\pi$-torsion and $f_i$-torsion coincide in
$\widehat{A'_i}$ if $A$ is $R$-flat. Since $I$ is $\pi$-open in
$A$, it contains a power of $\pi$. Since $f_i$ generates
$I\widehat{A'_i}$, the $f_i$-torsion is contained in the
$\pi$-torsion. But $A_i$ is $R$-flat since $A$ is $R$-flat, so
$\widehat{A_i}=\widehat{A_i'}/(f_i-\mathrm{torsion})$ is $R$-flat,
i.e. has no $\pi$-torsion.

The remainder of the statement is clear.   \end{proof}
\begin{prop}\label{iso}
Let $\X$ be a special formal $R$-scheme, and let $h:\X'\rightarrow
\X$ be an admissible blow-up with center $\mathcal{I}$. The
induced morphism of rigid varieties $h_\eta:\X'_\eta\rightarrow
\X_\eta$ is an isomorphism.
\end{prop}
\begin{proof}
 We may assume that $\X$ is affine, say
$\X=\mathrm{Spf}\,A$, with $J$ as largest ideal of definition. Let
$I=(f_1,\ldots,f_p)$ be a $\pi$-open ideal in $A$, and let
$h:\X'\rightarrow \X$ be the blow-up of $\X$ at $I$. We define
$B_n$ and $C_n$ as in Section \ref{generic}, for $n>0$.

We adopt the notation from Lemma \ref{blow-explicit}. Since the
admissible blow-up $h$ is adic, and the induced morphism
$V(J\mathcal{O}_{\X'})\rightarrow V(J)$ is of finite type, it
follows from \cite[7.2.2]{dj-formal} that the restriction of
$h_\eta:(\mathrm{Spf}\,\widehat{A_i})_\eta\rightarrow \X_\eta$
over $\mathrm{Sp}\,C_n$ is given by  $$h_{n,i}:
\mathrm{Sp}\,C_n\widehat{\otimes}_{A}\widehat{A_i}\rightarrow\mathrm{Sp}\,C_n$$
for each $i$ and each $n$.  However, the natural map
$$C_n\widehat{\otimes}_A\widehat{A'_i}\rightarrow C_n\widehat{\otimes}_{A}\widehat{A_i}$$
is an isomorphism since the $f_i$-torsion in $\widehat{A'_i}$ is
killed if we invert $\pi$ (because $f_i$ divides $\pi$ in $A'_i$,
as $f_i$ generates the open ideal $IA'_i$). Hence,
$\{h_{n,1},\ldots,h_{1,p}\}$ is nothing but the rational cover of
$\mathrm{Sp}\,C_n$ associated to the tuple $(f_1,\ldots,f_p)$ (see
\cite[8.2.2]{BGR}); note that these elements generate the unit
ideal in $C_n$, since $I$ contains a power of $\pi$, which is a
unit in $C_n$.   \end{proof}
\begin{definition}[Dilatation]\label{dilat}
Let $\X$ be a flat special formal $R$-scheme, and let
$\mathcal{I}$ be a coherent ideal sheaf on $\X$ containing $\pi$.
Consider the admissible blow-up $h:\X'\rightarrow \X$ with center
$\mathcal{I}$. If $\mathfrak{U}$ is the open formal subscheme of
$\X'$ where $\mathcal{I}\mathcal{O}_{\X'}$ is generated by $\pi$,
we call $\mathfrak{U}\rightarrow \X$ the dilatation of $\X$ with
center $\mathcal{I}$.
\end{definition}
\begin{prop}\label{dilatprop}
Let $\X$ be a flat special formal $R$-scheme, let $\mathcal{I}$ be
a coherent ideal sheaf on $\X$ containing $\pi$, and let $Z$ be a
closed subscheme of $\X_0$. The dilatation of $\widehat{\X/Z}$
with center $\widehat{\mathcal{I}/Z}$ is the formal completion
along $Z$ of the dilatation of $\X$ with center $\mathcal{I}$.

If, moreover, $\X$ is $stft$ over $R$ and $(\X,\mathcal{I})$ has
an algebraic model $(X,I)$, then the dilatation of $\widehat{X/Z}$
at $\widehat{\mathcal{I}/Z}$ is the formal completion along $Z$ of
the dilatation of $X$ at $I$ (as defined in \cite[3.2/1]{neron}).
\end{prop}
\begin{proof}
This is clear from the definition.   \end{proof}
\begin{prop}\label{univ-dilat}
Let $\X$ be a flat special formal $R$-scheme, and let
$\mathcal{I}$ be a coherent ideal sheaf on $\X$ containing $\pi$.
Let $h:\mathfrak{U}\rightarrow \X$ be the dilatation with center
$\mathcal{I}$, and denote by $\mZ$ the closed formal subscheme of
$\X_s$ defined by $\mathcal{I}$. The dilatation $\mathfrak{U}$ is
a flat special formal $R$-scheme, and
$h_s:\mathfrak{U_s}\rightarrow \X_s$ factors through $\mZ$. The
induced morphism $h_\eta:\mathfrak{U}_\eta\rightarrow \X_\eta$ is
an open immersion.

If $g:\mathfrak{V}\rightarrow \X$ is any morphism of flat special
formal $R$-schemes such that $g_s:\mathfrak{V}_s\rightarrow \X_s$
factors through $\mZ$, then there is a unique morphism of formal
$R$-schemes $i:\mathfrak{V}\rightarrow \mathfrak{U}$ such that
$g=h\circ i$.

If $\mathcal{I}$ is open, then $\mathfrak{U}$ is $stft$ over $R$.
\end{prop}
\begin{proof}
It is clear that $h_s$ factors through $\mZ$. The morphism
$\mathfrak{U}_\eta\rightarrow \X_\eta$ is an open embedding
because $\mathfrak{U}$ is an open formal subscheme of the blow-up
$\X'$ of $\X$ at $\mathcal{I}$, and $\X'_\eta\rightarrow \X_\eta$
is an isomorphism by Proposition \ref{iso}.

Since $g_s$ factors through $\mZ$, we have
$\mathcal{I}\mathcal{O}_{\mV}=(\pi)$.
In particular, by flatness of $\mV$, the ideal
$\mathcal{I}\mathcal{O}_{\mV}$ is invertible, and by the universal
property of the blow-up, $g$ factors uniquely through a morphism
$i:\mV\rightarrow \X'$ to the blow-up $\X'\rightarrow \X$ at
$\mathcal{I}$. The image of $\mathfrak{V}$ in $\X'$ is necessarily
contained in $\mathfrak{U}$ since $\pi$ generates
$\mathcal{I}\mathcal{O}_{\mV}$: if $v$ were a closed point of
$\mathfrak{V}$ mapping to a point $x$ in $\X'\setminus
\mathfrak{U}$, then we could write $\pi=a\cdot f$ in
$\mathcal{O}_{\X',x}$ with $f\in \mathcal{I}\mathcal{O}_{\X'}$ and
$a$ not a unit. Thus yields $\pi=i^*a\cdot i^*f$ in
$\mathcal{O}_{\mathfrak{V},v}$, but since $i^*f$ belongs to
$\mathcal{I}\mathcal{O}_{\mV}$, we have also $i^*f=b\cdot \pi$ in
$\mathcal{O}_{\mathfrak{V},v}$, so $\pi=c\cdot \pi$ with $c$ not a
unit, and $0=(1-c)\cdot \pi$. Since $1-c$ is invertible in
$\mathcal{O}_{\mathfrak{V},v}$, this would mean that $\pi=0$ in
$\mathcal{O}_{\mathfrak{V},v}$, which contradicts flatness of
$\mathfrak{V}$ over $R$.

Finally, assume that $\mathcal{I}$ is open. Then the ideal
$\mathcal{I}\mathcal{O}_\mU$ is open, and by definition of the
dilatation, it is generated by $\pi$. This implies that
$\mathcal{I}\mathcal{O}_{\mU}$ is an ideal of definition, and that
$\mU$ is $stft$.   \end{proof}
\begin{prop}\label{dilat-comm}
Let $\X$ be a flat special formal $R$-scheme, let $U$ be a reduced
closed subscheme of $\X_0$, and denote by $\mU\rightarrow \X$ the
completion map of $\X$ along $U$. If we denote by $\mU'\rightarrow
\mU$ and $\X'\rightarrow \X$ the dilatations with center
$\mU_0=U$, resp. $\X_0$, then there exists a unique morphism of
formal $R$-schemes $\mU'\rightarrow \X'$ such that the square
$$\begin{CD}
\mU'@>>>\mU
\\@VVV @VVV
\\\X'@>>> \X
\end{CD}$$
commutes, and this morphism $\mU'\rightarrow \X'$ is the
dilatation of $\X'$ with center $\X'_s\times_{\X_0}U$.
\end{prop}
\begin{proof}
The existence of such a unique morphism $\mU'\rightarrow \X'$
follows immediately from Proposition \ref{univ-dilat}; to show
that this is the dilatation with center $\X'_s\times_{\X_0}U$, it
suffices to check that $\mU'\rightarrow \X'$ satisfies the
universal property in Proposition \ref{univ-dilat}. Let
$h:\mZ\rightarrow \X'$ be a morphism of flat special formal
$R$-schemes such that $h_s:\mZ_s\rightarrow \X'_s$ factors through
$\X'_s\times_{\X_0}U$. This means that the composed morphism
$\mZ_s\rightarrow \X_0$ factors through $U$, and hence,
$\mZ\rightarrow \X$ factors through a morphism $g:\mZ\rightarrow
\mU$. Moreover, again by Proposition \ref{univ-dilat}, $g$ factors
uniquely through a morphism $f:\mZ\rightarrow \mU'$.

On the other hand, let $f':\mZ\rightarrow \mU'$ be another
morphism of formal $R$-schemes such that $h$ is the composition of
$f'$ with $\mU'\rightarrow \X'$. Then the compositions of $f$ and
$f'$ with $\mU'\rightarrow \mU$ coincide, so $f=f'$ by the
uniqueness property in Proposition \ref{univ-dilat} for the
dilatation $\mU'\rightarrow \mU$.   \end{proof}

Berthelot's construction of the generic fiber of a special formal
$R$-scheme can be restated in terms of dilatations. Let
$\mathcal{J}$ be an ideal of definition of $\X$, and for any
integer $e> 0$, consider the dilatation
$$h^{(e)}:\mathfrak{U}^{(e)}\rightarrow \X$$ with center
$(\pi,\mathcal{J}^e)$. The formal $R$-scheme $\mU^{(e)}$ is $stft$
over $R$, and $h^{(e)}_\eta$ is an open immersion. By the
universal property of the dilatation, $h^{(e)}$ can be decomposed
uniquely as
$$\begin{CD} \mathfrak{U}^{(e)}@>h^{(e,e')}>> \mathfrak{U}^{(e')}@>h^{(e')}>> \X\end{CD}$$
for any pair of integers $e'\geq e\geq 0$. Moreover, $h^{(e,e')}$
induces an open immersion
$h^{(e,e')}_\eta:\mU^{(e)}_\eta\hookrightarrow \mU^{(e')}_\eta$.

\begin{lemma}\label{imag}
The image of the open immersion
$$h^{(e)}_\eta:\mU^{(e)}_\eta\rightarrow \X_\eta$$ consists of the
points $x\in \X_\eta$ such that $|f(x)|\leq |\pi|$ for each
element $f$ of the stalk $(\mathcal{J}^e)_{sp(x)}$.
\end{lemma}
\begin{proof}
Let $R'$ be a finite extension of $R$, and let $\psi$ be a section
in $\X(R')$. Then by the universal property in Proposition
\ref{univ-dilat}, $\psi$ lifts to a section in $\mU^{(e)}(R')$ iff
$(\mathcal{J}^e,\pi)\mathcal{O}_{\Spf R'}$ (the pull-back through
$\psi$) is generated by $\pi$. If we denote by $x$ the image of
the morphism $\psi_\eta$ in $\X_\eta$, this is equivalent to
saying that $|f(x)|\leq |\pi|$ for all $f$ in
$(\mathcal{J}^e)_{sp(x)}$.   \end{proof}

\begin{prop}
The set $\{\mU^{(e)}_\eta\,|\,e>0\}$ is an admissible cover of
$\X_\eta$.
\end{prop}
\begin{proof}
Let $Y=\mathrm{Sp}\, B$ be any affinoid variety over $K$, endowed
with a morphism of rigid $K$-varieties $\varphi:Y\rightarrow
\X_\eta$. We have to show that the image of $\varphi$ is contained
in $\mU^{(e)}_\eta$, for $e$ sufficiently large. We may assume
that $\X$ is affine, say $\X=\Spf A$. Since $\mathcal{J}$ is an
ideal of definition on $\X$, we have $|f(x)|<1$ for any $f\in
\mathcal{J}(\X)$ and any $x\in \X_\eta$. Since $\mathcal{J}$ is
finitely generated, and by the Maximum Modulus Principle
\cite[6.2.1.4]{BGR}, there exists a value $e>0$ such that for each
element $f$ in $\mathcal{J}(\X)^e$ and each point $y$ of $Y$,
$|f(\varphi(y))|<|\pi|$. By Lemma \ref{imag}, this implies that
the image of $\varphi$ is contained in $\mU_\eta^{(e)}$.
\end{proof} Hence, we could have \textit{defined} $\X_\eta$ as the
limit of the direct system $(\mU^{(e)}_\eta,h^{(e,e')}_\eta)$ in
the category of rigid $K$-varieties.

\begin{remark}
If we assume that $\X=\mathrm{Spf}\,A$ is affine and that
$\mathcal{J}$ is the largest ideal of definition on $\X$, then the
dilatation $\mathfrak{U}^{(e)}\rightarrow \X$ is precisely the
morphism $\mathrm{Spf}\,B_n\rightarrow \X$ introduced at the
beginning of Section \ref{generic}. This can be seen by using the
explicit description in Lemma \ref{blow-explicit}.
\end{remark}

\subsection{Irreducible components of special formal schemes}
If $\X$ is a special formal $R$-scheme, the underlying topological
space $|\X|=|\X_0|$, and even the scheme $\X_0$, reflect rather
poorly the geometric properties of $\X$. For instance, if $A$ is a
special $R$-algebra, $|\mathrm{Spf}\,A|$ can be irreducible even
when $\mathrm{Spec}\,A$ is reducible (e.g. for $A=R[[x,y]]/(xy)$,
where $(\mathrm{Spf}\,A)_0=\mathrm{Spec}\,k$). Conversely,
$|\mathrm{Spf}\,A|$ can be reducible even when $A$ is integral
(e.g. for $A=R\{x,y\}/(\pi-xy)$).

Therefore, a more subtle definition of the irreducible components
of $\X$ is needed. We will use the normalization map
$\widetilde{\X}\rightarrow \X$ constructed in \cite{conrad} (there
normalization was already used to define the irreducible
components of a rigid variety). The normalization map is a finite
morphism of special formal $R$-schemes.

\begin{definition}
Let $\X$ be a special formal $R$-scheme, and let
$\widetilde{\X}\rightarrow \X$ be a normalization map. We say that
$\X$ is irreducible if $|\widetilde{\X}|$ is connected.
\end{definition}

\begin{definition}\label{irred}
Let $\X$ be a special formal $R$-scheme, and let
$h:\widetilde{\X}\rightarrow \X$ be a normalization map. We denote
by $\widetilde{\X}_i$, $i=1,\ldots,r$ the connected components of
$\widetilde{\X}$ (defined topologically), and by $h_i$ the
restriction of $h$ to $\widetilde{\X}_i$.

For each $i$, we denote by $\X_i$ the reduced closed subscheme of
$\X$ defined by the kernel of the natural map
$\psi_i:\mathcal{O}_{\X}\rightarrow
(h_i)_*\mathcal{O}_{\widetilde{\X}_i}$.
We call $\X_i$, $i=1,\ldots,r$ the irreducible components of $\X$.
\end{definition}

Note that $(h_i)_*\mathcal{O}_{\widetilde{\X}_i}$ is coherent, by
finiteness of $h$, so the kernel of $\psi_i$ is a coherent ideal
sheaf on $\X$ and $\X_i$ is well-defined.

In the affine case, the irreducible components of $\X$ correspond
to the minimal prime ideals of the ring of global sections, as one
would expect:

\begin{lemma}\label{irred-affine}
Let $\X=\Spf A$ be an affine special formal $R$-scheme, and denote
by $P_i$, $i=1,\ldots,r$ the minimal prime ideals of $A$. Then the
irreducible components of $\X$ are given by $\Spf A/P_i$,
$i=1,\ldots,r$.
\end{lemma}
\begin{proof}
If $A\rightarrow \widetilde{A}$ is a normalization map, then
$\widetilde{A}=\prod_{i=1}^{r}\widetilde{A/P_i}$.  Hence,
$\widetilde{\X}_i=\Spf (\widetilde{A/P_i})$ are the connected
components of $\widetilde{\X}$, and $\X_i=\Spf A/P_i$ for each
$i$.   \end{proof}

\begin{lemma}
With the notations of Definition \ref{irred}, the morphism
$h_i:\widetilde{\X}_i\rightarrow \X_i$ induced by $h$ is a
normalization map for each $i$. Hence, $\X_i$ is irreducible for
each $i$.
\end{lemma}
\begin{proof}
Fix an index $i$ in $\{1,\ldots,r\}$, let $\mU=\Spf A$ be an open
affine formal subscheme of $\X$, and denote by
$P_j,\,j=1,\ldots,q$ the minimal prime ideals in $A$. Since
normalization commutes with open immersions \cite[1.2.3]{conrad},
$\widetilde{\X}_i\cap h^{-1}(\mU)$ is a union of connected
components of $\widetilde{\mU}$, i.e. it is of the form
$\mathrm{Spf}\,\prod_{j\in J}\widetilde{A/P_j}$ for some subset
$J$ of $\{1,\ldots,q\}$. Then by definition, $\X_i\cap \mU$ is the
closed formal subscheme defined by the ideal $P=\cap_{j\in J}P_j$,
and the minimal prime ideals of $A/P$ are precisely the images of
the ideals $P_j,\,j\in J$. This means that $h^{-1}(\mU)\cap
\widetilde{\X}_i\rightarrow \mU\cap \X_i$ is a normalization map.
  \end{proof}
%

There is an important pathology in comparison to the scheme case:
a non-empty open formal subscheme of an irreducible special formal
$R$-scheme is not necessarily irreducible, as is shown by the
following example. Put $A=R\{x,y\}/(xy-\pi)$, and $\X=\Spf A$.
Then $\X$ is irreducible since $A$ is a domain. However, if we
denote by $O$ the point of $\X$ defined by the open ideal
$(\pi,x,y)$, then $\X\setminus \{O\}$ is disconnected.

\begin{prop}\label{norm-comm}
Let $Y$ be either an excellent Noetherian scheme over $R$ or a
special formal $R$-scheme. In the first case, let $\mathcal{I}$ be
a coherent ideal sheaf on $Y$ such that the $\mathcal{I}$-adic
completion $\mY$ of $Y$ is a special formal $R$-scheme; in the
second, let $\mathcal{I}$ be any open coherent ideal sheaf on $Y$,
and denote again by $\mY$ the $\mathcal{I}$-adic completion of
$Y$.

If $h:\widetilde{Y}\rightarrow Y$ is a normalization map, then its
$\mathcal{I}$-adic completion
$\widehat{h}:\widetilde{\mY}\rightarrow \mY$ is also a
normalization map.
\end{prop}
\begin{proof}
We may assume that $Y$ is affine, say $Y=\Spec A$ (resp. $\Spf
A$), so that $\mathcal{I}$ is defined by an ideal $I$ of $A$.
Moreover, we may assume that $A$ is reduced. Denote by
$\widehat{A}$ the $I$-adic completion of $A$. Since
$\widetilde{A}$ is finite over $A$, $B:=\widetilde{A}\otimes_A
\widehat{A}$ is the $I$-adic completion of $\widetilde{A}$. Hence,
by \cite[1.2.2]{conrad}, it suffices to show that $B$ is normal,
and that $B/P$ is reduced for all minimal prime ideals $P$ of $A$.

Normality of $B$ follows from excellence of $\widetilde{A}$. Let
$P$ be any minimal prime ideal of $A$. Then $\widetilde{A}/P$ is
reduced, and hence, so is $B/P$, again by excellence of
$\widetilde{A}$ (see \cite[7.8.3(v)]{ega4.2}).   \end{proof}
\begin{cor}\label{irred-compl}
We keep the notation of Lemma \ref{norm-comm}, and we denote by
$Z_1,\ldots,Z_q$ the connected components of the closed subscheme
$Z$ of $Y$ defined by $\mathcal{I}$. If $Y_1,\ldots,Y_r$ are the
irreducible components of $Y$, then the irreducible components of
$\mY$ are given by the irreducible components of
$\widehat{Y_i/Z_j}$ for $i=1,\ldots,r$ and $j=1,\ldots,q$ (where
$\widehat{Y_i/Z_j}$ may be empty for some $i,j$).
\end{cor}
\begin{proof}
Denote by $h_i:\widetilde{Y}_i\rightarrow Y$ the restriction of
$h$ to $\widetilde{Y}_i$, for each $i$, and by
$\widehat{h_i}:\mU_i\rightarrow \widehat{Y/Z}$ its
$\mathcal{I}$-adic completion. By exactness of the completion
functor, the kernel of
$$\mathcal{O}_{\widehat{Y/Z}}\rightarrow
(\widehat{h_i})_*\mathcal{O}_{\mU_i}$$ is the defining ideal sheaf
of the completion $\widehat{Y_i/Z}$, for each $i$.
\end{proof}
\subsection{Strict normal crossings and resolution of
singularities}\label{strict} \begin{definition}\label{sncd} A
special formal $R$-scheme $\X$ is regular if, for any point $x$ of
$\X$, the local ring $\mathcal{O}_{\X,x}$ is regular (it suffices
to check this at closed points).

Let $\X$ be a regular special formal $R$-scheme. We say that a
coherent ideal sheaf $\mathcal{I}$ on $\X$ is a strict normal
crossings ideal, if the following conditions hold:

(1) there exists at each point $x$ of $\X$ a regular system of
local parameters $$(x_0,\ldots,x_m)$$ in $\mathcal{O}_{\X,x}$ such
that $\mathcal{I}_x$ is generated by $\prod_{j=0}^{m}(x_j)^{M_j}$
for some tuple $M$ in $\N^{\{0,\ldots,m\}}$.

(2) if $\mE$ is the closed formal subscheme of $\X$ defined by
$\mathcal{I}$, then the irreducible components of $\mE$ are
regular.
\end{definition}

\begin{lemma}\label{def-equiv}
If (1) holds, then condition (2) is equivalent to the condition
that $\mathcal{O}_{\mE_i,x}$ is a domain for each irreducible
component $\mE_i$ of $\mE$ and each point $x$ of $\mE_i$.
\end{lemma}
\begin{proof}
 If $\mE_i$ is regular, then $\mathcal{O}_{\mE_i,x}$ is regular
and hence a domain. Conversely, assume that
$\mathcal{O}_{\mE_i,x}$ is a domain. Using the notations in
condition (1) of Definition \ref{sncd}, we choose an open affine
neighbourhood $\mU=\Spf A$ of $x$ in $\X$ such that
$x_0,\ldots,x_m$ are defined on $\mU$ and such that $\mathcal{I}$
is generated by $\prod_{j=0}^{m}(x_j)^{M_j}$ on $\mU$. Denote by
$\mathfrak{M}$ the maximal ideal of $A$ defining $x$.

Since $\mathcal{O}_{\X,x}/(x_j)$ is regular for each $j$, so is
$(A/(x_j))_{\mathfrak{M}}$ (these Noetherian local rings have
isomorphic $\mathfrak{M}$-adic completions by the proof of
\cite[1.2.1]{conrad}). Hence, shrinking $\mU$, we may assume that
$A/(x_j)$ is a domain for each $j$. Then the irreducible
components of $\mE\cap \mU$ are defined by $x_j=0$ for
$j=0,\ldots,m$, $M_j\neq 0$, and since normalization commutes with
open immersions \cite[1.2.3]{conrad}, $\mE_i\cap \mU$ is a union
of such components. However, since $\mathcal{O}_{\mE_i,x}$ is a
domain, we see that it is of the form $\mathcal{O}_{\X,x}/(x_j)$
for some $j$. In particular, $\mathcal{O}_{\mE_i,x}$ is regular.
%
%
  \end{proof}

Hence, if $\mE$ is a scheme, condition (2) follows from condition
(1) (since any local ring of an irreducible scheme is a domain).


 We say that a closed formal subscheme $\mE$ of a regular special formal $R$-scheme $\X$ is a
strict normal crossings divisor if its defining ideal sheaf is a
strict normal crossings ideal. We say that a special formal
$R$-scheme $\mY$ has strict normal crossings if $\mY$ is regular
and the special fiber $\mY_s$ is a strict normal crossings
divisor, i.e. if the ideal sheaf $\pi\mathcal{O}_{\mY}$ is a
strict normal crossings ideal.

Now let $\X$ be a regular special formal $R$-scheme, and let
$\mathcal{I}$ be a strict normal crossings ideal on $\X$, defining
a closed formal subscheme $\mE$ of $\X$. We can associate to each
irreducible component $\mE_i$ of $\mE$ a multiplicity $m(\mE_i)$
as follows. Choose any point $x$ on $\mE_i$, denote by
$\mathfrak{P}_{i}$
the defining ideal sheaf of $\mE_i$, and by $\mathfrak{P}_{i,x}$ its stalk at $x$.
\begin{lemma}\label{dvr}
The ring $(\mathcal{O}_{\X,x})_{\mathfrak{P}_{i,x}}$ is a DVR.
\end{lemma}
\begin{proof}
The ring $(\mathcal{O}_{\X,x})_{\mathfrak{P}_{i,x}}$ is regular
since $\X$ is regular, so it suffices to show that
$\mathfrak{P}_{i,x}$ is principal. However, if $(x_0,\ldots,x_m)$
is a regular system of local parameters
 in $\mathcal{O}_{\X,x}$ such that
$\mathcal{I}_x$ is generated by $\prod_{i=0}^{m}(x_i)^{M_i}$, then
we've seen in the proof of Lemma \ref{def-equiv} that
$\mathfrak{P}_{i,x}$ is generated by $x_j$ at $x$, for some index
$j$.   \end{proof}

 We
define the multiplicity $m(\mE_i,x)$ of $\mE_i$ at $x$ as the
length of the $(\mathcal{O}_{\X,x})_{\mathfrak{P}_{i,x}}$-module
$(\mathcal{O}_{\mE,x})_{\mathfrak{P}_{i,x}}$.

\begin{lemma-definition}
The multiplicity $m(\mE_i,x)$ does not depend on the point $x$.
Therefore, we denote it by $m(\mE_i)$, and we call it the
multiplicity of $\mE_i$ in $\mE$.

If $x$ is any point of $\mE_i$, and if $\mathcal{I}_x$ is
generated by $\prod_{j=0}^{m}x_j^{M_j}$ in $\mathcal{O}_{\X,x}$,
with $(x_0,\ldots,x_m)$ a regular system of local parameters in
$\mathcal{O}_{\X,x}$, then $\mathfrak{P}_{i,x}$ is generated by
$x_j$ for some index $j$, and $m(\mZ_i)=M_j$.
\end{lemma-definition}
\begin{proof}
We've seen in the proof of Lemma \ref{dvr} that
$\mathfrak{P}_{i,x}$ is generated by $x_j$ for some index $j$, and
that $(\mathcal{O}_{\X,x})_{\mathfrak{P}_{i,x}}$ is a DVR with
uniformizing parameter $x_j$, so clearly $m(\mE_i,x)=M_j$.

Moreover, there exists an open neighbourhood $\mU$ of $x$ in $\X$
such that $x_0,\ldots,x_m$ are defined on $\mU$, and such that
$\mathfrak{P}_i$ is generated by $x_j$ on $\mU$ and $\mathcal{I}$
by $\prod_{j=0}^{m}x_j^{M_j}$. This shows that $m(\mE_i,y)=M_j$
for each point $y$ in a sufficiently small neighbourhood of $x$.
Hence, $m(\mE_i,y)$ is locally constant on $\mE_i$, and therefore
constant since $\mE_i$ is connected.   \end{proof}

If $\X$ is a regular special formal $R$-scheme and $\mE$ is strict
normal crossings divisor, then we write $\mE=\sum_{i\in
I}N_i\mE_i$ to indicate that $\mE_i$, $i\in I$, are the
irreducible components of $\mE$, and that $N_i=m(\mE_i)$ for each
$i$. We say that $\mE$ is a tame strict normal crossings divisor
if the multiplicities $N_i$ are prime to the characteristic
exponent of the residue field $k$ of $R$. We say that $\X$ has
tame strict normal crossings if $\X$ is regular and $\X_s$ is a
tame strict normal crossings divisor.

If $Z$ is a separated $R$-scheme of finite type, $Z$ is regular,
and its special fiber $Z_s$ is a (tame) strict normal crossings
divisor (in the classical sense), then we say that $Z$ has (tame)
strict normal crossings.

For any non-empty subset $J$ of $I$, we define $\mE_J:=\cap_{i\in
J}\mE_i$ (i.e. $\mE_J$ is defined by the sum of the defining ideal
sheaves of $\mE_i$, $i\in J$), $E_J:=\cap_{i\in J} (\mE_i)_0$ and
$E_J^{o}:=E_J\setminus (\cup_{i\notin J}(\mE_i)_0)$. Moreover, we
put $$m_J:=gcd\{N_i\,|\,i\in J\}$$ If $i\in I$, we write $E_i$
instead of $E_{\{i\}}=(\mE_i)_0$. Note that $\mE_J$ is regular and
$E_J=(\mE_J)_0$ for each non-empty subset $J$ of $I$.

\begin{example}\label{notirr}
Consider the special formal $R$-scheme
$$\X=\mathrm{Spf}\,R[[x,y]]/(\pi-x^{N_1}y^{N_2})$$
Then $\X_s=\mathrm{Spf}\,k[[x,y]]/(x^{N_1}y^{N_2})$, and we get
$\X_s=N_1\mE_1+N_2\mE_2$ with $\mE_1=\mathrm{Spf}\,k[[y]]$ and
$\mE_2=\mathrm{Spf}\,k[[x]]$. Note that $E_1^o=E_2^o=\emptyset$,
while $E_{\{1,2\}}$ is a point (the maximal ideal $(\pi,x,y)$).

The varieties $E_i$ are not necessarily irreducible. Consider, for
instance, the smooth special formal $R$-scheme
$$\X=\mathrm{Spf}\,R[[x]]\{y,z\}/(x-yz)$$
Then $\X_s$ is the formal $k$-scheme
$\mathrm{Spf}\,k[[x]][y,z]/(x-yz)$ which is irreducible, since
$k[[x]][y,z]/(x-yz)$ has no zero-divisors. However,
$\X_0=\mathrm{Spec}\,k[y,z]/(yz)$ is reducible.   \end{example}

If $m_J$ is prime to the characteristic exponent of $k$, we
construct an \'etale cover $\widetilde{E}_J^o$ of $E_J^o$ as
follows: we can cover $E_J^o$ by affine open formal subschemes
$\mU=\mathrm{Spf}\,V$ of $\X$ such that $\pi=uv^{m_J}$ with
$u,v\in V$ and $u$ a unit. We put
$$\widetilde{\mU}=\mathrm{Spf}\,V[T]/(uT^{m_J}-1)$$
The restrictions of $\mU$ over $E_J^o$ glue together to an \'etale
cover $\widetilde{E}_J^o$ of $E_J^o$.

If $\X$ is a $stft$ formal $R$-scheme, then all these definitions
coincide with the usual ones (see \cite{NiSe}). In particular,
$\mE_i=E_i$ is a regular $k$-variety, and $N_i$ is the length of
the local ring of the $k$-scheme $\X_s$ at the generic point of
$E_i$.

\begin{definition}
Let $\X$ be a regular special formal $R$-scheme with strict normal
crossings, with $\X_s=\sum_{i\in I}N_i\mE_i$, and let $J$ be a
non-empty subset of $I$. We say that an integer $d>0$ is
$J$-linear if there exist integers $\alpha_j>0$, $j\in J$, with
$d=\sum_{j\in J}\alpha_j N_j$. We say that $d$ is $\X_s$-linear if
$d$ is $J$-linear for some non-empty subset $J$ of $I$ with
$|J|>1$ and $E_J^o\neq \emptyset$.
\end{definition}

\begin{lemma}\label{linear}
Let $\X$ be regular special formal $R$-scheme with strict normal
crossings.
 There exists a sequence of admissible
blow-ups $$\pi^{(j)}:\X^{(j+1)}\rightarrow \X^{(j)},\quad
j=0,\ldots,r-1$$ such that
\begin{itemize}
\item $\X^{(0)}=\X$, \item the special fiber of $\X^{(j)}$ is a
strict normal crossings divisor
$$\X_s^{(j)}=\sum_{i\in I^{(j)}} N_i^{(j)}\mE_i^{(j)},$$ \item $\pi^{(j)}$
is the formal blow-up with center $\mE^{(j)}_{J^{(j)}}$, for some
 subset $J^{(j)}$ of $I^{(j)}$, with $|J^{(j)}|>1$, \item $d$ is not
$\X_s^{(r)}$-linear.
\end{itemize}
\end{lemma}
\begin{proof}
The proof of \cite[5.17]{NiSe2} carries over verbatim to this
setting.   \end{proof}

\begin{definition}
A resolution of singularities of a generically smooth flat
$R$-variety $X$ (resp. a generically smooth, flat special formal
$R$-scheme), is a proper morphism of flat $R$-varieties (resp. a
morphism of flat special formal $R$-schemes) $h:X'\rightarrow X$,
such that $h$ induces an isomorphism on the generic fibers, and
such that $X'$ is regular, with as special fiber a strict normal
crossings divisor $X'_s$. We say that the resolution $h$ is
\textit{tame} if $X'_s$ is a tame strict normal crossings divisor.
%
\end{definition}

\begin{lemma}\label{completion}
Let $A$ be a special $R$-algebra, let $X$ be a $stft$ formal
$R$-scheme, and let $Z$ be a closed subscheme of $X_s$. Finally,
let $U$ be a Noetherian scheme, and let $V$ be a closed subscheme
of $U$.

If $\mathfrak{M}$ is an open prime ideal of $A$, defining a point
$x$ of $\mathrm{Spf}\,A$, then the local morphism
$A_{\mathfrak{M}}\rightarrow \mathcal{O}_{\mathrm{Spf}\,A,x}$
induces an isomorphism on the completions (w.r.t. the respective
maximal ideals)
$$\widehat{A}_{\mathfrak{M}}\cong\widehat{\mathcal{O}}_{\mathrm{Spf}\,A,x} $$
If $x$ is a point of $Z$, then the local morphism
$\mathcal{O}_{X,x}\rightarrow \mathcal{O}_{\widehat{X/Z},x} $
induces an isomorphism on the completions
$$\widehat{\mathcal{O}}_{X,x}\cong\widehat{\mathcal{O}}_{\widehat{X/Z},x} $$
If $x$ is a point of $V$, then the local morphism
$\mathcal{O}_{U,x}\rightarrow \mathcal{O}_{\widehat{U/V},x} $
induces a canonical isomorphism on the completions
$$\widehat{\mathcal{O}}_{U,x}\cong\widehat{\mathcal{O}}_{\widehat{U/V},x} $$
\end{lemma}
\begin{proof}
The first point is shown in the proof of \cite[1.2.1]{conrad}. As
for the second, if $J$ is the defining ideal sheaf of $Z$ on $X$,
then
$$(\mathcal{O}_{X,x})/(J^{n})\cong \mathcal{O}_{X/J^n,x}\cong
(\mathcal{O}_{\widehat{X/Z},x})/(J^{n})$$ for each $n\geq 1$. The
proof of the third point is analogous.   \end{proof}

\begin{lemma}\label{reg}
Let $A$ be a special $R$-algebra, let $X$ be a separated scheme of
finite type over $R$, or a $stft$ formal $R$-scheme, and let $Z$
be a closed subscheme of $X_s$. Let $U$ be a Noetherian
$R$-scheme, and let $V$ be a closed subscheme of $U$.
\begin{enumerate}
\item $\mathrm{Spec}\,A$ is regular iff $\mathrm{Spf}\,A$ is
regular. Moreover, $X$ is regular at the points of $Z$ iff
$\widehat{X/Z}$ is regular, and $U$ is regular at the points of
$V$ iff $\widehat{U/V}$ is regular. \item If
$(\mathrm{Spec}\,A)_s$ is a strict normal crossings divisor, then
$(\mathrm{Spf}\,A)_s$ is a strict normal crossings divisor.
Moreover, if $X_s$ is a strict normal crossings divisor at the
points of $Z$, then $(\widehat{X/Z})_s$ is a strict normal
crossings divisor\footnote{The converse implication is false, as
is seen by taking a regular flat formal curve $X$ over $R$ whose
special fiber $X_s$ is an irreducible curve with a node $x$, and
putting $Z=\{x\}$.}. The same holds if we replace $X$ by $U$ and
$Z$ by $V$, and if we assume that $\widehat{U/V}$ is a special
formal $R$-scheme and that $U$ is excellent. \item
$\mathrm{Spf}\,A$ is generically smooth, iff $(A\otimes_R
K)_{\mathfrak{M}}$ is geometrically regular over $K$, for each
maximal ideal $\mathfrak{M}$ of $A\otimes_R K$. Moreover, if $X$
is generically smooth, then $\widehat{X/Z}$ is generically smooth.
\item If $K$ is perfect, any regular special formal $R$-scheme
$\X$ is generically smooth.
\end{enumerate}
\end{lemma}
\begin{proof}
Regularity of a local Noetherian ring is equivalent to regularity
of its completion \cite[17.1.5]{ega4.1}, so (1) follows from Lemma
\ref{completion}, the fact that regularity can be checked at
maximal ideals, and the fact that any maximal ideal of an adic
topological ring is open. Point (2) follows from the fact that,
for any local Noetherian ring $S$, a tuple $(x_0,\ldots,x_m)$ in
$S$ is a regular system of local parameters for $S$ iff it is a
regular system of local parameters for $\widehat{S}$. The only
delicate point is that we have to check if condition (2) in
Definition \ref{sncd} holds for $(\widehat{X/Z})_s$ and
$(\widehat{U/V})_s$. This, however, follows from Corollary
\ref{irred-compl}.

Now we prove (3). By \cite[2.8]{formrigIII}, smoothness of
$\X_\eta$ is equivalent to geometric regularity of
$\mathcal{O}_{\X_\eta,x}$ over $K$, for each point $x$ of
$\X_\eta$. By \cite[7.1.9]{dj-formal}, $x$ corresponds canonically
to a maximal ideal $M$ of $A\otimes_R K$, and the completions of
$\mathcal{O}_{\X_\eta,x}$ and $(A\otimes_R K)_M$ are isomorphic.
We can conclude by  using the same arguments as in
\cite[2.4(3)]{NiSe}. To prove the second part of (3), we may
assume that $X$ is a $stft$ formal $R$-scheme: if $X$ is a
generically smooth separated scheme of finite type over $R$, then
its $\pi$-adic completion $\widehat{X}$ is generically smooth by
\cite[2.4(3)]{NiSe}, and we have
$\widehat{X/Z}=\widehat{\widehat{X}/Z}$. If $X$ is a $stft$ formal
$R$-scheme, (3) follows from the fact that
$(\widehat{X/Z})_{\eta}$ is canonically isomorphic to the tube
$]Z[$, which is an open rigid subvariety of $X_\eta$.

For point (4), we may assume that $\X=\mathrm{Spf}\,A$ is affine.
It suffices to show that $(A\otimes_R K)_{\mathfrak{M}}$ is
geometrically regular over $K$ for each maximal ideal
$\mathfrak{M}$, by point (3). But $A$ is regular by (1), and since
$K$ is perfect, $(A\otimes_R K)_{\mathfrak{M}}$ is geometrically
regular over $K$.   \end{proof}

\begin{prop}\label{affineres}
If $k$ has characteristic zero, any affine generically smooth flat
special formal $R$-scheme $\X=\mathrm{Spf}\,A$ admits a 
resolution of singularities by means of admissible blow-ups.
\end{prop}
\begin{proof}
By Temkin's resolution of singularities for quasi-excellent
schemes of characteristic zero \cite{temkin-resol},
$\mathrm{Spec}\,A$ admits a resolution of singularities
$Y\rightarrow \Spec A$ by means of blow-ups whose centers contain
a power of $\pi$. Completing w.r.t. an ideal of definition $I$ of
$A$ yields a resolution $h:\mY\rightarrow \Spf A$ by means of
admissible blow-ups, by Lemma \ref{reg}(2).   \end{proof}

\subsection{Etale morphisms of special formal schemes}
We define \textit{\'etale} and \textit{adic \'etale} morphisms of
formal $R$-schemes as in \cite[2.6]{formal1}. A local homomorphism
of local rings $(A,\mathfrak{M})\rightarrow (B,\mathfrak{N})$ is
called unramified if $\mathfrak{N}=\mathfrak{M}B$ and
$B/\mathfrak{N}$ is separable over $A/\mathfrak{M}$. We recall the
following criterion.

\begin{lemma}\label{etalecrit}
Let $h:\mathfrak{Y}\rightarrow \mathfrak{X}$ be a morphism of
pseudo-finite type of Noetherian adic formal schemes, and let $y$
be a point of $\mathfrak{Y}$. Then the following properties are
equivalent:
\begin{enumerate}
\item The local homomorphism
$h^*:\mathcal{O}_{\mathfrak{X},h(y)}\rightarrow
\mathcal{O}_{\mathfrak{Y},y}$ is flat and unramified. \item The
local homomorphism
$\widehat{h^*}:\widehat{\mathcal{O}}_{\mathfrak{X},h(y)}\rightarrow
\widehat{\mathcal{O}}_{\mathfrak{Y},y}$ is flat and unramified.
\item $h$ is \'etale at $y$. \end{enumerate}
\end{lemma}
\begin{proof}
Use \cite{formal2},\,(3.1),\,(4.5)\,,(6.5).   \end{proof}

In (2), the completions can be taken either w.r.t. the adic
topologies on $\X$ and $\mY$, or w.r.t. the topologies defined by
the respective maximal ideals.


\begin{lemma}\label{tamecov}
Let $\X$ be a regular special formal $R$-scheme (or a regular
$R$-variety), such that $\X_s$ is a strict normal crossings
divisor. Let $h:\mathfrak{Y}\rightarrow \X$ be an \'etale morphism
of adic formal schemes. Then $\mathfrak{Y}$ is regular, and
$\mathfrak{Y}_s$ is a strict normal crossings divisor. It is tame
if $\X_s$ is tame.
\end{lemma}
\begin{proof}
Let $y$ be a closed point on $\mathfrak{Y}_s$, and put $x=h(y)$.
Take a regular system of local parameters $(x_0,\ldots,x_m)$ in
$\mathcal{O}_{\X,x}$, such that $\pi=u\prod_{i=0}^m x_i^{N_i}$,
with $u$ a unit.

By Lemma \ref{etalecrit}, $h^*:\mathcal{O}_{\X,x}\rightarrow
\mathcal{O}_{\mathfrak{Y},y}$ is flat and unramified. In
particular,
 $(h^{*}x_0,\ldots,h^{*}x_m)$ is a regular system of local
parameters in $\mathcal{O}_{\mathfrak{Y},y}$. It satisfies
$$\pi=h^{*}u\cdot\prod_{i=0}^m (h^{*}x_i)^{N_i}$$

Finally, it is clear that the irreducible components of $\mY_s$
are the connected components of the regular closed formal
subschemes $\mY_s\times_{\X_s} \mE_j$, with $\mE_j$, $j\in J$ the
irreducible components of $\X_s$.   \end{proof}


\section{Computation of nearby cycles on formal
schemes}\label{sec-nearby}
\subsection{Algebraic covers}\label{mumford}

\begin{definition}\label{mumdiag}
Let $\X$ be a flat special formal $R$-scheme. A nice algebraizable
cover for $\X$ at a closed point $x$ of $\X_0$ is a surjective
finite adic \'etale morphism of special formal $R$-schemes
$\mZ\rightarrow \mU$, with $\mU$ a Zariski-open neighbourhood of
$\X$, such that $\mZ_0/\mU_0$ is a tame \'etale covering, and such
that each point of $\mZ$ has a Zariski-open neighbourhood which is
isomorphic to the formal completion of a regular $R$-variety $Z$
with tame strict normal crossings, along a closed subscheme of
$Z_s$.
\end{definition}

\begin{prop}\label{rev}
We assume that $k$ is perfect. Let $\X$ be a regular special
formal $R$-scheme with tame strict normal crossings. Then $\X$
admits a nice algebraizable cover at any closed point $x$ of
$\X_0$.
\end{prop}
\begin{proof}
We may assume that $\X$ is affine, say $\X=\mathrm{Spf}\,A$, and
that there exist elements $x_0,\ldots,x_m$ in $A$ with
$\pi=u\prod_{i=0}^{m}x_i^{N_i}$, with $u$ a unit and $N_i\in \N$,
and such that $(x_0,\ldots,x_m)$ is a regular system of local
parameters on $\X$ at $x$. Put $d:=gcd(N_0,\ldots,N_m)$, and
consider the finite \'etale morphism
$$g:\mZ:=\mathrm{Spf}\,A[T]/(uT^d-1)\rightarrow \X$$
 Since $\X_s$ is tame, $d$ is prime to
the characteristic exponent $p$ of $k$, and $\mZ_0$ is a tame
\'etale covering of $\X_0$.

 By Bezout's Lemma,
we can find integers $a_j,\,j=0,\ldots,m$, such that
$d=\sum_{j=0}^{p}a_jN_j$.
 On $\mZ$, we
have $\pi=\prod_{j=0}^{m}(T^{-a_j}x_j)^{N_j}$. The sections
$z_j:=T^{-a_j}x_j$ on $\mZ$ define a morphism of formal schemes
over $\mathrm{Spec}\,R$
$$h:\mZ\rightarrow Y:=\mathrm{Spec}\,R[y_0,\ldots,y_m]/(\pi-\prod_{j=0}^{m}
y_j^{N_j})$$ given by $h^*(y_i)=z_i$, mapping the fiber over $x$
to the origin. Moreover, $(z_0,\ldots,z_m)$ is a regular system of
local parameters at any point $z$ of $\mZ$ lying over $x$, and by
Lemma \ref{etalecrit}, $h$ is \'etale at $z$. Hence, shrinking
$\X$, we may assume that $h$ is \'etale on $\mZ$.

By \cite[7.12]{formal2}, if $z$ is any point of $\mZ$ lying over
$x$, there exists a Zariski-open neighbourhood of $z$ in $\mZ$
which is a formal completion of an adic \'etale $Y$-scheme $Z$
along a closed subscheme of $Z_s$; $Z$ is automatically a regular
$R$-variety with tame strict normal crossings, by Lemma
\ref{tamecov}.
  \end{proof}

In particular, every regular special formal $R$-scheme with tame
strict normal crossings is algebraizable locally w.r.t. the
\'etale topology. Combining this with Proposition \ref{affineres},
we see that, if $k$ has characteristic zero and $\X$ is a
generically smooth affine special formal $R$-scheme, then there
exist a morphism of special formal $R$-schemes $h:\X'\rightarrow
\X$ such that $h_\eta$ is an isomorphism, and an \'etale cover
$\{\mU_i\}$ of $\X'$ by algebraizable special formal $R$-schemes
$\mU_i$.

\subsection{Computation of the nearby cycles on tame
strict normal crossings} We can use the constructions in the
preceding section to generalize Grothendieck's computation of the
tame nearby cycles on a tame strict normal crossings divisor
\cite[Exp. I]{sga7a}, to the case of special formal $R$-schemes
\cite{berk-vanish2}. Let $\X$ be a regular special formal
$R$-scheme, such that the special fiber $\X_s$ is a tame strict
normal crossings divisor $\sum_{i\in I} N_i \mE_i$. For any
non-empty subset $J$ of $I$, we denote by $M_J$ the kernel of the
linear map
$$\Z^{J}\rightarrow \Z:(z_j)\mapsto \sum_{j\in J}N_jz_j.$$

\begin{prop}\label{nearby}
Suppose that $k$ is algebraically closed. Let $\X$ be a regular
special formal $R$-scheme, such that $\X_s$ is a tame strict
normal crossings divisor $\sum_{i\in I} N_i \mE_i$. Let $M$ be a
torsion ring, with torsion orders prime to the characteristic
exponent of $k$. For each non-empty subset $J\subset I$, and each
$i\geq 0$, the $i$-th cohomology sheaf of tame nearby cycles
$R^i\psi_{\eta}^t(M)$ associated to $\X$, is tamely lisse on
$E_J^o$. Moreover, for each $i>0$, and each point $x$ on $E_J^o$,
there are canonical isomorphisms
$$R^i\psi_{\eta}^t(M)_x=R^0\psi_{\eta}^t(M)_x\otimes \bigwedge^i
M_J^{\vee}$$ and $R^0\psi_{\eta}^t(M)_x\cong M^{F_J},$ where $F_J$
is a set of cardinality $m_J$, on which $G(K^t/K)$ acts
transitively.
%
\end{prop}
\begin{proof}
 By \cite{berk-vanish2}, Cor. 2.3, whenever
$h:\mY\rightarrow \X$ is an adic \'etale morphism of special
formal $R$-schemes, we have
$$R^i\psi^t_{\eta}(M|_{\mY_\eta})\cong
h_0^{*}R^i\psi_\eta^t(M|_{\X_{\eta}})$$ Hence, by Proposition
\ref{rev}, we may suppose that $\X$ is the formal completion of a
regular $R$-variety $X$ with tame strict normal crossings, along a
closed subscheme of $X_s$.
 By Berkovich' Comparison Theorem
\cite[5.1]{Berk-vanish}, it suffices to prove the corresponding
statements for $X$ instead of $\X$. The computation of the fibers
was done in \cite[Exp. I]{sga7a}. To see that $R^i\psi_\eta^t(M)$
is tamely lisse on $E_J^o$, one can argue as follows: by the
arguments in the proof of Proposition \ref{rev}, one can reduce to
the case
$$X=\Spec R[x_0,\ldots,x_m]/(\pi-x_0^{N_0}\cdot\ldots\cdot
x_q^{N_q})$$ with $q\leq m$, $N_j>0$ for each $j$, and with
$E_J^o$ defined by $x_0=\ldots=x_q=0$. By smooth base change, one
reduces to the case where $m=q$, and $E_J^o$ is the origin. Now
the statement is trivial.   \end{proof}
\begin{corollary}\label{computation}
Suppose that $k$ is algebraically closed. Let $\X$ be a regular
special formal scheme over $R$, such that the special fiber $\X_s$
is a tame strict normal crossings divisor $\sum_{i\in I} N_i
\mE_i$.  Let $x$ be a closed point on $E_J^o$, for some non-empty
subset $J\subset I$, and let $e>0$ be an integer. Let $\varphi$ be
a topological generator of the tame geometric monodromy group
$G(K^t/K)$.
\begin{itemize}
\item $Tr(\varphi^e\,|\,R\psi^t_{\eta}(\Q_{\ell})_x)=0$, if
$|J|>1$, or if $J$ is a singleton $\{i\}$ and $N_i\nmid e$, \item
$Tr(\varphi^e\,|\,R\psi^t_{\eta}(\Q_{\ell})_x)=N_i$, if $J=\{i\}$
and $N_i|e$.
\end{itemize}
\end{corollary}
\begin{corollary}
If $k$ is algebraically closed, of characteristic zero, then
$R\psi^t_\eta(\Q_\ell)$ is constructible on $\X_0$, for any
generically smooth special formal $R$-scheme $\X$, and the action
of $G(K^{t}/K)$ is continuous.
\end{corollary}
\begin{proof}
This follows from Proposition \ref{affineres}, Proposition
\ref{nearby}, and Corollary \cite[2.3]{berk-vanish2}.
  \end{proof}
\section{Motivic integration on special formal schemes}\label{sec-mot}
Throughout this section, we assume that $k$ is perfect.
A possible approach to define motivic integration on special
formal $R$-schemes $\X$, would be to introduce the Greenberg
scheme $Gr^R(\X)$ of $\X$ (making use of the fact that $V(J)$ is
of finite type over $R$ if $J$ is an ideal of definition on $\X$)
and to generalize the constructions in \cite{sebag1} and
\cite{motrigid} to this setting.
 We will take a shortcut, instead, making use of appropriate $stft$
models for special formal $R$-schemes. The theory of motivic
integration on $stft$ formal $R$-schemes was developed in
\cite{sebag1}, and this theory was used to define motivic
integrals of differential forms of maximal degree on smooth
quasi-compact rigid varieties in \cite{motrigid}. The
constructions were refined to a relative setting in \cite{NiSe2},
and extended to so-called ``bounded'' rigid varieties in
\cite{NiSe-weilres}.

\begin{definition}
Let $\X$ be a special formal $R$-scheme. A N\'eron smoothening for
$\X$ is a morphism of special formal $R$-schemes $\mY\rightarrow
\X$ such that $\mY$ is adic smooth over $R$, and such that
$\mY_\eta\rightarrow \X_\eta$ is an open embedding satisfying
$\mY_\eta(K^{sh})=\X_\eta(K^{sh})$.
\end{definition}

\begin{remark}
If $\X$ is $stft$ over $R$, we called this a \textbf{weak} N\'eron
smoothening in \cite{NiSe2}, to make a distinction with a stronger
variant of the definition. Since this distinction is irrelevant
for our purposes, we omit the adjective ``weak'' from our
terminology.   \end{remark}

Any generically smooth $stft$ formal $R$-scheme $\X$ admits a
N\'eron smoothening $\mY\rightarrow \X$, by \cite[3.1]{formner}.
Moreover, by \cite[6.1]{NiSe2}, the class $[\mY_s]$ of $\mY_s$ in
$\mathcal{M}_{\X_0}/(\LL-[\X_0])$ does not depend on $\mY$.
Following \cite{motrigid}, we called this class the motivic Serre
invariant of $\X$, and denoted it by $S(\X)$. We will generalize
this result to special formal $R$-schemes.

\begin{lemma}\label{unibounded}
If $\X$ is a flat special formal $R$-scheme, and
$h:\mathfrak{Y}\rightarrow \X$ is the dilatation with center
$\X_0$, then $\mathfrak{Y}$ is $stft$ over $R$, and $h$ induces an
open embedding $\mathfrak{Y}_\eta\rightarrow \X_\eta$ with
$\mY_\eta(K^{sh})=\X_\eta(K^{sh})$.
\end{lemma}
\begin{proof}
By the universal property in Proposition \ref{univ-dilat}, any
$R^{sh}$-section of $\X$ lifts uniquely to $\mathfrak{Y}$.
\end{proof}
\begin{prop}
Any generically smooth special formal $R$-scheme $\X$ admits a
N\'eron smoothening.
\end{prop}
\begin{proof}
We may assume that $\X$ is flat over $R$. Take $\mY$ as in Lemma
\ref{unibounded}; it is generically smooth, since $\mY_\eta$ is
open in $\X_\eta$. If $\mY'\rightarrow \mY$ is a N\'eron
smoothening of $\mY$, then the composed morphism $\mY'\rightarrow
\X$ is a N\'eron smoothening for $\X$.   \end{proof}

\begin{lemma}\label{dilat-change}
Let $\X$ be a flat, generically smooth $stft$ formal $R$-scheme,
and let $U$ be a closed subscheme of $\X_s$. If we denote by
$\mY\rightarrow \X$ the dilatation with center $U$, then the image
of $S(\mY)$ under the forgetful morphism
$$\mathcal{M}_{\mY_0}/(\LL-[\mY_0])\rightarrow \mathcal{M}_{U}/(\LL-[U])$$
 coincides with the image of $S(\X)$
under the base change morphism
$$\mathcal{M}_{\X_0}/(\LL-[\X_0])\rightarrow \mathcal{M}_{U}/(\LL-[U])$$
Likewise, if $\omega$ is a differential form of maximal degree
(resp. a gauge form) on $\X_\eta$, then the image of
$\int_{\mY}|\omega|$ coincides with the image of
$\int_{\X}|\omega|$ in $\widehat{\mathcal{M}}_{U}$, resp.
$\mathcal{M}_{U}$.
\end{lemma}
\begin{proof}
If we denote by $h:\mY'\rightarrow \X$ be the blow-up of $\X$ with
center $U$, then $\mY$ is, by definition, an open formal subscheme
of $\mY'$. If $g:\mZ'\rightarrow \mY'$ is a N\'eron smoothening,
and if we put $\mZ=g^{-1}(\mY)$, then by the universal property of
the dilatation in Proposition \ref{univ-dilat}, the induced open
embedding of Greenberg schemes $Gr^R(\mZ)\rightarrow Gr^R(\mZ')$
is an isomorphism onto the cylinder $(h_s\circ g_s\circ
\theta_0)^{-1}(U)$ (here $\theta_0$ denotes the truncation
morphism $Gr^R(\mZ')\rightarrow \mZ'_s$).   \end{proof} We recall
that the motivic integral $\int_{\X}|\omega|$, with $\X$
generically smooth and $stft$ over $R$, was defined in \cite[\S
6]{NiSe2}, refining the construction in \cite[4.1.2]{motrigid}
\begin{definition}\label{def-serre}
Let $\X$ be a generically smooth, flat special formal $R$-scheme,
and let $h:\X'\rightarrow \X$ be the dilatation with center
$\X_0$. We define the motivic Serre invariant $S(\X)$ of $\X$ by
$$S(\X):=S(\X')\mbox{ in }\mathcal{M}_{\X_0}/(\LL-[\X_0])$$
If $\omega$ is a differential form of maximal degree (resp. a
gauge form) on $\X_\eta$, then we put
$$\int_{\X}|\omega|:=\int_{\X'}|\omega|$$
in $\widehat{\mathcal{M}}_{\X_0}$, resp. $\mathcal{M}_{\X_0}$.

If $\X$ is a generically smooth special formal $R$-scheme, we
denote by $\X^{flat}$ its maximal flat closed subscheme (obtained
by killing $\pi$-torsion) and we put $S(\X)=S(\X^{flat})$ and
$$\int_{\X}|\omega|:=\int_{\X^{flat}}|\omega|$$
\end{definition}
It follows from Lemma \ref{dilat-change} that this definition
coincides with the usual one if $\X$ is $stft$ over $R$. Note that
even in this case, the dilatation $\X'\rightarrow \X$ is not
necessarily an isomorphism, since $\X_s$ might not be reduced. In
fact, $\X'$ might be empty, for instance if $\X$ has strict normal
crossings with all multiplicities $>1$.

\begin{prop}\label{int-dilat}
Let $\X$ be a generically smooth special formal $R$-scheme. If
$\mY\rightarrow \X$ is a N\'eron smoothening for $\X$, then
$$S(\X)=[\mY_0] \quad \in \mathcal{M}_{\X_0}/(\LL-[\X_0])$$ For any
differential form of maximal degree (resp. gauge form) $\omega$ on
$\X_\eta$, we have
$$\int_{\X}|\omega|=\int_{\mY}|\omega|$$
in $\widehat{\mathcal{M}}_{\X_0}$, resp. $\mathcal{M}_{\X_0}$.
\end{prop}
\begin{proof}
We may assume that $\X$ is flat over $R$; let $\X'\rightarrow \X$
be the dilatation with center $\X_0$. By the universal property of
the dilatation in Proposition \ref{univ-dilat} and the fact that
$\mY_s$ is reduced, any N\'eron smoothening $\mY\rightarrow \X$
factors through a morphism of $stft$ formal $R$-schemes
$\mY\rightarrow \X'$, and by Lemma \ref{unibounded}, this is again
a N\'eron smoothening. So the result follows from (the proof of)
\cite[6.11]{NiSe2} (see also \cite{NiSe-weilres} for an addendum
on the mixed dimension case).   \end{proof}

We showed in \cite{NiSe-weilres} that the generic fiber $\X_\eta$
of a generically smooth special formal $R$-scheme $\X$ is a
so-called ``bounded'' rigid variety over $K$ (this also follows
immediately from Lemma \ref{unibounded}), and we defined the
motivic Serre invariant $S(\X_\eta)$ of $\X_\eta$, as well as
motivic integrals $\int_{\X_\eta}|\omega|$ of differential forms
$\omega$ on $\X_\eta$ of maximal degree.
\begin{prop}\label{specialize}
If $\X$ is a generically smooth special formal $R$-scheme, then
$S(\X_\eta)$ is the image of $S(\X)$ under the forgetful morphism
$$\mathcal{M}_{\X_0}/(\LL-[\X_0])\rightarrow \mathcal{M}_k/(\LL-1)$$ If $\omega$ is a gauge form on $\X_\eta$,
then $S(\X)$ is the image of $\int_{\X}|\omega|$ in
$\mathcal{M}_{\X_0}/(\LL-[\X_0])$. If $\omega$ is a differential
form of maximal degree (resp. a gauge form) on $\X_\eta$, then
$\int_{\X_\eta}|\omega|$ is the image of $\int_{\X}|\omega|$ under
the forgetful morphism $\widehat{\mathcal{M}}_{\X_0}\rightarrow
\widehat{\mathcal{M}}_{k}$, resp. $\mathcal{M}_{\X_0}\rightarrow
\mathcal{M}_k$.
\end{prop}
\begin{proof}
This is clear from the definitions, and the corresponding
properties for $stft$ formal $R$-schemes \cite[6.4]{NiSe2}.
\end{proof}

\begin{definition}\label{def-volpoin}
Suppose that $k$ has characteristic zero. Let $\X$ be a
generically smooth special formal $R$-scheme. Let $\omega$ be a
gauge form on $\X_\eta$. We define the volume Poincar\'e series of
$(\X,\omega)$ by
$$S(\X,\omega;T):=\sum_{d>0}\left(\int_{\X(d)}|\omega(d)|\right)T^d\
\mathrm{in}\ \mathcal{M}_{\X_0}[[T]]$$
\end{definition}
\begin{remark}
The motivic Serre invariant $S(\X)$ and the motivic integral
$\int_{\X}|\omega'|$ (for any differential form $\omega'$ of
maximal degree on $\X_\eta$) are independent of the choice of
uniformizer $\pi$. The volume Poincar\'e series, however, depends
on the choice of $\pi$, or more precisely, on the $K$-fields
$K(d)$. If $k$ is algebraically closed, then $K(d)$ is the unique
extension of degree $d$ of $K$, up to $K$-isomorphism, and
$S(\X,\omega;T)$ is independent of the choice of $\pi$. See also
the remark following Definition \ref{motvol2}.   \end{remark}

\begin{prop}\label{tube-volume}
 Let $\X$ be a
generically smooth special formal $R$-scheme, and let $U$ be a
locally closed subscheme of $\X_0$. Denote by $\mathfrak{U}$ the
formal completion of $\X$ along $U$. Then $S(\mathfrak{U})$ is the
image of $S(\X)$ under the base change morphism
$$\mathcal{M}_{\X_0}/(\LL-[\X_0])\rightarrow \mathcal{M}_U/(\LL-[U])$$

If $\omega$ is a gauge form on $\X_\eta$, then
$\int_{\mU}|\omega|$ is the image of $\int_{\X}|\omega|$ under the
base change morphism
$$\mathcal{M}_{\X_0}\rightarrow \mathcal{M}_U$$ (the analogous
statement holds if $\omega$ is merely a differential form of
maximal degree).

If, moreover, $k$ has characteristic zero, then
$S(\mathfrak{U},\omega;T)$ is the image of $S(\X,\omega;T)$ under
the base-change morphism
$$\mathcal{M}_{\X_0}[[T]]\rightarrow
\mathcal{M}_U[[T]]$$

In particular, if $\X$ is $stft$ over $R$, then $S(\mathfrak{U})$
and $S(\mU,\omega;T)$ coincide with the invariants with support
$S_U(\X)$ and $S_U(\X,\omega;T)$ defined in \cite{NiSe}.
%
\end{prop}
\begin{proof}
We may assume that $\X$ is flat. If $U$ is open in $\X_0$, then
these results are clear from the definitions, since dilatations
commute with flat base change; so we may suppose that $U$ is a
reduced closed subscheme of $\X_0$. Let $h:\X'\rightarrow \X$ be
the dilatation with center $\X_0$, and let $\mU'\rightarrow \mU$
be the dilatation with center $\mU_0=U$. By Proposition
\ref{dilat-comm}, there exists a unique morphism of $stft$ formal
$R$-schemes $\mU'\rightarrow \X'$ such that the square
$$\begin{CD}
\mU'@>>>\mU
\\@VVV @VVV
\\\X'@>>> \X
\end{CD}$$
commutes, and $\mU'\rightarrow \X'$ is the dilatation with center
$\X'_s\times_{\X_0} U$. Now we can conclude by Lemma
\ref{dilat-change}.   \end{proof}
\begin{cor}\label{strat}
If $\{U_i,\,i\in I\}$ is a finite stratification of $\X_0$ into
locally closed subsets, and $\mU_i$ is the formal completion of
$\X$ along $U_i$, then \begin{eqnarray*} S(\X)&=&\sum_{i\in I}
S(\mU_i) \\\int_{\X}|\omega|&=&\sum_{i\in I}\int_{\mU_i}|\omega|
\\ S(\X,\omega;T)&=&\sum_{i\in I}S(\mU_i,\omega;T)
\end{eqnarray*}
(we applied the forgetful morphisms $\mathcal{M}_{U_i}\rightarrow
\mathcal{M}_{\X_0}$ to the right-hand sides).
\end{cor}

\begin{definition}
Let $\X$ be a special formal $R$-scheme. A \textbf{special}
N\'eron smoothening for $\X$ is a morphism of special formal
$R$-schemes $\mY\rightarrow \X$ such that $\mY$ is smooth over
$R$, and such that $\mY_\eta\rightarrow \X_\eta$ is an open
embedding satisfying $\mY_\eta(K^{sh})=\X_\eta(K^{sh})$.
\end{definition}
In particular, any N\'eron smoothening is a special N\'eron
smoothening (which makes the terminology somewhat paradoxical).
\begin{prop}\label{special-neron}
If $\X$ is a generically smooth flat special formal $R$-scheme,
then there exists a composition of admissible blow-ups
$\mY\rightarrow \X$ such that $Sm(\mY)\rightarrow \X$ is a special
N\'eron smoothening for $\X$.
\end{prop}
\begin{proof}
Let $\mathcal{I}$ be the largest ideal of definition for $\X$, and
let $h:\X'\rightarrow \X$ be the admissible blow-up with center
$\mathcal{I}$. We denote by $\mU$ the open formal subscheme of
$\X'$ where $\mathcal{I}\mathcal{O}_{\X'}$ is generated by $\pi$,
i.e. $\mU\rightarrow \X$ is the dilatation with center $\X_0$, and
$\mU$ is $stft$ over $R$.

By \cite[3.1]{formner} and \cite[2.5]{formrigI}, there exists an
admissible blow-up $g:\mY'\rightarrow \mU$ such that
$Sm(\mY')\rightarrow \mU$ is a N\'eron smoothening. We will show
that this blow-up extends to an admissible blow-up $\mY\rightarrow
\X'$. Choose an integer $j>0$ such that $\pi^j$ is contained in
the center $\mathcal{J}$ of $g$. By \cite[9.4.7]{ega1}, the
pull-back of $\mathcal{J}$ to $V(\mathcal{I}^j\mathcal{O}_{\mU})$
extends to a coherent ideal sheaf $\mathcal{J}'$ on
$V(\mathcal{I}^j\mathcal{O}_{\X'})$; we'll denote again by
$\mathcal{J}'$ the corresponding coherent ideal sheaf on $\X'$.
Since formal blow-ups commute with open embeddings, the admissible
blow-up $\mY\rightarrow \X'$ with center $\mathcal{J}'$ extends
$g$.

Finally, let us show that the composed morphism
$Sm(\mY)\rightarrow \X$ is a special N\'eron smoothening. It
suffices to show that the natural map
$Sm(\mY)_\eta(K^{sh})\rightarrow \X_\eta(K^{sh})$ is surjective.
By Lemma \ref{unibounded}, $\mU_\eta(K^{sh})\rightarrow
\X_\eta(K^{sh})$ is surjective, and since $Sm(\mY')\rightarrow
\mU$ is a N\'eron smoothening, $Sm(\mY')_\eta(K^{sh})\rightarrow
\mU_\eta(K^{sh})$ is surjective; but $\mY'$ is an open formal
subscheme of $\mY$, so the result follows.   \end{proof}
\begin{prop}\label{smooth-serre}
If $\X$ is a smooth special formal $R$-scheme, then
$$S(\X)=[\X_0]$$
in $\mathcal{M}_{\X_0}/(\LL-[\X_0])$.
\end{prop}
\begin{proof}
Stratifying $\X_0$ in regular pieces, we might as well assume that
$\X_0$ is regular from the start, by Corollary \ref{strat}.
Moreover, we may suppose that $\X=\mathrm{Spf}\,A$ is affine, and
that $\X_0$ is defined by a regular sequence
$(\pi,x_1,\ldots,x_q)$ in $A$. The dilatation of $\X$ with center
$\X_0$ is given by
$$\mY:=\mathrm{Spf}\,A\{T_1,\ldots,T_q\}/(x_i-\pi
T_i)_{i=1,\ldots,q}\rightarrow \X$$ Now $\mY$ is flat, and
$$\mY_s=\mY_0=\mathrm{Spec}\,\left(A/(\pi,x_1,\ldots,x_q)\right)[T_1,\ldots,T_q]\cong \X_0\times_k \A^{q}_k$$
Since $k$ is perfect and $\X_0$ is regular, $\mY_0$ is smooth over
$k$, and hence $\mY$ is smooth over $R$, and
$$S(\X)=S(\mY)=[\mY_0]=[\X_0]$$ in
$\mathcal{M}_{\X_0}/(\LL-[\X_0])$.   \end{proof}
\begin{corollary}\label{serreneron}
If $h:\mY\rightarrow \X$ is a special N\'eron smoothening, then
$$S(\X)=[\mY_0]\in \mathcal{M}_{\X_0}/(\LL-[\X_0])$$
\end{corollary}

\section{Computation of Serre invariants and motivic
integrals}\label{sec-comput} Throughout this section, we assume
that $k$ is perfect.
\subsection{Serre invariants of the ramifications}
If $\X$ is a special formal $R$-scheme, we denote by
$\widetilde{\X}\rightarrow \X$ the normalization morphism
\cite{conrad}.
\begin{theorem}\label{neronsmooth}
Let $\X$ be a regular special formal $R$-scheme, such that $\X_s$
is a tame strict normal crossings divisor $\sum_{i\in I}N_i\mE_i$,
and let $d>0$ be an integer, prime to the characteristic exponent
of $k$. If $d$ is not $\X_s$-linear, then
$$h:Sm(\widetilde{\X(d)})\rightarrow \X(d)$$ is a special N\'eron
smoothening. Moreover, if we put
$$\widetilde{E}(d)_i^o=(\widetilde{\X(d)})_0\times_{\X_0}E_i^o$$
for each $i$ in $I$, then
$$Sm(\widetilde{\X(d)})_0=\bigsqcup_{N_i|d}\widetilde{E}(d)_i^o$$
and the $E_i^o$-variety $\widetilde{E}(d)_i^o$ is canonically
isomorphic to $\widetilde{E}_i^o$ (defined in Section
\ref{strict}).
\end{theorem}
\begin{proof}
The fact that
$$Sm(\widetilde{\X(d)})_0=\bigsqcup_{N_i|d}\widetilde{E}(d)_i^o$$
and that the $E_i^o$-variety $\widetilde{E}(d)_i^o$ is canonically
isomorphic to $\widetilde{E}_i^o$, can be proven exactly as in
\cite[4.4]{NiSe}.

 Since $\X(d)_\eta$ is smooth and, a fortiori, normal, and
normalization commutes with taking generic fibers
\cite[2.1.3]{conrad}, $h_\eta$ is an open embedding. Since $d$ is
not $\X_s$-linear, the obvious generalization of
\cite[5.15]{NiSe2} implies that
$$\left(Sm(\widetilde{\X(d)})\right)_\eta(K(d)^{sh})=\X(d)_\eta(K(d)^{sh})$$
i.e. $h$ is a special N\'eron smoothening.   \end{proof}
\begin{cor}\label{computserre}
Let $\X$ be a regular special formal $R$-scheme, such that $\X_s$
is a strict normal crossings divisor $\sum_{i\in I}N_i\mE_i$, and
let $d>0$ be an integer, prime to the characteristic exponent of
$k$.

 If $d$ is not $\X_s$-linear, then
 $$S(\X(d))=\sum_{i\in I, N_i|d}[\widetilde{E}^o_i]$$ in
 $\mathcal{M}_{\X_0}/(\LL-[\X_0])$.
\end{cor}
\begin{proof}
If $\X'$ is the completion of $\X$ along $\sqcup_{N_i|d}E_i^o$,
then $S(\X(d))=S(\X'(d))$, since $d$ is not $\X_s$-linear, by a
straightforward generalization of \cite[5.15]{NiSe2}. Hence, we
may as well assume that $\X'=\X$. In this case, since $d$ is prime
to the characteristic exponent of $k$, $\X_s$ is tame, so we can
use Theorem \ref{neronsmooth} and Corollary \ref{serreneron} to
conclude.   \end{proof}
\subsection{Order of a top form at a section}\label{ordersect}
This subsection is a straightforward generalization of \cite[\S
6.2]{NiSe}. Let $\X$ be a generically smooth special formal
$R$-scheme, of pure relative dimension $m$. Let $R'$ be a finite
extension of $R$, of ramification index $e$, and denote by $K'$
its quotient field.

\begin{definition}\label{orderideal}
For any element $\psi$ of $\X(R')$, and any ideal sheaf
$\mathcal{I}$ on $\X$, we define $ord(\mathcal{I})(\psi)$ as the
length of the $R'$-module $R'/\psi^*\mathcal{I}$.
\end{definition}

We recall that the length of the zero module is $0$, and the
length of $R'$ is $\infty$.

For any element $\psi$ of $\X(R')$, the $R'$-module
$M:=(\psi^{*}\Omega^m_{\X/R})/(\mathrm{torsion})$ is a free rank
$1$ sublattice of the rank $1$ $K'$-vector space
$(\psi_\eta)^{*}(\Omega^m_{\X_\eta/K})$.

\begin{definition}
For any global section $\omega$ of $\Omega^m_{\X/R}$ and any
section $\psi$ in $\X(R')$, we define the order
$ord(\omega)(\psi)$ of $\omega$ at $\psi$ as follows: we choose an
integer $a\geq 0$ such that $\omega':=\pi^a\psi_\eta^*(\omega)$
belongs to the
 sublattice $M$ of $(\psi_\eta)^{*}(\Omega^m_{\X_\eta/K})$, and we
 put $$ord(\omega)(\psi)=\mathrm{length}_{R'}(M/R'\omega')-e.a$$
 This definition does not depend on $a$.
\end{definition}

If $e=1$, this definition coincides with the one given in
\cite[4.1]{motrigid}. It only depends on the completion of $\X$ at
$\psi(0)\in \X_0$.
%
%
If
$\omega$ is a gauge form on $\X_\eta$,
$ord(\omega)(\psi)$ is finite.


Now let $h:\mZ\rightarrow \X$ be a morphism of generically smooth
special formal $R$-schemes, both of pure relative dimension $m$.
Let $R'$ be a finite extension of $R$, and fix a section $\psi$ in
$\mZ(R')$. The canonical morphism
$$h^*\Omega^m_{\X/R}\rightarrow \Omega^m_{\mZ/R}$$
induces a morphism of free rank $1$ $R'$-modules
$$(\psi^*h^*\Omega^m_{\X/R})/(\mathrm{torsion})\rightarrow (\psi^*\Omega^m_{\mZ/R})/(\mathrm{torsion})$$
We define $ord(Jac_h)(\psi)$ as the length of its cokernel.

If $\Omega^m_{\mZ/R}/(\mathrm{torsion})$ is a locally free rank
$1$ module over $\mathcal{O}_{\mZ}$, we define the Jacobian ideal
sheaf $\mathcal{J}ac_h$ of $h$ as the annihilator of the cokernel
of the morphism
$$h^*\Omega^m_{\X/R}\rightarrow \Omega^m_{\mZ/R}/(\mathrm{torsion})$$
and we have $ord(Jac_h)(\psi)=ord(\mathcal{J}ac_h)(\psi)$. The
following lemmas are proved as their counterparts \cite{NiSe},
Lemma 6.4-5.

\begin{lemma}\label{jac}
Let $h:\mZ\rightarrow \X$ be a morphism of generically smooth
special formal $R$-schemes, both of pure relative dimension $m$.
Let $R'$ be a finite extension of $R$.
 For any global section $\omega$ of $\Omega^m_{\X/K}$, and any
 $\psi'\in \mZ(R')$,
$$ord(h^*\omega)(\psi)=ord(\omega)(h(\psi))+ord(Jac_h)(\psi)$$
\end{lemma}

\begin{lemma}\label{arcchange}
Let $e\in \N^{\ast}$ be prime to the characteristic exponent of
$k$. Let $\X$ be a regular special formal $R$-scheme with tame
strict normal crossings.
 Let $\omega$ be a global section of
$\Omega^m_{\X_\eta/K}$. We denote by $\widetilde{\omega(e)}$ the
pullback of $\omega$ to the generic fiber of $\widetilde{\X(e)}$.
Let $R'$ be a finite extension of $R(e)$, and let $\psi(e)$ be a
section in $Sm(\widetilde{\X(e)})(R')$. If we denote by $\psi$ its
image in $\X(R')$, then
$$ord(\omega)(\psi)=ord(\widetilde{\omega(e)})(\psi(e))$$
\end{lemma}

\subsection{Order of a gauge form on a smooth formal $R$-scheme}\label{sec-ordergauge} Let $\X$
be a smooth special formal $R$-scheme of pure relative dimension
$m$, and let $\omega$ be a $\X$-\textit{bounded} gauge form on
$\X_\eta$ (see Definition \ref{mbounded}). We denote by
$Irr(\X_0)$ the set of irreducible components of $\X_0$ (note that
$\X_0$ is not always smooth over $k$; see Example \ref{notirr}).
Let $C$ be an irreducible component of $\X_0$, and let $\xi$ be
its generic point. The local ring $\mathcal{O}_{\X,\xi}$ is a UFD,
since $\X$ is regular. Since $\X$ is smooth over $R$,
$(\Omega^m_{\X/R})_{\xi}$ is a free rank $1$ module over
$\mathcal{O}_{\X,\xi}$ by \cite[4.8]{formal1} and
\cite[5.10]{formal2}, and $\pi$ is irreducible in
$\mathcal{O}_{\X,\xi}$.

\begin{definition}\label{def-ord}
Let $A$ be a UFD, and let $a$ be an irreducible element of $A$.
Let $N$ be a free $A$-module of rank one, and let $n$ be an
element of $N$. We choose an isomorphism of $A$-modules $A\cong
N$, and we define $ord_a n$ as follows: if $n\neq 0$, $ord_a n$ is
the largest $q\in \N$ such that $a^q|n$ in $A$. If $n=0$, we put
$ord_a n=\infty$. This definition does not depend on the choice of
isomorphism $A\cong N$.
\end{definition}
\begin{definition}
Let $\X$ be a smooth special formal $R$-scheme of pure relative
dimension $m$, let $C$ be an irreducible component of $\X_0$, and
denote by $\xi$ its generic point. If $\omega$ is a $\X$-bounded
$m$-form on $\X_\eta$, we can choose $b\in \N$ such that
$\pi^b\omega$ extends to a section $\omega'$ of
$(\Omega^m_{\X/R})_{\xi}$, and we put
$$ord_{C}\omega:=ord_{\pi}\omega'-b $$
This definition does not depend on $b$.
\end{definition}

\begin{lemma}\label{unit}
Let $\X$ be a smooth connected special formal $R$-scheme, and let
$f$ be an element of $\mathcal{O}_{\X}(\X)$. If $f$ is a unit on
$\X_\eta$, and $f$ is not identically zero on $\X_0$, then $f$ is
a unit on $\X$.
\end{lemma}
\begin{proof}
We may assume that $\X=\mathrm{Spf}\,A$ is affine. By the
correspondence between maximal ideals of $A\otimes_R K$ and points
of $\X_\eta$ explained in \cite[7.1.9]{dj-formal}, we see that $f$
is a unit in $A\otimes_R K$, so there exists an element $q$ in $A$
and an integer $i\geq 0$ such that $fq=\pi^i$. We may assume that
either $i=0$ (in which case $f$ is a unit), or $q$ is not
divisible by $\pi$ in $A$. Since $\X$ is smooth, $\pi$ is a prime
in $A$, so $\pi$ divides $f$ if $f$ is not a unit; this
contradicts the hypothesis that $f$ does not vanish identically on
$\X_0$.   \end{proof}
\begin{lemma}\label{ord-section}
Let $\X$ be a smooth special formal $R$-scheme, and let $\omega$
be a $\X$-bounded gauge form on $\X_\eta$. Let $R'$ be a finite
unramified extension of $R$, and consider a section $\psi\in
\X(R')$. If $C$ is an irreducible component of $\X_0$ containing
$\psi(0)$, then
$$ord_{C}(\omega)=ord(\omega)(\psi)$$
\end{lemma}
\begin{proof}
We denote by $\xi$ the generic point of $C$. Multiplying with
powers of $\pi$, we may assume that $\omega$ is defined on $\X$.
Moreover, we may assume that there exists a section $\omega_0$ in
$\Omega^m_R(\X/R)$ which generates $\Omega^m_{\X/R}$ at each point
of $\X$, and we write $\omega=f\omega_0$ with $f\in
\mathcal{O}_{\X}(\X)$. Dividing by an appropriate power of $\pi$,
we may assume that $\pi\nmid f$ in $\mathcal{O}_{\X,\xi}$ and
$ord_C(\omega)=0$. Since $\omega$ is gauge on $\X_\eta$, $f$ is a
unit on $\X_\eta$, so by Lemma \ref{unit}, $f$ is a unit on $\X$.
Hence, $$ord(\omega)(\psi)=ord_{\pi}\psi^*(f)=0$$   \end{proof}
\begin{cor}
If $C_1$ and $C_2$ are irreducible components of the same
connected component $C$ of $\X_0$, then
$ord_{C_1}(\omega)=ord_{C_2}(\omega)$, and we denote this value by
$ord_C(\omega)$.
%
\end{cor}
\begin{proof}
We can always find a section $\psi$ in $\X(R')$ with $R'/R$ finite
and unramified, and with $\psi(0)\in C_1\cap C_2$.
\end{proof}
\begin{cor}\label{ord-complete}
If $\X_0$ is connected, and $Z$ is a locally closed subset of
$\X_0$, then
$$ord_{C}(\omega)=ord_{\X_0}(\omega)$$ for any connected component
$C$ of $Z$, where the left hand side is computed on the completion
$\widehat{\X/Z}$.
\end{cor}
\begin{lemma}\label{ord-dilat}
Let $\X$ be a smooth connected special formal $R$-scheme, and let
$Z$ be a regular closed subscheme of $\X_0$. If we denote by
$h:\mY\rightarrow \X$ the dilatation with center $Z$, and by $c$
the codimension of $Z$ in $\X$, then for any $\X$-bounded gauge
form $\omega$ on $\X_\eta$,
$$ord_{\X_0}(\omega)=ord_{\mY_0}(\omega)+c-1$$
\end{lemma}
\begin{proof}
By Corollary \ref{ord-complete} and Proposition \ref{dilatprop},
we may assume that $\X_0=Z$. Now the result follows from Lemma
\ref{ord-section} and Lemma \ref{jac}, since
$\mathcal{J}ac_h=(\pi^{c-1})$.
  \end{proof}

\begin{prop}\label{comput-smooth}
Let $\X$ be a smooth special formal $R$-scheme, of pure relative
dimension $m$, and denote by $\mathcal{C}(\X_0)$ the set of
connected components of $\X_0$. For any $\X$-bounded gauge form
$\omega$ on $\X_\eta$,
$$\int_{\X}|\omega|=\LL^{-m}\sum_{C\in
\mathcal{C}(\X_0)}[C]\LL^{-ord_C(\omega)}$$ in
$\mathcal{M}_{\X_0}$.
\end{prop}
\begin{proof}
By Corollaries \ref{strat} and \ref{ord-complete}, we may assume
that $\X_0$ is regular and connected. By definition,
$$\int_{\X}|\omega|=\int_{\mY}|\omega|$$
where $\mY\rightarrow \X$ is the dilatation with center $\X_0$. In
the proof of Proposition \ref{smooth-serre}, we saw that $\mY$ is
smooth, and $\mY_0=[\X_0]\LL^{c-1}$, with $c$ the codimension of
$\X_0$ in $\X$. Hence, we can conclude by Lemma \ref{ord-dilat}.
  \end{proof}
\begin{cor}\label{volume-smooth}
Let $\X$ be a generically smooth special formal $R$-scheme ofpure
relative dimension $m$, and let $\omega$ be a gauge form on
$\X_\eta$. If $\mY\rightarrow \X$ is a special N\'eron
smoothening, then
$$\int_{\X}|\omega|=\LL^{-m}\sum_{C\in
\mathcal{C}(\mY_0)}[C]\LL^{-ord_C(\omega)}$$ in
$\mathcal{M}_{\X_0}$.
\end{cor}

\section{The trace formula}
Let $F$ be any field, let $Z$ be a variety over $F$, and let $A$
be an abelian group. The abelian group $C(Z,A)$ of constructible
$A$-functions on $Z$ is the subgroup of the abelian group of
functions of sets $Z\rightarrow A$, consisting of mappings of the
form
$$\varphi=\sum_{S\in \mathcal{S}}a_S.\mathbf{1}_{S}$$
where $\mathcal{S}$ is a finite stratification of $Z$ into
constructible subsets, and $a_S\in A$ for $S\in \mathcal{S}$, and
where $\mathbf{1}_{S}$ denotes the characteristic function of $S$.
Note that a constructible function on $Z$ is completely determined
by its values on the set of closed points $Z^o$ of $Z$. If $A$ is
a ring, $C(Z,A)$ carries a natural ring structure. If $A=\Z$, we
call $C(Z,A)$ the ring of constructible functions on $Z$, and we
denote it by $C(Z)$.

For any constructible $A$-function
$$\varphi=\sum_{S\in\mathcal{S}}a_S.\mathbf{1}_{S}$$ on $Z$ as above, we can define its
integral w.r.t. the Euler characteristic as follows:
\begin{eqnarray*}
\int_{Z}\varphi d\chi&:=&\sum_{S\in\mathcal{S}}a_S\chi_{top}(S)
\end{eqnarray*}
If the group operation on $A$ is written multiplicatively, we
write $\int_Z^{\times}$ instead of $\int_Z$. The calculus of
integration with respect to the Euler characteristic was, to our
knowledge, first introduced in \cite{viro}.

\begin{lemma}\label{sheaf}
Suppose that $F$ is algebraically closed, and let $Z$ be a variety
over $F$. Let $\ell$ be a prime number, invertible in $F$, and let
$\mathcal{L}$ be a tamely constructible $\Q_\ell^s$-adic sheaf on
$Z$. Suppose that a finite cyclic group $G$ with generator $g$
acts on $\mathcal{L}$. We denote by $Z^{o}$ the set of closed
points on $Z$.
\begin{enumerate}
\item the mapping $$Tr(g\,|\,\mathcal{L}_*):Z^{o}\rightarrow
\Q_{\ell}^s:x\mapsto Tr(g\,|\,\mathcal{L}_x)$$ defines a
constructible $\Q_\ell^s$-function on $Z$, and
$$Tr(g\,|\,\oplus_{i\geq 0}H^i_c(Z,\mathcal{L}))=\int_{Z}Tr(g\,|\,\mathcal{L}_*)d\chi$$
\item the mapping $$\zeta(g\,|\,\mathcal{L}_*;T):Z^{o}\rightarrow
\Q_{\ell}^s[[T]]^{\times}:x\mapsto \zeta(g\,|\,\mathcal{L}_x;T)$$
defines a constructible $\Q_\ell^s[[T]]^{\times}$-function on $Z$,
and
$$\zeta(g\,|\,\oplus_{i\geq 0}H^i_c(Z,\mathcal{L});T)=\int^{\times}_{Z}\zeta(g\,|\,\mathcal{L}_*;T)d\chi$$
\end{enumerate}

\end{lemma}
\begin{proof}
First, we prove (1). By additivity of $H_{c}(.\,)$, we may suppose
that $Z$ is normal and $\mathcal{L}$ is tamely lisse on $Z$. In
this case, $Tr(g\,|\,\mathcal{L}_*)$ is constant on $Z$, and the
result follows from \cite[5.1]{NiSe}. Now (2) follows from the
identity \cite[1.5.3]{weil1}
$$det(Id-T.M\,|\,V)^{-1}=\mathrm{exp}(\sum_{d>0}Tr(M^d\,|\,V)\frac{T^d}{d})$$
for any endomorphism $M$ on a finite dimensional vector space $V$
over a field of characteristic zero.   \end{proof}
\begin{lemma}\label{laumon}
Let $G$ be a finite group, let $F$ be an algebraically closed
field, and fix a prime $\ell$, invertible in $F$. If
$f:Y\rightarrow X$ is a morphism of separated $F$-schemes of
finite type, and $\mathcal{L}$ is a constructible $E[G]$-sheaf on
$Y$, for some finite extension $E$ of $\Q_\ell$, then
$$f_*[\mathcal{L}]=f_![\mathcal{L}]\ \mathrm{in}\ K_0(X;E[G])$$
In particular, if $X=\mathrm{Spec}\,F$, then
$$Tr(g\,|\,\oplus_{i\geq 0}H^i(Y,\mathcal{L}))=Tr(g\,|\,\oplus_{i\geq 0}H^i_c(Y,\mathcal{L}))$$
for each element $g$ of $G$.
\end{lemma}
\begin{proof}
If $G$ is the trivial group, then this is a well-known theorem of
Laumon's \cite{laumon-euler}. His proof carries over verbatim to
the case where $G$ is any finite group.
%
  \end{proof}
\begin{corollary}\label{illusie}
Let $G=<g>$ be a finite cyclic group. Let $F$ be an algebraically
closed field, let $U$ be a variety over $F$, and let $\mathcal{L}$
be a tamely constructible $\Q^s_\ell[G]$-sheaf on $U$, for any
prime $\ell$ invertible in $F$. Then
\begin{eqnarray*}
Tr(g\,|\,\oplus_{i\geq
0}H^i(U,\mathcal{L}))&=&Tr(g\,|\,\oplus_{i\geq
0}H^i_c(U,\mathcal{L}))
\\&=&\int_U Tr(g\,|\,\mathcal{L}_*)\,d\chi
\end{eqnarray*}
\end{corollary}
\begin{proof}
The first equality follows from Lemma \ref{laumon}, while the
equality between the second and the third expression follows from
Lemma \ref{sheaf}.
  \end{proof} The following theorem is a broad generalization
of \cite[5.4]{NiSe}. It implies, in particular, that the
assumptions that $X$ is algebraic and $Z$ is proper, are
superfluous in the statement of \cite[5.4]{NiSe}.
\begin{theorem}[Trace formula]\label{trace}
Assume that $k$ is perfect. Let $\varphi$ be a topological
generator of the tame geometric monodromy group $G(K^t/K^{sh})$.
Let $\X$ be a generically smooth special formal $R$-scheme, and
suppose that $\X$ admits a tame resolution of singularities
$h:\mY\rightarrow \X$, with $\mY_s=\sum_{i\in I} N_i \mE_i$.  For
any integer $d>0$, prime to the characteristic exponent of $k$, we
have
$$\chi_{top}\left(S(\X_\eta(d))\right)
=Tr(\varphi^d\,|\,H(\,\overline{\X_\eta}\,))=\sum_{N_i |d} N_i
\chi_{top}(E_i^o)$$
\end{theorem}
\begin{proof}
We may assume that $\mY=\X$, and that $k$ is algebraically closed
(since the motivic Serre invariant is clearly compatible with
unramified extensions of the base $R$). The equality
$$\chi_{top}\left(S(\X_\eta(d))\right)
=\sum_{N_i |d} N_i \chi_{top}(E_i^o)$$ can be proven as in
\cite[5.4]{NiSe}: the expression holds if $d$ is not
$\X_s$-linear, by Corollary \ref{computserre} and the fact that
$\widetilde{E}_i^o$ is a degree $N_i$ finite \'etale cover of
$E_i$. Moreover, the right hand side of the equality does not
change under blow-ups with center $\mE_J$ with $\emptyset \neq
J\subset I$, by the obvious generalization of \cite[5.2]{NiSe}, so
the expression holds in general by Lemma \ref{linear}.

So it suffices to prove that
$$Tr(\varphi^d\,|\,H(\,\overline{\X_\eta}\,))=\sum_{N_i |d} N_i
\chi_{top}(E_i^o)$$ By \cite[2.3(ii)]{berk-vanish2}, there is for
each $i\geq 0$ a canonical isomorphism
$$H^i(\X_0,R\psi_\eta^t(\Q_\ell|_{\X_\eta}))\cong
H^i(\,\overline{\X_\eta}\,,\Q_\ell)$$ and hence,
\begin{equation}\label{trace-eq}
Tr(\varphi^d\,|\,H(\,\overline{\X_\eta}\,))=Tr(\varphi^d\,|\,H(\X_0,R\psi_\eta^t(\Q_\ell|_{\X_\eta})))
\end{equation}
By our local computation in Proposition \ref{nearby} we can filter
$R^j\psi^t_\eta(\Q_{\ell}|_{\X_\eta})$ by constructible subsheaves
which are stable under the monodromy action and such that the
action of $\varphi$ on successive quotients has finite order.
Hence, we can apply Lemma \ref{illusie}. Combined with the
computation in Corollary \ref{computation}, we obtain the required
equality.
\end{proof}
\begin{cor}
If $k$ has characteristic zero, and $\X$ is a generically smooth
special formal $R$-scheme, then
$$\chi_{top}\left(S(\X_\eta(d))\right)
=Tr(\varphi^d\,|\,H(\,\overline{\X_\eta}\,))$$ for any integer
$d>0$. In particular,
$$\chi_{top}\left(S(X(d))\right)
=Tr(\varphi^d\,|\,H(\,\overline{X}\,))$$ for any smooth
quasi-compact rigid variety $X$ over $K$.
\end{cor}
\begin{proof}
Since generically smooth affine special formal $R$-schemes admit a
resolution of singularities if $k$ has characteristic zero, by
Proposition \ref{affineres}, we can cover $\X$ by a finite family
of open affine formal subschemes, such that
$$\chi_{top}\left(S(\mV(d))\right)
=Tr(\varphi^d\,|\,H(\,\overline{\mV_\eta}\,))$$ whenever $\mV$ is
an intersection of members of this cover. But both sides of this
equality are additive w.r.t. $\mV$ (for the right hand side, use
equation (\ref{trace-eq}) and Lemma \ref{illusie}), which yields
the result for $\mV=\X$.   \end{proof}
\begin{cor}
If $k$ is an algebraically closed field of characteristic zero,
and $X_K^{an}$ is the analytification of a proper smooth variety
$X_K$ over $K$, then
$$\chi_{top}\left(S(X^{an}_K(d))\right)
=Tr(\varphi^d\,|\,\oplus_{i\geq 0} H^i(X_K\times_K K^s,\Q_\ell))$$
\end{cor}
\begin{proof}
If $X$ is any flat proper $R$-model for $X_K$, then $X_K^{an}$ is
canonically isomorphic to the generic fiber $X_\eta$ of the
$\pi$-adic completion $\widehat{X}$, by \cite[0.3.5]{bert}.
Moreover, by \cite[7.5.4]{Berk-etale}, there is a
canonical isomorphism
$$H^i(X_K\times_K K^s,\Q_\ell)\cong H^i(\overline{X_\eta})$$ for each $i\geq 0$.
  \end{proof}


As we observed in \cite[\S 5]{NiSe}, some tameness condition is
necessary in the statement of the trace formula: if $R$ is the
ring of Witt vectors $W(\mathbb{F}_p^s)$, and $\X$ is
$\mathrm{Spf}\,R[x]/(x^p-\pi)$, then $S(\X)=0$ while the trace of
$\varphi$ on the cohomology of $\overline{\X_\eta}$ is $1$. It
would be interesting to find an intrinsic tameness condition on
$\X_\eta$ under which the trace formula holds.

It would also be interesting to find a proof of the trace formula
which does not rely on explicit computation, and which does not
use resolution of singularities. One could use the following
strategy. There's no harm in assuming that $k$ is algebraically
closed. After admissible blow-up, we may suppose that the
$R$-smooth part $Sm(\X)$ of $\X$ is a weak N\'eron model for $\X$,
by Proposition \ref{special-neron}. On $Sm(\X)_0$, the tame nearby
cycles are trivial, so the trace of $\varphi$ on the cohomology of
$Sm(\X)_\eta$ yields $\chi_{top}(S(\X))$. Hence, it suffices to
prove that the trace of $\varphi$ on the cohomology of
$R\psi^t_\eta(\Q_\ell|_{\X_\eta})|_Y$ vanishes, where $Y$ denotes
the complement of $Sm(\X)_0$ in $\X_0$. We'll assume the following
result: for each $i\geq 0$, there is a canonical
$G(K^t/K)$-equivariant isomorphism
\begin{equation}\label{hypo} H^i(Y,R\psi^t_\eta(\Q_\ell|_{\X_\eta})|_Y)\cong
H^i(\overline{]Y[},\Q_\ell)\end{equation} This is known if $\X$ is
algebraizable, by Berkovich' comparison theorem
\cite[3.1]{berk-vanish2}. If $K$ has characteristic zero and $\X$
is $stft$ over $R$, then it can be proven in the same way as
Huber's result \cite[3.15]{huber-finite} (which is the
corresponding result for Berkovich' functor $R\Theta_K$ from
\cite[\S 2]{berk-vanish2}).

Assuming (\ref{hypo}), everything reduces to the following
assertion: if $\X_\eta$ is the smooth generic fiber of a special
formal $R$-scheme and satisfies $\X_\eta(K)=\emptyset$ and an
appropriate tameness condition (in particular if $k$ has
characteristic zero), then
$$Tr(\varphi\,|\,H(\overline{\X_\eta}))=0$$ At this point, I don't know how to prove
this without making an explicit computation on a resolution of
singularities.

%

\section{The volume Poincar\'e series, and the motivic
volume}\label{sec-volume} Throughout this section, we assume that
$k$ has characteristic zero.
\subsection{Order of a gauge form along strict normal
crossings}\label{ordstrict} Throughout this subsection, $\X$
denotes a regular special formal $R$-scheme of pure relative
dimension $m$, whose special fiber is a strict normal crossings
divisor $\X_s=\sum_{i\in I}N_i\mE_i$.
\begin{definition}\label{omegai}
For each $i\in I$, and each point $x$ of $E_i$, we denote by
$\mathfrak{P}_{i,x}$ the (not necessarily open) prime ideal in
$\mathcal{O}_{\X,x}$ corresponding to $\mE_i$, and we define
$\mathcal{O}_{\X,\mE_i,x}$ as the localization
$$\mathcal{O}_{\X,\mE_i,x}=(\mathcal{O}_{\X,x})_{\mathfrak{P}_{i,x}}$$
Moreover, we introduce the $\mathcal{O}_{\X,\mE_i,x}$-module
$$\Omega_{\X,\mE_i,x}:=(\Omega^m_{\X/R,x})_{\mathfrak{P}_{i,x}}/(\mathcal{O}_{\X,\mE_i,x}-\mathrm{torsion})$$
\end{definition}

By Lemma \ref{dvr}, $\mathcal{O}_{\X,\mE_i,x}$ is a discrete
valuation ring. Note that the valuation of $\pi$ in
$\mathcal{O}_{\X,\mE_i,x}$ equals $N_i$.

If $\X$ is $stft$ over $R$, then $E_i=\mE_i$, and if we denote by
$\xi_i$ the generic point of $E_i$, then $\mathfrak{P}_{i,\xi_i}$
is the maximal ideal of $\mathcal{O}_{\X,\xi_i}$. So
$\mathcal{O}_{\X,\mE_i,\xi_i}=\mathcal{O}_{\X,\xi_i}$, and
$\Omega_{\X,\mE_i,\xi_i}$ is the module $\Omega_i$ considered in
\cite[6.7]{NiSe}. Note that in the general case, $E_i$ is not
necessarily irreducible (see Example \ref{notirr}).

\begin{lemma}\label{etco}
If $h:\mY\rightarrow \X$ is an \'etale morphism, and
$\mathfrak{C}$ is a connected component of $h^{-1}(\mE_i)$, then
the natural map
$$h^*\Omega^m_{\X/R}\rightarrow \Omega^m_{\mY/R}$$
induces an isomorphism
$\Omega_{\X,\mE_i,x}\otimes_{\mathcal{O}_{\X,\mE_i,x}}\mathcal{O}_{\mY,\mathfrak{C},y}\cong
\Omega_{\mY,\mathfrak{C},y}$ for each point $y$ on
$C=\mathfrak{C}_0$ and with $x=h(y)$.
\end{lemma}
\begin{proof}
Since $h$ is \'etale,
$$h^*\Omega^m_{\X/R}\rightarrow \Omega^m_{\mY/R}$$ is an
isomorphism by \cite[4.10]{formal1}. Let $\mathfrak{P}'_y$ be the
prime ideal in $\mathcal{O}_{\mY,y}$ defining $\mathfrak{C}$.
 Since $h$ is
\'etale, the local morphism $h^{*}:\mathcal{O}_{\X,x}\rightarrow
\mathcal{O}_{\mY,y}$ is a flat, unramified monomorphism by Lemma
\ref{etalecrit}, and by localization, so is
$$\mathcal{O}_{\X,\mE_i,x}\rightarrow
\mathcal{O}_{\mY,\mathfrak{C},y}$$  The isomorphism
$$\Omega^m_{\X/R,x}\otimes_{\mathcal{O}_{\X,x}}\mathcal{O}_{\mY,y}\cong
\Omega^m_{\mY/R,y}$$ localizes to an isomorphism
$$(\Omega^m_{\X/R,x})_{\mathfrak{P}_{i,x}}\otimes_{\mathcal{O}_{\X,\mE_i,x}}\mathcal{O}_{\mY,\mathfrak{C},y}\cong
(\Omega^m_{\mY/R,y})_{\mathfrak{P}'_y}$$ which, by flatness of
$\mathcal{O}_{\X,\mE_i,x}\rightarrow
\mathcal{O}_{\mY,\mathfrak{C},y}$, induces an isomorphism
$$\Omega_{\X,\mE_i,x}\otimes_{\mathcal{O}_{\X,\mE_i,x}}\mathcal{O}_{\mY,\mathfrak{C},y}\cong
\Omega_{\mY,\mathfrak{C},y}$$   \end{proof}

\begin{lemma}\label{rank}
For each $i\in I$ and each point $x$ of $E_i$, the
$\mathcal{O}_{\X,\mE_i,x}$-module $\Omega_{\X,\mE_i,x}$ is free of
rank one.
\end{lemma}
\begin{proof}
Since $\Omega_{\X,\mE_i,x}$ is finite over
$\mathcal{O}_{\X,\mE_i,x}$ and torsion-free, and
$\mathcal{O}_{\X,\mE_i,x}$ is a PID, the module
$\Omega_{\X,\mE_i,x}$ is free over $\mathcal{O}_{\X,\mE_i,x}$. Let
us show that its rank equals $1$. By Lemma \ref{etco}, we may pass
to an \'etale cover and assume that there exists a regular system
of local parameters $(x_0,\ldots,x_n)$ in $\mathcal{O}_{\X,x}$
with $\mathfrak{P}_{i,x}=(x_i)$ and
$\pi=\prod_{i=0}^{n}x_i^{N_i}$. Deriving this expression, we see
that $\Omega_{\X,\mE_i,x}$ is generated by $dx_0\wedge \ldots
\wedge\widehat{dx_i}\wedge \ldots \wedge dx_n$.   \end{proof}
%

Note that the natural map $\Omega^m_{\X/R}(\X)\rightarrow
\Omega_{\X,\mE_i,x}$ factors through
$$\Omega^m_{\X/R}(\X)/(\pi-\mathrm{torsion})$$ since
$\Omega_{\X,\mE_i,x}$ has no torsion.

\begin{definition}\label{order}
Fix $i\in I$ and let $x$ be a point of $E_i$. For any $$\omega \in
\Omega^m_{\X/R}(\X)/(\pi-\mathrm{torsion})$$ we define the order
of $\omega$ along $\mE_i$ at $x$ as the length of the
$\mathcal{O}_{\X,\mE_i,x}$-module
$$\Omega_{\X,\mE_i,x}/(\mathcal{O}_{\X,\mE_i,x}\cdot \omega)$$
 and we denote it by $ord_{\mE_i,x}\omega$.

If $\omega$ is a $\X$-bounded $m$-form on $\X_\eta$, there exists
an integer $a\geq 0$ and an affine open formal subscheme $\mU$ of
$\X$ containing $x$, such that $\pi^a\omega$ belongs to
$$\Omega^m_{\X/R}(\mU)/(\pi-\mathrm{torsion})\subset
\Omega^m_{\X/R}(\mU)\otimes_R K$$ We define the order of $\omega$
along $\mE_i$ at $x$ as
$$ord_{\mE_i,x}\omega:=ord_{\mE_i,x}(\pi^a\omega)-aN_i$$
\end{definition}

This definition does not depend on $a$. If $\X$ is smooth, it
coincides with the one given in Section \ref{sec-ordergauge} in
the following sense: if $\X$ is connected and $x$ is any point of
$\X_0$, then $ord_{\X_0}\omega=ord_{\X_s,x}\omega$.
\begin{lemma}\label{cons-lemma}
Fix $i$ in $I$, and let $x$ and $y$ be points of $E_i$ such that
$y$ belongs to the Zariski-closure of $\{x\}$. For any
$\X$-bounded $m$-form $\omega$ on $\X_\eta$,
$$ord_{\mE_i,x}\omega=ord_{\mE_i,y}\omega$$
\end{lemma}
\begin{proof}
We may suppose that $\omega\in \Omega^m_{\X/R}(\X)$. The natural
localization map $\mathcal{O}_{\X,y}\rightarrow
\mathcal{O}_{\X,x}$ induces a flat, unramified local homomorphism
$\mathcal{O}_{\X,\mE_i,y}\rightarrow \mathcal{O}_{\X,\mE_i,x}$,
and $$\Omega_{\X,\mE_i,x}\cong
\Omega_{\X,\mE_i,y}\otimes_{\mathcal{O}_{\X,\mE_i,y}}\mathcal{O}_{\X,\mE_i,x}$$
Hence, we can conclude by the following algebraic property: if
$g:A\rightarrow A'$ is a flat, unramified morphism of discrete
valuation rings, if $M$ is a free $A$-module of rank $1$, and $m$
is an element of $M$, then the length of the $A$-module $M/(Am)$
equals the length of the $A'$-module $(M\otimes A')/(A'm)$.
Indeed: fixing an isomorphism of $A$-modules $A\cong M$, the
length of $A/(Am)$ is equal to the valuation of $m$ in $A$.
%
  \end{proof}
\begin{cor}\label{cons}
Fix $i$ in $I$. If $\omega$ is a $\X$-bounded $m$-form on
$\X_\eta$, then the function
$$E_i\rightarrow \Z:x\mapsto ord_{\mE_i,x}\omega$$
is constant on $E_i$.
\end{cor}
\begin{definition}
For any $i\in I$ and any $\X$-bounded $m$-form on $\X_\eta$, we
define the order of $\omega$ along $\mE_i$ by
$$ord_{\mE_i}\omega:=ord_{\mE_i,x}\omega$$
where $x$ is any point on $E_i$.
\end{definition}
By Corollary \ref{cons}, this definition does not depend on the
choice of $x$.

\begin{lemma}\label{ord0}
Let $h:\mY\rightarrow \X$ be an \'etale morphism of special formal
$R$-schemes. Let $\mathfrak{C}$ be a connected component of
$h^{-1}(\mE_i)$. For any $\X$-bounded $m$-form $\omega$ on
$\X_\eta$,
$$ord_{\mE_i}\omega=ord_{\mathfrak{C}} h^{*}_\eta\omega$$
\end{lemma}
\begin{proof}
This follows immediately from Lemma \ref{etco} and the algebraic
argument in the proof of Lemma \ref{cons-lemma}.   \end{proof}

%

For each $i\in I$, we denote by $\mathcal{I}_{\mE_i}$ the defining
ideal sheaf of $\mE_i$ in $\X$. For any finite extension $R'$ of
$R$ and any $\psi\in \X(R')$, we denote
$ord(\mathcal{I}_{\mE_i})(\psi)$ by $ord_{\mE_i}(\psi)$ (see
Definition \ref{orderideal}). If $R'$ has ramification degree $e$
over $R$, and the closed point $\psi(0)$ of the section $\psi$ is
contained in $E_i^o$, then the equality
$\pi=x_i^{N_i}\cdot\mathrm{(unit)}$ in $\mathcal{O}_{\X,\psi(0)}$
implies that $ord_{\mE_i}(\psi)=e/N_i$.

The following results are proven exactly as their counterparts
\cite[6.11-13]{NiSe}.
\begin{lemma}\label{arc}
Fix a non-empty subset $J$ of $I$.
Let $R'$
be a finite extension of $R$, 
and let
$\psi$ be an element of $\X(R')$, such that its closed point
$\psi(0)$ lies on
$E_J^o$. 
For any $\X$-bounded gauge form $\omega$ on $\X_\eta$,
$$ord(\omega)(\psi)=\sum_{i\in J}ord_{\mE_i}(\psi)(ord_{\mE_i}\omega-1)+\max_{i\in J}\{ord_{\mE_i}(\psi)\}$$
\end{lemma}

\begin{prop}\label{ord1}
 Let $\omega$ be a $\X$-bounded gauge form on $\X_\eta$. Take a subset $J$ of $I$, with $|J|>1$,
and $E_J^o\neq \emptyset$. Let $h:\X'\rightarrow \X$ be the formal
blow-up with center $\mE_J$, and denote by $\mE'_0$ its
exceptional component. We have
$$ord_{\mE'_0}\omega=\sum_{i\in J}ord_{\mE_i}\omega$$
\end{prop}
\begin{prop}\label{ord2}
 Let
$\omega$ be a $\X$-bounded gauge form on $\X_\eta$. Fix an integer
$e>0$. Denote by $\widetilde{\omega(e)}$ the pull-back of $\omega$
to the generic fiber of $\widetilde{\X (e)}$. For each $i\in I$,
with $N_i|e$, and each connected component $C$ of
$Sm(\widetilde{\X (e)})_0\times_{\X_0}E_i$, we have
$$ord_{C}(\widetilde{\omega(e)})= (e/N_i).ord_{\mE_i}\omega$$
where the left hand side is computed on the smooth special formal
$R$-scheme $Sm(\widetilde{\X (e)})$.
\end{prop}

\subsection{Volume Poincar\'e series}
\begin{theorem}\label{compvolume-d}
Let $\X$ be a regular special formal $R$-scheme of pure relative
dimension $m$, whose special fiber is a strict normal crossings
divisor $\X_s=\sum_{i\in I}N_i\mE_i$. Let $\omega$ be a
$\X$-bounded gauge form on $\X_\eta$, and put
$\mu_i=ord_{\mE_i}\omega$ for each $i\in I$. Then for any integer
$d>0$,
$$\int_{\X(d)}|\omega(d)|=\LL^{-m}\sum_{\emptyset\neq J\subset
I}(\LL-1)^{|J|-1}[\widetilde{E}_J^{o}](\sum_{\stackrel{k_i\geq
1,i\in J}{\sum_{i\in J} k_iN_i=d}}\LL^{-\sum_i k_i \mu_i}\,)\
\mbox{in}\ \mathcal{M}_{\X_0}$$
\end{theorem}
\begin{proof}
We'll show how the proof of the corresponding statement in
\cite[7.6]{NiSe} can be generalized.

First, suppose that $d$ is not $\X_0$-linear. Then it follows from
Theorem \ref{neronsmooth}, Corollary \ref{volume-smooth} and
Proposition \ref{ord2} that
\begin{eqnarray*}
\int_{\X}|\omega(d)|&=&\LL^{-m}\sum_{N_i|d}[\widetilde{E}_i^{o}]\LL^{-
d\mu_i/N_i}
\\&=&\LL^{-m}\sum_{\emptyset\neq J\subset
I}(\LL-1)^{|J|-1}[\widetilde{E}_J^{o}](\sum_{\stackrel{k_i\geq
1,i\in J}{\sum_{i\in J} k_iN_i=d}}\LL^{-\sum_i
k_i\mu_i}\,)\,,\qquad\qquad(*)\end{eqnarray*}

By Lemma \ref{linear}, it suffices to show that the expression
$(*)$ is invariant under formal blow-ups with center $\mE_J$,
$|J|>1$. This can be done as in \cite[7.6]{NiSe}, using an
immediate generalization of the local computation in
\cite[7.5]{NiSe}.
\end{proof}
\begin{corollary}\label{serrat}
Let $\X$ be a generically smooth special formal $R$-scheme, of
pure relative dimension $m$. Suppose that $\X$ admits a resolution
of singularities $\X'\rightarrow \X$, with special fiber
$\X'_s=\sum_{i\in I} N_i \mE_i$. Let $\omega$ be a $\X$-bounded
gauge form on $\X_\eta$.

 The volume
Poincar\'e series $S(\X,\omega;T)$ is rational over
$\mathcal{M}_{\X_0}$. In fact, if we put
$\mu_i:=ord_{\mE_i}\omega$, then the series is given explicitly by
$$S(\X,\omega;T)=\LL^{-m}\sum_{\emptyset\neq J\subset I}(\LL-1)^{|J|-1}[\widetilde{E}_J^{o}]
\prod_{i\in J}\frac{\LL^{-\mu_i}T^{N_i}}{1-\LL^{-\mu_i}T^{N_i}} \
\mathrm{in}\ \mathcal{M}_{\X_0}[[T]]$$
\end{corollary}

By Proposition \ref{affineres}, any affine generically smooth
special formal $R$-scheme admits a resolution of singularities. By
the additivity of the motivic integral, we obtain an expression
for the volume Poincar\'e series in terms of a finite atlas of
local resolutions. In particular, we obtain the following result.

\begin{cor}\label{rationality}
Let $\X$ be a generically smooth special formal $R$-scheme, of
pure relative dimension $m$. Let $\omega$ be a $\X$-bounded gauge
form on $\X_\eta$.
 The volume
Poincar\'e series $S(\X,\omega;T)$ is rational over
$\mathcal{M}_{\X_0}$. More precisely, there exists a finite subset
$S$ of $\Z\times \N^{\ast}$ such that $S(\X,\omega;T)$ belongs to
the subring
$$\mathcal{M}_{\X_0}\left[\frac{\LL^{a}T^b}{1-\LL^{a}T^b} \right]_{(a,b)\in S} $$
of $\mathcal{M}_{\X_0}[[T]]$.
\end{cor}

\subsection{The Gelfand-Leray form and the local singular
series}\label{sec-gelfand} Let $\X$ be a special formal
$R$-scheme, of pure relative dimension $m$. Then $\X$ is also a
formal scheme of pseudo-finite type over $\mathrm{Spec}\,k$, in
the terminology of \cite{formal1}, and the sheaves of continuous
differential forms $\Omega^{i}_{\X/k}$ are coherent, by
\cite[3.3]{formal1}.

Consider the morphism of coherent $\mathcal{O}_\X$-modules
$$\begin{CD}i:\Omega^m_{\X/k}@>d\pi\wedge >>
\Omega^{m+1}_{\X/k}:\omega\mapsto d\pi\wedge\omega\end{CD}$$ Since
$d\pi\wedge \Omega^{m-1}_{\X/k}$ is contained in its kernel, and
$$\Omega^m_{\X/k}/\left(d\pi\wedge \Omega^{m-1}_{\X/k}\right)\cong
\Omega^{m}_{\X/R}$$ by \cite[3.10]{formal1}, $i$ descends to a
morphism of coherent $\mathcal{O}_{\X}$-modules
$$\begin{CD}i:\Omega^m_{\X/R}@>d\pi\wedge >>
\Omega^{m+1}_{\X/k}\end{CD}$$ We've seen in Section \ref{generic}
that any coherent module $\mathcal{F}$ on $\X$ induces a coherent
module $\mathcal{F}_{rig}$ on $\X_\eta$, and this correspondence
is functorial. Hence, $i$ induces a morphism of coherent
$\mathcal{O}_{\X_\eta}$-modules
$$\begin{CD}i:\Omega^m_{\X_\eta/K}@>d\pi\wedge >>
(\Omega^{m+1}_{\X/k})_{rig}\end{CD}$$ by \cite[7.1.12]{dj-formal}.

\begin{definition}
If $\X$ is a special formal $R$-scheme, the Koszul complex of $\X$
is by definition the complex of coherent
$\mathcal{O}_{\X_\eta}$-modules
$$\begin{CD}
0@>>> \mathcal{O}_{\X_\eta}@>d\pi\wedge >>
(\Omega^1_{\X/k})_{rig}@>d\pi\wedge
>>(\Omega^2_{\X/k})_{rig}@>d\pi\wedge >>\ldots
\end{CD}$$
\end{definition}
\begin{lemma}\label{coker}
If $\X$ is a special formal $R$-scheme, and $i>0$ is an integer,
then there exists a canonical exact sequence of
$\mathcal{O}_{\X_\eta}$-modules
$$\begin{CD}
(\Omega^{i-1}_{\X/k})_{rig}@>d\pi\wedge >>
(\Omega^i_{\X/k})_{rig}@>>> \Omega^i_{\X_\eta/K}@>>> 0
\end{CD}$$
\end{lemma}
\begin{proof}
Since the functor $(.)_{rig}$ is exact by Proposition
\ref{rigmodule}, we get a canonical exact sequence of
$\mathcal{O}_{\X_\eta}$-modules
$$\begin{CD}
(\Omega^{i-1}_{\X/k})_{rig}@>d\pi\wedge >>
(\Omega^i_{\X/k})_{rig}@>>> (\Omega^i_{\X/k}/d\pi\wedge
\Omega^{i-1}_{\X/k})_{rig}@>>> 0
\end{CD}$$
But $\Omega^i_{\X/k}/d\pi\wedge \Omega^{i-1}_{\X/k}\cong
\Omega^i_{\X/R}$, so we can conclude by \cite[7.1.12]{dj-formal}.
  \end{proof}
\begin{lemma}\label{algebraic}
Let $X$ be a variety over $k$, and consider a morphism
$f:X\rightarrow \A^1_k=\mathrm{Spec}\,k[\pi]$ which is smooth of
pure relative dimension $m$ over the torus
$\mathbb{G}_m=\mathrm{Spec}\,k[\pi,\pi^{-1}]$. If we denote by
$\X$ the $\pi$-adic completion of $f$, then the Koszul complex of
$\X$ is exact, and
$$\begin{CD}i:\Omega^{m}_{\X_\eta/K}@>d\pi\wedge >>(\Omega^{m+1}_{\X/k})_{rig}\end{CD}$$ is an isomorphism.
\end{lemma}
\begin{proof}
Put $X'=X\times_{\A^1_k}\mathbb{G}_m$. Since $f$ is smooth over
$\mathbb{G}_m$,
$$\begin{CD}0@>>>\mathcal{O}_{X'}@>d\pi\wedge >>\Omega^{1}_{X'/k}@>d\pi\wedge >>\Omega^{2}_{X'/k}@>d\pi\wedge >>\ldots
\end{CD}$$ is exact, and hence the cokernels of the inclusion maps
$$d\pi\wedge \Omega^{i-1}_{X/k}\rightarrow \mathrm{ker}\,(d\pi\wedge :\Omega^i_{X/k}\rightarrow \Omega^{i+1}_{X/k})$$ are
$\pi$-torsion modules. Taking $\pi$-adic completions and using
\cite[1.9]{formal1}, we see that the cokernels of the maps
$$d\pi\wedge \Omega^{i-1}_{\X/k}\rightarrow \mathrm{ker}\,(d\pi\wedge :\Omega^i_{\X/k}\rightarrow \Omega^{i+1}_{\X/k})$$ are
$\pi$-torsion modules, so they vanish by passing to the generic
fiber. We can conclude by exactness of $(.)_{rig}$ (Proposition
\ref{rigmodule}) that the Koszul complex of $\X$ is exact.
Moreover, $\Omega^{m+2}_{\X/k}$ is the $\pi$-adic completion of
$\Omega^{m+2}_{X/k}$; hence, it is $\pi$-torsion, and
$(\Omega^{m+2}_{\X/k})_{rig}=0$. By Lemma \ref{coker}, this
implies that
$$\begin{CD}i:\Omega^{m}_{\X_\eta/K}@>d\pi\wedge >>(\Omega^{m+1}_{\X/k})_{rig}\end{CD}$$ is an isomorphism.
  \end{proof}
%
\begin{lemma}\label{rig}
Let $h:\X\rightarrow \mY$ be a morphism of special formal
$R$-schemes, such that $h_\eta:\X_\eta\rightarrow \mY_\eta$ is
\'etale. If the Koszul complex of $\X$ is exact, then the natural
map
$$\varphi:h_\eta^*((\Omega^i_{\mY/k})_{rig})=(h^*\Omega^{i}_{\mY/k})_{rig}\rightarrow
(\Omega^{i}_{\X/k})_{rig}$$ is an isomorphism of coherent
$\mathcal{O}_{\X_\eta}$-modules, for each $i\geq 0$. If, moreover,
$h_\eta$ is surjective, then the Koszul complex of $\mY$ is exact.
\end{lemma}
\begin{proof}
%
We proceed by induction on $i$. For $i=0$, the statement is clear,
so assume $i>0$. We put $\Omega^{-1}_{\X/k}=0$ and
$\Omega^{-1}_{\mY/k}=0$.
Now consider the commutative diagram
$$\begin{CD}
(h^*\Omega^{i-2}_{\mY/k})_{rig}@>d\pi\wedge >>
(h^*\Omega^{i-1}_{\mY/k})_{rig}@>d\pi\wedge >>
(h^*\Omega^{i}_{\mY/k})_{rig}@>>>h^*_\eta(\Omega^i_{\mY/K})@>>> 0
\\ @| @| @VVV @|
\\ (\Omega^{i-2}_{\X/k})_{rig}@>d\pi\wedge >> (\Omega^{i-1}_{\X/k})_{rig}
@>d\pi\wedge >> (\Omega^{i}_{\X/k})_{rig}@>>>
\Omega^i_{\X_\eta/K}@>>> 0
\end{CD}$$
The bottom row is exact by exactness of the Koszul complex on $\X$
and by Lemma \ref{coker}. The upper row is exact except maybe at
$(h^*\Omega^{i-1}_{\mY/k})_{rig}$, by Lemma \ref{coker} (applied
to $\mY$) and flatness of $h_\eta$. The first and second vertical
arrows are isomorphisms by the induction hypothesis, and the
fourth one is an isomorphism since $h_\eta$ is \'etale
\cite[2.6]{formrigIII}. Now a diagram chase shows that the third
vertical arrow is an isomorphism as well. If $h_\eta$ is also
surjective, then by faithful flatness the Koszul complex of $\mY$
is exact since the complex
$$\begin{CD}
0@>>> h_\eta^*\mathcal{O}_{\mY_\eta}@>d\pi\wedge >>
h_\eta^*(\Omega^1_{\mY/k})_{rig}@>d\pi\wedge
>>h_\eta^*(\Omega^2_{\mY/k})_{rig}@>d\pi\wedge >>\ldots
\end{CD}$$ is isomorphic to the Koszul complex of $\X$ and hence
exact.
  \end{proof}

\begin{prop}\label{gelfand}
If $\X$ is a generically smooth special formal $R$-scheme of pure
relative dimension $m$, then the Koszul complex of $\X$ is exact,
and
$$\begin{CD}i:\Omega^{m}_{\X_\eta/K}@>d\pi\wedge >>(\Omega^{m+1}_{\X/k})_{rig}\end{CD}$$ is an isomorphism.
\end{prop}
\begin{proof}
We may assume that $\X$ is affine, say $\X=\mathrm{Spf}\,A$. We
use the notation of Section \ref{generic}.  The morphism of
special formal $R$-schemes
$$\mY:=\mathrm{Spf}\,B_n\rightarrow \X$$ induces an open embedding on
the generic fibers. By Lemma \ref{coker} and Lemma \ref{rig}, it
suffices to show that the Koszul complex of $\mY$ is exact and
that $(\Omega^{m+2}_{\mY/k})_{rig}=0$. Hence, we may as well
assume that $A$ is topologically of finite type over $R$. By
resolution of singularities (Proposition \ref{affineres}) and the
proof of Proposition \ref{rev}, and again applying Lemma
\ref{rig}, we may assume that $\X=\mathrm{Spf}\,A$ is endowed with
an \'etale morphism of formal $R$-schemes
$$f:\mZ\rightarrow
\mathrm{Spf}\,R\{x_0,\ldots,x_m\}/(\pi-\prod_{i=0}^{m}x_i^{N_i})$$
with $N_i\in \N$.
%
Since $f$ is \'etale, it is enough to prove the result for
$$\X=\mathrm{Spf}\,R\{x_0,\ldots,x_m\}/(\pi-\prod_{i=0}^{m}x_i^{N_i})$$
Now we can conclude by Lemma \ref{algebraic}.   \end{proof}
\begin{cor}\label{bd-iso}
If, moreover, $\X$ is affine, then the natural map
$$\begin{CD}\Omega^m_{\X/R}(\X)\otimes_R K@>d\pi\wedge >>
\Omega^{m+1}_{\X/k}(\X)\otimes_R K\end{CD}$$ is an isomorphism,
and it fits in a commutative diagram
$$\begin{CD}
\Omega^m_{\X/R}(\X)\otimes_R K@>d\pi\wedge
>>\Omega^{m+1}_{\X/k}(\X)\otimes_R K
\\@VVV @VVV
\\\Omega^{m}_{\X_\eta/K}(\X_\eta)@>d\pi\wedge >>(\Omega^{m+1}_{\X/k})_{rig}(\X_\eta)
\end{CD}$$
where the vertical arrows are injections and the horizontal arrows
are isomorphisms.
\end{cor}
\begin{proof}
This follows immediately from Lemma \ref{bounded-iso}.
\end{proof}
\begin{definition}[Gelfand-Leray form]\label{def-gelfand}
If $\X$ is a generically smooth special formal $R$-scheme of pure
relative dimension $m$, and if $\omega$ is an element of
$\Omega^{m+1}_{\X/k}(\X)$, then we denote by $\omega/d\pi$ the
inverse image of $\omega$ under the isomorphism
$$\begin{CD}i:\Omega^{m}_{\X_\eta/K}@>d\pi\wedge >>(\Omega^{m+1}_{\X/k})_{rig}\end{CD}$$ and we call it the
Gelfand-Leray form associated to $\omega$.
\end{definition}
\begin{remark}
Let us compare this definition with the construction made in
\cite[\S 9.2]{NiSe}. Let $X$ be a smooth irreducible variety over
$k$, of dimension $m+1$, and let $f:X\rightarrow
\A^1_k=\mathrm{Spec}\,k[t]$ be a morphism of $k$-varieties, smooth
over the torus $\mathrm{Spec}\,k[t,t^{-1}]$. Let $\omega$ be a
gauge form on $X$, and denote by $V$ the complement in $X$ of the
special fiber $X_s$ of $f$. In \cite[\S 9.2]{NiSe}, we constructed
a relative form $\omega/d\pi$ in $\Omega^m_{V/\A^1_k}(V)$ and we
defined the Gelfand-Leray form as the element of
$\Omega^m_{\X_\eta/K}(\X_\eta)$ induced by $\omega/d\pi$. It is
obvious from the constructions that this form coincides with our
Gelfand-Leray form associated to the element of
$\Omega^{m+1}_{\X/k}(\X)$ obtained from $\omega$ by completion.
  \end{remark}
\begin{cor}\label{gaugeform}
If $\X$ is a regular flat special formal $R$-scheme of pure
relative dimension $m$, and if $\omega$ is an element of
$\Omega^{m+1}_{\X/k}(\X)$, then $\omega/d\pi$ is $\X$-bounded. If,
moreover, $\omega$ is a gauge form on $\X$ (i.e. a nowhere
vanishing section of $\Omega^{m+1}_{\X/k}(\X)$), then
$\omega/d\pi$ is a bounded gauge form on $\X_\eta$.
\end{cor}
\begin{proof}
Boundedness follows from Corollary \ref{bd-iso}.  Now suppose that
$\omega$ is a gauge form on $\X$. We may assume that
$\X=\mathrm{Spf}\,A$ is affine. The fact that $\omega$ is gauge
means that $\omega\notin \mathfrak{M}\Omega^{m+1}_{\X/k}(\X)$ for
each prime ideal $\mathfrak{M}$ of $A$; using
\cite[7.1.9]{dj-formal}, we see that this implies that de image of
$\omega$ in $(\Omega^{m+1}_{\X/k})_{rig}(\X_\eta)$ does not vanish
at any point $x$ of $\X_\eta$. Hence, since the map $i$ of
Proposition \ref{gelfand} is an isomorphism of coherent
$\mathcal{O}_{\X_\eta}$-modules, $\omega/d\pi$ is gauge.
\end{proof}

\begin{lemma}\label{gelfand-pullback}
If $h:\mY\rightarrow \X$ is a morphism of generically smooth
special formal $R$-schemes, both of pure relative dimension $m$,
and if $\omega$ is a global section of $\Omega^{m+1}_{\X/k}$, then
$$(h^*\omega)/d\pi=(h_\eta)^*(\omega/d\pi)$$ in
$\Omega^m_{\mY_\eta/K}(\mY_\eta)$.
\end{lemma}
\begin{proof}
It suffices to show that
$$d\pi\wedge((h_\eta)^*\alpha)
=(h_\eta)^*(d\pi\wedge\alpha)$$ for any $m$-form $\alpha$ on
$\X_\eta$; substituting $\alpha$ by $\omega/d\pi$ yields the
result. For any $i\geq 0$, the square
$$\begin{CD}
h^*\Omega^{i}_{\X/k}@>>> \Omega^{i}_{\mY/k}
\\ @Vd\pi\wedge  VV @VVd\pi\wedge  V
\\h^*\Omega^{i+1}_{\X/k}@>>> \Omega^{i+1}_{\mY/k}
\end{CD}$$
commutes, and therefore
$$\begin{CD}
h^*\Omega^{m}_{\X/R}@>>> \Omega^{m}_{\mY/R}
\\ @Vd\pi\wedge  VV @VVd\pi\wedge  V
\\h^*\Omega^{m+1}_{\X/k}@>>> \Omega^{m+1}_{\mY/k}
\end{CD}$$
commutes. We can conclude by passing to the generic fiber.
\end{proof}

\begin{lemma}\label{perfect}
If $\X$ is a separated formal scheme of pseudo-finite type over
$\mathrm{Spec}\,F$, with $F$ a perfect field, and $\X$ is regular,
then $\X$ is smooth over $F$.
\end{lemma}
\begin{proof}
Let $x$ be a closed point of $\X_0$, and let $(x_0,\ldots,x_m)$ be
a regular system of local parameters on $\X$ at $x$. These define
a morphism of formal schemes of pseudo-finite type over
$\mathrm{Spec}\,F$
$$h:\mU\rightarrow \A^{m+1}_F$$
on some open neighbourhood $\mU$ of $x$ in $\X$. Since
$\A^{m+1}_F$ is smooth over $F$, it suffices to show that $h$ is
\'etale at $x$. This follows immediately from Lemma
\ref{etalecrit} and the fact that $F$ is perfect.   \end{proof}
\begin{lemma}\label{globaldiff}
Let $\X$ be a flat regular special formal $R$-scheme, of pure
relative dimension $m$. Then $\X$ is smooth over $k$, of pure
dimension $m+1$, and $\Omega^{m+1}_{\X/k}$ is a locally free sheaf
of rank $1$ on $\X$.
\end{lemma}
\begin{proof}
The fact that $\X$ is smooth over $k$ follows from Lemma
\ref{perfect}, since $k$ has characteristic zero. The fact that it
has pure dimension $m+1$ follows from the flatness of $\X$ over
$R$.

By \cite[4.8]{formal1}, the sheaf of continuous differential forms
$\Omega^1_{\X/k}$ is locally free; by \cite[5.10]{formal2}, it has
rank $m+1$.   \end{proof}

\begin{cor}\label{motvol-welldef2}
If $\X$ is a regular special formal $R$-scheme, then we can cover
$\X$ by open formal subschemes $\mU$ such that $\mU_\eta$ admits a
$\mU$-bounded gauge form.
\end{cor}
\begin{proof}
By Lemma \ref{globaldiff}, we can cover $\X$ by open formal
subschemes $\mU$ such that $\Omega^{m+1}_{\X/k}\cong
\mathcal{O}_{\mU}$. By Corollary \ref{gaugeform}, each $\mU_\eta$
admits a $\mU$-bounded gauge form.
  \end{proof}


If $h:\X'\rightarrow \X$ is a morphism of smooth formal
$k$-schemes of pseudo-finite type of pure dimension $m+1$, we
define the Jacobian ideal sheaf of $h$ as the annihilator of the
cokernel of the natural map of locally free rank one
$\mathcal{O}_{\X'}$-modules
$$\psi:h^{*}\Omega^{m+1}_{\X/k}\rightarrow \Omega^{m+1}_{\X'/k}$$
and we denote this ideal sheaf by $\mathcal{J}ac_{h/k}$ to
distinguish it from the Jacobian ideal sheaf $\mathcal{J}ac_h$
from Section \ref{ordersect}.

\begin{lemma}\label{jachk}
Let $h:\X'\rightarrow \X$ be a morphism of regular flat special
formal $R$-schemes, both of pure relative dimension $m$, such that
$h_\eta$ is \'etale. Then
 the
Jacobian ideal sheaf $\mathcal{J}ac_{h/k}$ is invertible, and
contains a power of $\pi$.
\end{lemma}
\begin{proof}
By Lemma \ref{rig} and exactness of the functor $(.)_{rig}$, we
see that $coker(\psi)_{rig}=0$. This means that
$\mathcal{J}ac_{h/k}$ contains a power of $\pi$, by Corollary
\ref{tor}.

Both $h^*\Omega^{m+1}_{\X/k}$ and $\Omega^{m+1}_{\X'/k}$ are line
bundles on $\X'$, by Lemma \ref{globaldiff}. Covering $\X'$ by
sufficiently small affine open formal subschemes
$\mU=\mathrm{Spf}\,A$, we may assume that they are trivial; let
$\omega$ and $\omega'$ be generators for $h^*\Omega^{m+1}_{\X/k}$
resp. $\Omega^{m+1}_{\X'/k}$. Then we can write
$\psi(\omega)=f\omega'$ on $\mU$ with $f$ in $A$, and $f$
generates the ideal sheaf $\mathcal{J}ac_{h/k}$ on $\mU$.
%
%
  \end{proof}

Let $\X$ be a regular special formal $R$-scheme, whose special
fiber is a strict normal crossings divisor $\sum_{i\in
I}N_i\mE_i$. Let $\mathcal{J}$ be an invertible ideal sheaf on
$\X$,  and fix $i\in I$. One can show as in Lemma \ref{cons-lemma}
that the length of the $\mathcal{O}_{\X,\mE_i,x}$-module
$\mathcal{O}_{\X,\mE_i,x}/\mathcal{J}\mathcal{O}_{\X,\mE_i,x}$ is
independent of the point $x$ of $E_i$. We call this value the
multiplicity of $\mathcal{J}$ along $\mE_i$.
\begin{definition}
Let $h:\X'\rightarrow \X$ be a morphism of regular flat special
formal $R$-schemes, both of pure relative dimension $m$, such that
$h_\eta$ is \'etale. If $\X'_s$ is a strict normal crossings
divisor $\sum_{i\in I}N_i\mE_i$, then we denote by $\nu_i-1$ the
multiplicity of $\mathcal{J}ac_{h/k}$ along $\mE_i$, and we write
$K_{\X'/\X}=\sum_{i\in I}(\nu_i-1)\mE_i$
\end{definition}
\begin{lemma}\label{ordgelf}
Let $h:\X'\rightarrow \X$ be a morphism of regular flat special
formal $R$-schemes, both of pure relative dimension $m$, such that
$h_\eta$ is \'etale. If $\X'_s$ is a strict normal crossings
divisor $\sum_{i\in I}N_i\mE_i$, and if $$K_{\X'/\X}=\sum_{i\in
I}(\nu_i-1)\mE_i$$ then for any gauge form $\omega$ on $\X$ and
any $i\in I$,
$$ord_{\mE_i}(h^*_\eta(\omega/d\pi))=\nu_i-N_i$$
\end{lemma}
\begin{proof}
First of all, note that $h^*_\eta(\omega/d\pi)=(h^*\omega)/d\pi$
by Lemma \ref{gelfand-pullback}. Choose an index $i$ in $I$ and a
point $x'$ on $E_i$, and put $x=h(x')$. Shrinking $\X$ to an open
formal neighbourhood of $x$, we may suppose that there exists an
integer $a$ such that $\phi:=\pi^a(\omega/d\pi)$ belongs to
$$\Omega^m_{\X/R}/(\pi-\mathrm{torsion})\subset \Omega^m_{\X_\eta/K}(\X_\eta)$$
since $\omega/d\pi$ is $\X$-bounded by Lemma \ref{gaugeform}.

Consider the commutative diagram
$$\begin{CD}
h^*\Omega^m_{\X/R}@>d\pi\wedge >> h^*\Omega^{m+1}_{\X/k}
\\@VVV @VVV
\\ \Omega^m_{\X'/R}@>d\pi\wedge >>
\Omega^{m+1}_{\X'/k}\end{CD}$$ Since $\Omega^{m+1}_{\X/k}$ and
$\Omega^{m+1}_{\X'/k}$ are locally free, we get a commutative
diagram (using the notation in Section \ref{ordstrict})
$$\begin{CD}
h^*\Omega^m_{\X/R}\otimes \mathcal{O}_{\X',\mE_i,x'}@>d\pi\wedge
>> h^*\Omega^{m+1}_{\X/k}\otimes\mathcal{O}_{\X',\mE_i,x'}
\\@V\varphi VV @VV\psi V
\\ \Omega_{\X',\mE_i,x'}@>d\pi\wedge >>
\Omega^{m+1}_{\X'/k}\otimes \mathcal{O}_{\X',\mE_i,x'}\end{CD}$$

By definition, $$ord_{\mE_i}(h_\eta^*(\omega/d\pi))=
\mathrm{length}\left(\Omega_{\X',\mE_i,x'}/(\varphi(\phi)\cdot\mathcal{O}_{\X',\mE_i,x'})\right)
-aN_i$$ Since $d\pi\wedge\phi =\pi^a\omega$ and $\omega$ generates
$\Omega^{m+1}_{\X/k}$, we see that the
$\mathcal{O}_{\X',\mE_i,x'}$-module
$$(\Omega^{m+1}_{\X'/k}\otimes
\mathcal{O}_{\X',\mE_i,x'})/(\psi(d\pi\wedge
\phi)\mathcal{O}_{\X',\mE_i,x'})$$ has length $\nu_i-1+aN_i$. But
$\psi(d\pi\wedge \phi)=\varphi(d\pi\wedge \phi)$, so it suffices
to show that the cokernel of the lower horizontal map
$$\begin{CD}\Omega_{\X',\mE_i,x'}@>d\pi\wedge >> \Omega^{m+1}_{\X'/k}\otimes
\mathcal{O}_{\X',\mE_i,x'}\end{CD}$$ has length $N_i-1$. Since
this value does not change if we pass to an \'etale cover of $\X'$
whose image contains $x'$ (by the algebraic argument used in the
proof of Lemma \ref{cons-lemma}), we may assume that
$\pi=\prod_{j=0}^{m}x_j^{N_j}$ on $\X'$, with $(x_0,\ldots,x_m)$ a
regular sequence. Hence, taking differentials, we see that
$$\omega_0:=dx_0\wedge\ldots\wedge\widehat{dx_i}\wedge\ldots\wedge dx_m$$
generates $\Omega_{\X',\mE_i,x'}$, and it is clear that
$ord_{x_i}(d\pi\wedge \omega_0)=N_i-1$ in
$\Omega^{m+1}_{\X'/k}\otimes \mathcal{O}_{\X',\mE_i,x'}$ (see
Definition \ref{def-ord} for the notation $ord_{x_i}(.)$\,).
  \end{proof}
\begin{prop}\label{mellin}
Let $\X$ be a regular flat special formal $R$-scheme of pure
relative dimension $m$, and let $\omega$ be a gauge form in
$\Omega^{m+1}_{\X/k}(\X)$. The volume Poincar\'e series
$S(\X,\omega/d\pi;T)$ only depends on $\X$, and not on $\omega$.
In fact, if $\X'\rightarrow \X$ is any resolution of
singularities, with $\X'_s=\sum_{i\in I}N_i\mE_i$ and
$K_{\X'/\X}=\sum_{i\in I}(\nu_i-1)\mE_i$, then
$S(\X,\omega/d\pi;T)$ is given explicitly by
$$S(\X,\omega/d\pi;T)=\LL^{-m}\sum_{\emptyset\neq J\subset
I}(\LL-1)^{|J|-1}[\widetilde{E}_J^{o}] \prod_{i\in
J}\frac{\LL^{N_i-\nu_i}T^{N_i}}{1-\LL^{N_i-\nu_i}T^{N_i}} \
\mathrm{in}\ \mathcal{M}_{\X_0}[[T]]$$
\end{prop}
\begin{proof}
By additivity of the motivic integral, we may assume that $\X$ is
affine. Then $\X$ admits a resolution of singularities by
Proposition \ref{affineres}, and the expression for
$S(\X,\omega/d\pi;T)$ follows from Corollary \ref{serrat} and
Lemma \ref{ordgelf}. This expression is clearly independent of
$\omega$.   \end{proof}
\begin{remark}
The fact that $S(\X,\omega/d\pi;T)$ does not depend on $\omega$
follows already from the fact that $\omega/d\pi$ is independent of
$\omega$ up to multiplication with a unit on $\X$. Hence, for each
$d>0$, $\int_{\X(d)}|(\omega/d\pi)(d)|$ does not depend on
$\omega$. Beware that $S(\X,\omega/d\pi;T)$ depends on the choice
of $\pi$ if $k$ is not algebraically closed; see the remark
following Definition \ref{def-volpoin}.   \end{remark}
\begin{definition}\label{localsing}
If $\X$ is a regular special formal $R$-scheme, and $\X$ admits a
gauge form $\omega$, then we define the local singular series of
$\X$ by
$$F(\X;*):\N^{\ast}\rightarrow \mathcal{M}_{\X_0}:d\mapsto
\int_{\X(d)}|\frac{\omega}{d\pi}(d)|$$ This definition only
depends on $\X$, and not on $\omega$, by Proposition \ref{mellin}.

If $\X$ is any regular special formal $R$-scheme, we choose a
finite cover $\{\mU_i\}_{i\in I}$ of $\X$ by open formal
subschemes such that each $\mU_i$ admits a gauge form $\omega_i$,
and we define the local singular series of $\X$ by
$$F(\X;*):\N^{\ast}\rightarrow \mathcal{M}_{\X_0}:d\mapsto
\sum_{\emptyset\neq J\subset I}(-1)^{|J|+1}F(\cap_{j\in
J}\mU_j;d)$$ This definition does not depend on the chosen cover,
by additivity of the motivic integral.
 We define the
Weil generating series of $\X$ by
$$S(\X;T)=\sum_{d>0}F(\X;d)T^d\in\mathcal{M}_{\X_0}[[T]]$$
\end{definition}
In the terminology of \cite[4.4]{ClLo}, the Weil generating series
is the Mellin transform of the local singular series. If $\X$
admits a gauge form $\omega$, then by definition,
$S(\X;T)=S(\X,\omega/d\pi;T)$.

The term ``Weil generating series'' is justified by the fact that
$F(\X;d)$ can be seen as a measure for the set
$\cup_{K'}\X_\eta(K')$ where $K'$ varies over the unramified
extensions of $K(d)$. Moreover, we have the following immediate
corollary of the trace formula in Theorem \ref{trace}.
\begin{prop}
Let $\varphi$ be a topological generator of $G(K^s/K^{sh})$. If
$\X$ is a regular special formal $R$-scheme, then for any integer
$d>0$,
$$\chi_{top}(F(\X;d))=Tr(\varphi^d\,|\,H(\overline{\X_\eta}))$$
\end{prop}

\subsection{The motivic volume}
Let $\X$ be a generically smooth, special formal $R$-scheme, of
pure relative dimension $m$, and let $\omega$ be a $\X$-bounded
gauge form on $\X_\eta$.
 It is not
possible to associate a motivic Serre invariant to
$\X\widehat{\times}_R\widehat{R^s}$ in a direct way, since the
normalization $R^s$ of $R$ in $K^s$ is not a discrete valuation
ring. We will define a motivic object by taking a limit of motivic
integrals over finite ramifications of $\X$, instead.

\begin{definition}[\cite{GLM}, (2.8)]
There is a unique $\mathcal{M}_{\X_0}$-linear morphism
$$\lim_{T\to\infty}:\mathcal{M}_{\X_0}\left[\frac{\LL^{a}T^b}{1-\LL^{a}T^b} \right]_{(a,b)\in \Z\times \N^{\ast}}
\longrightarrow \mathcal{M}_{\X_0}$$ mapping $$\prod_{(a,b)\in
I}\frac{\LL^{a}T^b}{1-\LL^{a}T^b}$$ to
$(-1)^{|I|}=(-1)^{|I|}[\X_0]$, for each finite subset $I$ of
$\Z\times \N^{\ast}$. We call the image of an element its limit
for $T\to\infty$.
\end{definition}
\begin{prop}\label{separable}
Let $\X$ be a generically smooth, special formal $R$-scheme, of
pure relative dimension $m$, and let $\omega$ be a $\X$-bounded
gauge form on $\X_\eta$. The limit of $-S(\X,\omega;T)$ for $T\to
\infty$ is well-defined, and does not depend on $\omega$. If
$\X'\rightarrow \X$ is any resolution of singularities, with
$\X'_s=\sum_{i\in I}N_i \mE_i$, then this limit is given
explicitly by
$$\LL^{-m}\sum_{\emptyset\neq J\subset
I}(1-\LL)^{|J|-1}[\widetilde{E}_J^{o}]$$ in $\mathcal{M}_{\X_0}$.
\end{prop}
\begin{proof}
This follows immediately from the computation in Corollary
\ref{serrat}, and the observation preceding Corollary
\ref{rationality}.   \end{proof}

\begin{definition}\label{motvolume}
Let $\X$ be a generically smooth special formal $R$-scheme of pure
relative dimension, and assume that $\X_\eta$ admits a
$\X$-bounded gauge form.
 The
motivic volume
$$S(\X;\widehat{K^s})\in
\mathcal{M}_{\X_0}$$ is by definition the limit of
$-S(\X,\omega;T)$ for $T\to\infty$, where $\omega$ is any
$\X$-bounded gauge form on $\X_\eta$.

%
\end{definition}
%
If $h:\mY\rightarrow \X$ is a morphism of
generically smooth special formal $R$-schemes such that $h_\eta$
is an isomorphism, and if $\X_\eta$ admits a $\X$-bounded gauge
form, then it is clear from the definition that
$S(\X;\widehat{K^s})=S(\mY;\widehat{K^s})$ in
$\mathcal{M}_{\X_0}$.

In definition \ref{motvolume}, the condition that $\X_\eta$ admits
a gauge form can be avoided as follows.

\begin{definition-prop}\label{motvol-resolution}
If $\X$ is a generically smooth special formal $R$-scheme which
admits a resolution of singularities, then there exists a morphism
of special formal $R$-schemes $h:\mY\rightarrow \X$ such that
$h_\eta$ is an isomorphism, and such that $\mY$ has a finite open
cover $\{\mU_i\}_{i\in I}$ such that $\mU_i$ has pure relative
dimension and $(\mU_i)_\eta$ admits a $\mU_i$-bounded gauge form
for each $i$. Moreover, the value
$$S(\X;\widehat{K^s})=\sum_{\emptyset\neq J\subset
I}(-1)^{|J|+1}S(\cap_{i\in J}\mU_i;\widehat{K^s})\in
\mathcal{M}_{\X_0}$$ only depends on $\X$.
\end{definition-prop}
\begin{proof}
Since $\X$ admits a resolution of singularities,  we may assume
that $\X$ is regular and flat; now it suffices to put $\mY=\X$ and
to apply Lemma \ref{motvol-welldef2}.

The fact that the expression for $S(\X;\widehat{K^s})$ only
depends on $\X$, follows from the additivity of the motivic
integral, and the fact that we can dominate any two such morphisms
$h$ by a third by taking the fibered product.   \end{proof}
\begin{definition}\label{motvol2}
Let $\X$ be a generically smooth special formal $R$-scheme, and
take a finite cover $\{\mU_i\}_{i\in I}$ of $\X$ by affine open
formal subschemes. Then we can define the motivic volume
$S(\X;\widehat{K^s})$ by
$$S(\X;\widehat{K^s})=\sum_{\emptyset\neq J\subset
I}(-1)^{|J|+1}S(\cap_{i\in J}\mU_i;\widehat{K^s})$$ This
definition only depends on $\X$.
%
\end{definition}

Note that the terms $S(\cap_{i\in J}\mU_i;\widehat{K^s})$ are
well-defined, since each $\mU_i$ admits a resolution of
singularities by Proposition \ref{affineres} and hence,
Proposition-Definition \ref{motvol-resolution} applies.

\begin{remark}
Beware that the motivic volume $S(\X;\widehat{K^s})$ depends on
the choice of uniformizer $\pi$ in $R$ (more precisely, on the
fields $K(d)$), if $k$ is not algebraically closed. For instance,
if $k=\Q$ and $\X=\Spf R[x]/(x^2-2\pi)$ then
$S(\X;\widehat{K^s})=[\Spec \Q(\sqrt{2})]$ in $\mathcal{M}_k$,
while for $\X=\Spf R[x]/(x^2-\pi)$ we find $S(\X;\widehat{K^s})=2$
(to see that these are distinct elements of $\mathcal{M}_k$, look
at their \'etale realizations in the Grothendieck ring of
$\ell$-adic representations of $G(\overline{\Q}/\Q)$).

If $k$ is algebraically closed, the motivic volume is independent
of the choice of uniformizer, since $K(d)$ is the unique extension
of degree $d$ of $K$ in $K^s$; see the remark following Definition
\ref{def-volpoin}.

It is not hard to see that for any unramified extension $R'$ of
$R$, and any generically smooth special formal $R$-scheme $\X$,
the motivic volume of $\X'=\X\widehat{\times}_R R'$ is the image
of the motivic volume of $\X$ under the natural base change
morphism $\mathcal{M}_{\X_0}\rightarrow \mathcal{M}_{\X'_0}$.
\end{remark}

Now we define the motivic volume of a smooth rigid $K$-variety
$\X_\eta$ that can be realized as the generic fiber of a special
formal $R$-scheme $\X$. If $\X_\eta$ is quasi-compact, this
definition was given in \cite[8.3]{NiSe}: the image of
$S(\X;\widehat{K^s})$ under the forgetful morphism
$\mathcal{M}_{\X_0}\rightarrow \mathcal{M}_k$ only depends on
$\X_\eta$, and it was called the motivic volume
$S(\X_\eta;\widehat{K^s})$ of $\X_\eta$. I do not know if this
still holds if $\X_\eta$ is not quasi-compact; the problem is that
it is not clear if any two formal $R$-models of $\X_\eta$ can be
dominated by a third. Therefore, we need an additional technical
condition (which might be superfluous).

\begin{definition}
Let $\X$ be a special formal $R$-scheme, and suppose that
$\X_\eta$ is reduced. For any $i\geq 0$, a section of
$\Omega^i_{\X_\eta}(\X_\eta)$ is called a universally bounded
$i$-form on $\X_\eta$, if it is $\mY$-bounded for each formal
$R$-model $\mY$ of $\X_\eta$.
\end{definition}

I don't know an example of a bounded $i$-form which is not
universally bounded. If $\X_\eta$ is reduced, then an analytic
function on $\X_\eta$ is bounded iff it is universally bounded, by
Lemma \ref{bounded-function}. If $\X_\eta$ is quasi-compact, then
any differential form on $\X_\eta$ is universally bounded.

\begin{lemma}\label{universally}
If $\X$ is an affine special formal $R$-scheme, and $\X_\eta$ is
reduced, then any $\X$-bounded $i$-form $\omega$ on $\X_\eta$ is
universally bounded.
\end{lemma}
\begin{proof}
 Since it suffices
to prove that $\pi^a\omega$ is universally bounded, for some
integer $a$, we may suppose that $\omega$ belongs to the image of
the natural map
$$\Omega^i_{\X/R}(\X)\rightarrow \Omega^i_{\X_\eta/K}(\X_\eta)$$
by Lemma \ref{bounded-affine}. This means that we can write
$\omega$ as a sum of terms of the form $a_0(da_1\wedge\dots\wedge
da_i)$ with $a_0,\ldots,a_i$ regular functions on $\X$, and hence
$a_0,\ldots,a_i$ are bounded by $1$ on $\X_\eta$. If $\mY$ is any
formal $R$-model for $\X_\eta$, then by Lemma
\ref{bounded-function}, the functions $a_0,\ldots,a_i$ on
$\X_\eta$ are $\mY$-bounded; so $\omega$ is $\mY$-bounded.
\end{proof}

\begin{definition-prop}
Let $\X$ be any generically smooth special formal $R$-scheme, and
assume that $\X_\eta$ admits a universally bounded gauge form
$\omega$. The image of $S(\X;\widehat{K^s})$ under the forgetful
morphism $\mathcal{M}_{\X_0}\rightarrow \mathcal{M}_k$ only
depends on $\X_\eta$; we call it the motivic volume of $\X_\eta$,
and denote it by $S(\X_\eta;\widehat{K^s})$.
\end{definition-prop}
\begin{proof}
If $\X$ is any formal $R$-model for $\X_\eta$, then $\omega$ is
$\X$-bounded, and the image of $S(\X;\widehat{K^s})$ in
$\mathcal{M}_k$ coincides with
$$-\lim_{T\to \infty}\sum_{d>0}\left(\int_{\X(d)_\eta}|\omega(d)|\right)T^d$$ by
Proposition \ref{specialize}. Hence, it does not depend on $\X$
(and neither on $\omega$).   \end{proof}

In particular, this definition applies to the generic fiber of an
affine regular special formal $R$-scheme $\X$ that admits a gauge
form $\omega\in \Omega^{\mathrm{max}}_{\X/k}(\X)$, by Corollary
\ref{gaugeform} and Lemma \ref{universally}. Hence, for any
generically smooth special formal $R$-scheme $\X$, we can cover
$\X_\eta$ by a finite number of open rigid subvarieties
$U_i,\,i\in I$ such that $S(\cap_{i\in J}U_i;\widehat{K^s})$ is
defined for each non-empty subset $J$ of $I$, by Proposition
\ref{affineres}. However, it is not clear if the value
$\sum_{\emptyset \neq J\subset I}(-1)^{|J|+1}S(\cap_{i\in
J}U_i;\widehat{K^s})$ is independent of the chosen cover: if
$V_\ell,\ell \in L$ is another such cover, I do not know if
$V_\ell\cap U_i$ admits a universally bounded gauge form for all
$i$ and $\ell$.
%
%

If $\X$ is $stft$, we recover the definitions from \cite{NiSe}.

\begin{prop}
Let $\X$ be a generically smooth special formal $R$-scheme, and
$V$ a locally closed subset of $\X_0$. If we denote by
$\mathfrak{V}$ the formal completion of $\X$ along $V$, then
$S(\mathfrak{V};\widehat{K^s})$ coincides with the image of
$S(\X;\widehat{K^s})$under the base-change morphism
$\mathcal{M}_{\X_0}\rightarrow \mathcal{M}_V$.
\end{prop}
\begin{proof}
This follows immediately from Proposition \ref{tube-volume}.
\end{proof} If $\X$ is $stft$ over $R$, then we called in
\cite{NiSe} the image of $S(\X;\widehat{K^s})$ in
$\mathcal{M}_{V}$ the motivic volume of $\X$ with support in $V$,
and we denoted it by $S_V(\X;\widehat{K^s})$. The above
proposition shows that
$$S_V(\X;\widehat{K^s})=S(\mathfrak{V};\widehat{K^s})$$ in
$\mathcal{M}_{V}$. In particular, it only depends on
$\mathfrak{V}$, and not on the embedding in $\X$.

\begin{prop}\label{eulintr}
If $\X$ is a generically smooth special formal $R$-scheme, then
$$\chi_{top}(S(\X;\widehat{K^s}))=\chi_{\acute{e}t}(\overline{\X_\eta})$$
where $\chi_{\acute{e}t}$ is the Euler characteristic associated
to Berkovich' \'etale $\ell$-adic cohomology for non-archimedean
analytic spaces.

In particular, if $\X_\eta$ admits a universally bounded gauge
form, then
$$\chi_{top}(S(\X_\eta;\widehat{K^s}))=\chi_{\acute{e}t}(\overline{\X_\eta})$$
\end{prop}
\begin{proof}
Let $\varphi$ be a topological generator of the geometric
monodromy group $G(K^s/K^{sh})$. By definition,
$$\chi_{top}(S(\X;\widehat{K^s}))=-\lim_{T\to\infty}\left(\sum_{d>0}\chi_{top}(S(\X(d)_\eta))T^d\right)$$
Hence, by our Trace Formula in Theorem \ref{trace},
$$\chi_{top}(S(\X;\widehat{K^s}))=-\lim_{T\to\infty}\left(\sum_{d>0}Tr(\varphi^d\,|\,H(\overline{\X_\eta})\,)T^d\right)$$

Recall the identity \cite[1.5.3]{weil1}
\begin{eqnarray*}\label{weil1}
\sum_{d>0}Tr(F^d,V)T^d&=&T\frac{d}{dT}log(det(1-TF,V)^{-1})
\\&=&-\frac{T\frac{d}{dT}(det(1-TF,V))}{det(1-TF,V)}
\end{eqnarray*} for any endomorphism $F$ on a finite dimensional vector
space $V$ over a field of characteristic zero. Taking limits, we
get
$$
-\lim_{T\to\infty}\sum_{d>0}Tr(F^d,V)T^d=dim(V)$$ Applying this to
$F=\varphi$ and $V=H(\overline{\X_\eta})$ yields the result.
\end{proof}

\section{The analytic Milnor fiber}\label{Milnor} In this section, we prove that the
analytic Milnor fiber introduced in \cite{NiSe} determines a
singularity up to formal equivalence. We do not impose any
restriction on the residue field $k$.

\subsection{Branches of formal schemes}

\begin{definition}
 Let $\X$ be a flat special formal $R$-scheme, and
let $x$ be a closed point of $\X_0$. Consider the normalization
$\widetilde{\X}\rightarrow \X$, and let $x_1,\ldots,x_m$ be the
points on $\widetilde{\X}_0$ lying over $x$. We call the special
formal $R$-schemes
$$\mathrm{Spf}\,\widehat{\mathcal{O}}_{\widetilde{\X},x_1},\ldots,
\mathrm{Spf}\,\widehat{\mathcal{O}}_{\widetilde{\X},x_m}$$ the
branches of $\X$ at $x$.
\end{definition}

\begin{remark}
This notion should not be confused with the branches of the
special fiber $\X_s$ at $x$. For instance, if
$\X=\mathrm{Spf}\,R\{x,y\}/(xy-\pi)$, then
$\X_s=\mathrm{Spec}\,k[x,y]/(xy)$ has two branches at the origin,
while $\X$ is normal.   \end{remark}

\begin{prop}\label{prop-branch}
Let $\X$ be a flat special formal $R$-scheme, and let $x$ be a
closed point on $\X_0$. Suppose that $]x[$ is normal, and let
$x_1,\ldots,x_m$ be the points lying over $x$ in the normalization
$h:\widetilde{\X}\rightarrow \X$. There is a canonical isomorphism
$$\mathcal{O}_{]x[}(\,]x[\,)=\prod_{i=1}^m {\mathcal{O}}_{]x_i[}(\,]x_i[\,)$$
Moreover, $\prod_{i=1}^m
\widehat{\mathcal{O}}_{\widetilde{\X},x_i}$ is canonically
isomorphic to the subring of $\mathcal{O}_{]x[}(\,]x[\,)$
consisting of the analytic functions $f$ on $]x[$ with supremum
norm $|f|_{sup}\leq 1$.
\end{prop}
\begin{proof}
The map $h_\eta:\widetilde{\X}_\eta\rightarrow \X_\eta$ is a
normalization map by \cite[2.1.3]{conrad}, and so is its
restriction over $]x[$, by \cite[1.2.3]{conrad}. Hence, $h_\eta$
is an isomorphism over $]x[$, so it is clear that $]x[\cong
\sqcup_{i=1}^m ]x_i[$. Therefore, we may as well assume that $\X$
is normal. In this case, the result follows from
\cite[7.4.1]{dj-formal}.   \end{proof}
\begin{cor}\label{equivalent}
Let $\X$ and $\mY$ be flat special formal $R$-schemes, and let $x$
and $y$ be closed points of $\X_0$, resp. $\mY_0$. Then $]x[$ and
$]y[$ are isomorphic over $K$, iff the disjoint unions of the
branches of $(\X,x)$, resp. $(\mY,y)$ are isomorphic over $R$. In
particular, if $\X$ and $\mY$ are normal at $x$, resp. $y$,  then
the $R$-algebras $\widehat{\mathcal{O}}_{\X,x}$ and
$\widehat{\mathcal{O}}_{\mY,y}$ are isomorphic iff $]x[$ and $]y[$
are isomorphic over $K$.
\end{cor}
\begin{lemma}\label{existsection}
Let $\X$ be a smooth special formal $R$-scheme, let $R'$ be any
finite unramified extension of $R$, and denote by $k'$ its residue
field. The natural map
$$\X(R')\rightarrow \X_0(k')$$ is surjective.
\end{lemma}
\begin{proof}
By formal smoothness, the map $\X(R'/\pi^{n+1})\rightarrow
\X(R'/\pi^n)$ is surjective for each $n\geq 0$, so since $R'$ is
complete, we see that $\X(R')\rightarrow \X_0(k')$ is surjective.
  \end{proof}

\begin{cor}\label{smoothcrit}
Consider a flat special formal $R$-scheme $\X$ and a point $x$ of
$\X_0(k)$. Then $\X$ is smooth at $x$, iff $]x[$ is isomorphic to
an open unit polydisc
$$\mathbb{B}^m_{K}=(\Spf\,R[[x_1,\ldots,x_m]])_\eta$$ for some $m\geq 0$.
\end{cor}
\begin{proof}
Replacing $\X$ by its formal completion at $x$, we may as well
assume that $\X_0=\{x\}$.
 If $\X$ is smooth at $x$, then $\X(R)$ is non-empty by Lemma \ref{existsection}, and
hence $\X\cong \Spf R[[x_1,\ldots,x_m]]$ by \cite[3.1/2]{neron}.
Hence, $]x[$ is isomorphic to the open unit polydisc
$\mathbb{B}^m_{K}$.
For the converse implication, assume that $\X_\eta\cong
\mathbb{B}^m_{K}$. Then $\X$ is normal since $\mathbb{B}^m_{K}$ is
normal and connected; so we can apply Corollary \ref{equivalent}
to $(\X,x)$ and $(\mathrm{Spf}\,R[[x_1,\ldots,x_m]],0)$.
\end{proof}
\subsection{The analytic Milnor fiber}
In this section, we put $R=k[[\pi]]$ and $K=k((\pi))$.
 Let $f:X\rightarrow \mathrm{Spec}\,k[\pi]$ be a
morphism from a $k$-variety $X$ to the affine line, let $x$ be a
closed point on the special fiber $X_s$ of $f$, and assume that
$f$ is flat at $x$. Denote by $\widehat{X}$ the $\pi$-adic
completion of $f$; it is a $stft$ formal $R$-scheme. The tube
$\mathscr{F}_x:=\,]x[$ of $x$ in $\widehat{X}$ is canonically
isomorphic to the generic fiber of the flat special formal
$R$-scheme $\mathrm{Spf}\,\widehat{\mathcal{O}}_{X,x}$ (the
$R$-structure being given by $f$), by \cite[0.2.7]{bert}. In
\cite{NiSe}, we called $\mathscr{F}_x$ the analytic Milnor fiber
of $f$ at $x$, based on a topological intuition explained in
\cite[4.1]{NiSe3} and a cohomological comparison result: if $k=\C$
and $X$ is smooth at $x$, then the \'etale $\ell$-adic cohomology
of $\mathscr{F}_x$ corresponds to the singular cohomology of the
classical topological Milnor fiber of $f$ at $x$, by
\cite[9.2]{NiSe}. If $f$ has smooth generic fiber (e.g. when
$X-X_s$ is smooth and $k$ has characteristic zero), then the
analytic Milnor fiber $\mathscr{F}_x$ of $f$ at $x$ is a smooth
rigid variety over $K$.

The arithmetic and geometric properties of $\mathscr{F}_x$ reflect
the nature of the singularity of $f$ at $x$ (see for instance
Proposition \ref{arithmilnor}). We will see in Proposition
\ref{milnor} that $\mathscr{F}_x$ is even a complete invariant of
the formal germ of the singularity $(f,x)$, if $X$ is normal at
$x$.

%
\begin{definition}
Let $X$ and $Y$ be $k$-varieties, endowed with $k$-morphisms
$f:X\rightarrow \mathrm{Spec}\,k[\pi]$ and $g:Y\rightarrow
\mathrm{Spec}\,k[\pi]$.
 We say that $(f,x)$
and $(g,y)$ are formally equivalent if
$\widehat{\mathcal{O}}_{X,x}$ and $\widehat{\mathcal{O}}_{Y,y}$
are isomorphic as $R$-algebras (the $R$-algebra structures being
given by $f$, resp. $g$).
\end{definition}

\begin{prop}\label{milnor}
Let $X$ and $Y$ be irreducible $k$-varieties, and let
$f:X\rightarrow \mathrm{Spec}\,k[\pi]$ and $g:Y\rightarrow
\mathrm{Spec}\,k[\pi]$ be dominant morphisms.
Let $x$ and $y$ be closed points on the special fibers $X_s$,
resp. $Y_s$, and assume that $X$ and $Y$ are normal at $x$, resp.
$y$. The analytic Milnor fibers $\mathscr{F}_x$ and
$\mathscr{F}_y$ of $f$ at $x$, resp. $g$ at $y$, are isomorphic
over $K$, iff $(f,x)$ and $(g,y)$ are formally equivalent.

More precisely, the completed local ring
$\widehat{\mathcal{O}}_{X,x}$ is recovered, as a $R$-algebra, as
the algebra of analytic functions $f$ on $\mathscr{F}_x$ with
$|f|_{sup}\leq 1$.
\end{prop}
\begin{proof}
By Proposition \ref{prop-branch}, it suffices to show that
$\widehat{\mathcal{O}}_{X,x}$ and $\widehat{\mathcal{O}}_{Y,y}$
are normal and flat. Normality follows from normality of
$\mathcal{O}_{X,x}$ and $\mathcal{O}_{Y,y}$,  by excellence, and
flatness follows from the fact that $f$ and $g$ are flat.
\end{proof}

\begin{prop}
Let $X$ be any $k$-variety, let $f:X\rightarrow
\mathrm{Spec}\,k[\pi]$ be a morphism of $k$-varieties, let $x$ be
a $k$-rational point on the special fiber $X_s$ of $f$, and assume
that $f$ is flat at $x$. Then $f$ is smooth at $x$ iff
$\mathscr{F}_x$ is isomorphic to an open unit polydisc
$\mathbb{B}^m_K$ for some $m\geq 0$.
\end{prop}
\begin{proof}
This follows from Corollary \ref{smoothcrit} (smoothness of $f$ at
$x$ is equivalent to smoothness of
$\mathrm{Spf}\,\widehat{\mathcal{O}}_{X,x}$ over $R$.
\end{proof}

\begin{prop}\label{arithmilnor}
Let $X$ be a smooth irreducible $k$-variety, let $f:X\rightarrow
\mathrm{Spec}\,k[\pi]$ be a dominant morphism, and let $x$ be a
closed point of $X_s$ whose residue field $k_x$ is separable over
$k$. The following are equivalent:
\begin{enumerate}
\item the morphism $f$ is smooth at $x$,
\item the analytic Milnor fiber $\mathscr{F}_x$ of $f$ at $x$
contains a $K'$-rational point for some finite unramified
extension $K'$ of $K$.
\end{enumerate}
If $k$ is perfect, then each of the above statements is also
equivalent to

\vspace{2pt}
 $\ \mathrm{(3)}\ $the analytic Milnor
fiber $\mathscr{F}_x$ of $f$ at $x$ is smooth over $K$ and
satisfies $S(\mathscr{F}_x)\neq 0$.

\vspace{5pt} If $k$ has characteristic zero, and if we denote by
$\varphi$ a topological generator of the geometric monodromy group
$G(K^s/K^{sh})$, then each of the above statements is also
equivalent to

$\ \mathrm{(4)}\ $the analytic Milnor fiber $\mathscr{F}_x$ of $f$
at $x$ satisfies
\begin{equation*}Tr(\varphi\,|\,H(\mathscr{F}_x\widehat{\times}_K\widehat{K^s},\Q_\ell))\neq
0\end{equation*}

\vspace{5pt} If $k=\C$, and if we denote by $F_x$ the classical
topological Milnor fiber of $f$ at $x$ and by $M$ the monodromy
transformation on the graded singular cohomology
$$H_{sing}(F_x,\C)=\oplus_{i\geq 0}H^i_{sing}(F_x,\C)$$ then each
of the above statements is also equivalent to

$\ \mathrm{(5)}\ $the topological Milnor fiber $F_x$ of $f$ at $x$
satisfies
$$Tr(M\,|\,H_{sing}(F_x,\C))\neq
0$$
\end{prop}
\begin{proof}
The implication $(1)\Rightarrow  (2)$ (with $K'=k_x((\pi))$)
follows from Lemma \ref{existsection}. The implication
$(2)\Rightarrow (1)$ follows from \cite[3.1/2]{neron}: denote by
$R'$ the normalization of $R$ in $K'$. Since
$X_R:=X\times_{k[\pi]}R$ is regular and $R$-flat, the existence of
a $R'$-section through $x$ on the $R$-scheme $X_R$ implies
smoothness of $X_R$ at $x$. But the set $\mathscr{F}_x(K')$ is
canonically bijective to the set of $R'$-sections on $X_R$ through
$x$.

The implication $(4)\Rightarrow (3)$ follows from the trace
formula (Theorem \ref{trace}), $(3)\Rightarrow (2)$ is obvious,
and $(1)\Rightarrow (4)$ follows from \cite[3.5]{berk-vanish2} and
the triviality of the $\ell$-adic nearby cycles of $f$ at $x$.
Finally, the equivalence $(4)\Leftrightarrow (5)$ follows from the
comparison theorem \cite[9.2]{NiSe}.   \end{proof} The equivalence
of $(1)$ and $(5)$ (for $k=\C$) is a classical result by A'Campo
\cite{Acampo2}.

\section{Comparison to the motivic zeta
function}\label{sectioncompar} We suppose that $k$ has
characteristic zero, and we put $R=k[[\pi]]$. Let $X$ be a smooth,
irreducible $k$-variety, of dimension $m$, and consider a dominant
morphism $f:X\rightarrow \mathrm{Spec}\,k[\pi]$. The formal
$\pi$-adic completion of $f$ is a generically smooth, flat $stft$
formal $R$-scheme $\widehat{X}$. We denote by $X_\eta$ its generic
fiber.

\begin{definition}
We call $X_\eta$ the rigid nearby fiber of the morphism $f$. It is
a separated, smooth, quasi-compact rigid variety over
$K=k((\pi))$.
\end{definition}

\subsection{The monodromy zeta function}
\begin{definition}
Suppose that $k$ is algebraically closed. For any locally closed
subset $V$ of $X_s$, we define the monodromy zeta function of $f$
at $V$ by
$$\zeta_{f,V}(T)=\prod_{i\geq 0}\det
(1-T\varphi\,|\,H^i(\,]V[\widehat{\times}_K
\widehat{K^s},\Q_{\ell}))^{(-1)^{i+1}}\in \Q_{\ell}[[T]]$$ where
$\varphi$ is a topological generator of the Galois group
$G(K^s/K)$.
\end{definition}
\begin{lemma}\label{comparecoh}
If $k=\C$, and $x$ is a closed point of $X_s$, then
$\zeta_{f,x}(T)$ coincides with
$$\zeta(M\,|\,\oplus_{i\geq 0}H^i_{sing}(F_x,\Q);T)$$ where $F_x$ denotes the
topological Milnor fiber of $f$ at $x$, and $M$ is the monodromy
transformation on the graded singular cohomology space
$\oplus_{i\geq 0}H^i_{sing}(F_x,\Q)$.
\end{lemma}
\begin{proof}
This follows from the comparison result in \cite[9.2]{NiSe}.
\end{proof}

For $k=\C$, the function $\zeta_{f,x}(T)$ is known as the
monodromy zeta function of $f$ at $x$.
\subsection{Denef and Loeser's motivic zeta functions}\label{DL}
 As in \cite[p.1]{DLinvent}, we denote, for any integer
$d>0$, by $\mathcal{L}_d(X)$ the $k$-scheme representing the
functor
$$(k-algebras)\rightarrow (Sets)\,:\,A\mapsto X(A[t]/(t^{d+1}))$$
Following \cite[3.2]{DL3}, we denote by $\mathcal{X}_d$ and
$\mathcal{X}_{d,1}$ the $X_s$-varieties
\begin{eqnarray*}
\mathcal{X}_{d}&:=&\{\psi\in
\mathcal{L}_{d}(X)\,|\,ord_tf(\psi(t))=d\} \\
\mathcal{X}_{d,1}&:=&\{\psi\in
\mathcal{L}_{d}(X)\,|\,f(\psi(t))=t^d\,mod\,t^{d+1}\}
\end{eqnarray*} where the structural morphisms to
$X_s$ are given by reduction modulo $t$.

In \cite[3.2.1]{DL3}, the motivic zeta function $Z(f;T)$ of $f$ is
defined as
$$Z(f;T)=\sum_{d=1}^{\infty}[\mathcal{X}_{d,1}]\LL^{-md}T^{d}\in \mathcal{M}_{X_s}[[T]]$$
and the \textit{na\"{i}ve} motivic zeta function
$Z^{na\ddot{\i}ve}(T)$ is defined as
$$Z^{na\ddot{\i}ve}(f;T)=\sum_{d=1}^{\infty}[\mathcal{X}_{d}]\LL^{-md}T^{d}\in \mathcal{M}_{X_s}[[T]]$$
If $U$ is a locally closed subscheme of $X_s$, the local
(na\"{i}ve) motivic zeta function $Z_U(f;T)$ (resp.
$Z^{na\ddot{\i}ve}_U(f;T)$) with support in $U$ is obtained by
applying the base change morphism
$\mathcal{M}_{X_s}[[T]]\rightarrow \mathcal{M}_U[[T]]$.

Let $h:X'\rightarrow X$ be an embedded resolution for $f$, with
$X'_s=\sum_{i\in I}N_i E_i$, and with Jacobian divisor
$K_{X'|X}=\sum_{i\in I}(\nu_i-1)E_i$.

By \cite{DL3}, Theorem 3.3.1, we have
\begin{eqnarray*}Z(f;T)&=&\sum_{\emptyset\neq J\subset
I}(\LL-1)^{|J|-1}[\widetilde{E}_J^{o}] \prod_{i\in
J}\frac{\LL^{-\nu_i}T^{N_i}}{1-\LL^{-\nu_i}T^{N_i}} \ \mathrm{in}\
\mathcal{M}_{X_s}[[T]]
\\Z^{na\ddot{\i}ve}(f;T)&=&\sum_{\emptyset\neq J\subset
I}(\LL-1)^{|J|}[E_J^{o}] \prod_{i\in
J}\frac{\LL^{-\nu_i}T^{N_i}}{1-\LL^{-\nu_i}T^{N_i}} \ \mathrm{in}\
\mathcal{M}_{X_s}[[T]]\end{eqnarray*}
 Inspired by the $p$-adic case \cite{Denef5},
Denef and Loeser defined the motivic
 nearby cycles $\mathcal{S}_f$ by taking formally the
limit of $-Z(f;T)$ for $T\to\infty$, i.e.
$$\mathcal{S}_f= \sum_{\emptyset\neq J\subset I}(1-\LL)^{|J|-1}[\widetilde{E}_J^{o}]\in \mathcal{M}_{X_s}$$
For each closed point $x$ of $X_s$, they denote by
$\mathcal{S}_{f,x}$ the image of $\mathcal{S}_f$ under the base
change morphism $\mathcal{M}_{X_s}\rightarrow \mathcal{M}_x$, and
they called $\mathcal{S}_{f,x}$ the motivic Milnor fiber of $f$ at
$x$. This terminology is justified by the fact that, when $k=\C$,
 the mixed Hodge structure of $\mathcal{S}_{f,x}\in \mathcal{M}_{\C}$
coincides with the mixed Hodge structure on the cohomology of the
topological Milnor fiber of $f$ at $x$ (in an appropriate
Grothendieck group of mixed Hodge structures); see
\cite[4.2]{DL5}.

\begin{theorem}
Let $x$ be a closed point of $X_s$. The local zeta functions
$Z_x(f;T)$ and $Z_x^{na\ddot{\i}ve}(f;T)$, and the motivic Milnor
fiber $\mathcal{S}_{f,x}$, depend only on the rigid $K$-variety
$\mathscr{F}_x$, the analytic Milnor fiber of $f$ at $x$.
\end{theorem}
\begin{proof}
This follows from Proposition \ref{milnor}, since all these
invariants can be computed on the $R$-algebra
$\widehat{\mathcal{O}}_{X,x}$, i.e. they are invariant under
formal equivalence. To see this, note that any arc
$\psi:\mathrm{Spec}\,k'[[t]]\rightarrow X$ with origin $x$ factors
through a morphism of $k$-algebras
$\widehat{\mathcal{O}}_{X,x}\rightarrow k'[[t]]$, and that
$f(\psi)\in k'[[t]]$ is simply the image of $\pi$ under the
composition $$\begin{CD}k[[\pi]]@>f>>
\widehat{\mathcal{O}}_{X,x}@>\psi>> k'[[t]]\end{CD}$$
\end{proof}
\begin{remark}
In fact, the local zeta function $Z_x(f;T)$ carries additional
structure, coming from a $\widehat{\mu}(k)$-action on the spaces
$\mathcal{X}_{d,1}$ (see \cite[3.2.1]{DL3}); the resulting
$\widehat{\mu}(k)$-action on the motivic Milnor fiber
$\mathcal{S}_{f,x}$ captures the semi-simple part of the monodromy
action on the cohomology of the topological Milnor fiber, by
\cite[3.5.5]{DL3}. The same argument as above shows that
$\mathscr{F}_x$ completely determines the zeta function
$Z_x(f;T)\in \mathcal{M}^{\widehat{\mu}(k)}_x[[T]]$ and the
motivic Milnor fiber $\mathcal{S}_{f,x}\in
\mathcal{M}^{\widehat{\mu}(k)}_x$, where
$\mathcal{M}^{\widehat{\mu}(k)}_x$ is the localized Grothendieck
ring of varieties over $x$ with good $\widehat{\mu}(k)$-action
\cite[2.4]{DL3}.   \end{remark}

In Corollary \ref{comparlocal} and Theorem \ref{cycles}, we'll
realize $Z_x(f;T)$ and $\mathcal{S}_{f,x}$ explicitly in terms of
the analytic Milnor fiber $\mathscr{F}_x$.
\subsection{Comparison to the motivic zeta
function}\label{comparz}
We define the local singular series
associated to $f$ by
$$F(f;d)=F(\widehat{X};d)\ \in \mathcal{M}_{X_s}$$
(see Definition \ref{localsing}).
 We define the motivic Weil generating series associated to $f$ by
$$S(f;T):=S(\widehat{X};T)=\sum_{d>0}F(f;d)T^d\ \in
\mathcal{M}_{X_s}[[T]]$$ For any locally closed subscheme $U$ of
$X_s$, we define the motivic Weil generating series $S_U(f;T)$
with support in $U$ as the image of $S(f;T)$ under the base-change
morphism
$$\mathcal{M}_{X_s}[[T]]\rightarrow \mathcal{M}_{U}[[T]]$$
If we denote by $\mU$ the formal completion of $\widehat{X}$ along
$U$, $S_U(f;T)$ coincides with $S(\mU;T)$ by
Proposition-Definition \ref{tube-volume}; in particular, it
depends only on $\mU$.

We recall the following result \cite[9.10]{NiSe}.
\begin{theorem}\label{comparzeta}
 We
 have
$$S(f;T)=\LL^{-(m-1)}Z(f;\LL T)\in
\mathcal{M}_{X_s}[[T]]$$
\end{theorem}
\begin{cor}\label{comparlocal}
For any closed point $x$ on $X_s$,
$$\LL^{-(m-1)}Z_x(f;\LL T)=
S(\mathrm{Spf}\,\widehat{\mathcal{O}}_{X,x};T)=\sum_{d>0}\left(\int_{\mathscr{F}_x(d)}|\omega/d\pi(d)|\right)T^d
\in \mathcal{M}_{x}[[T]]$$ where $\omega$ is any gauge form on
$\mathrm{Spf}\,\widehat{\mathcal{O}}_{X,x}$
 and where we
view $\mathscr{F}_x$ as a rigid variety over $k_x((\pi))$, with
$k_x$ the residue field of $x$.
\end{cor}
 Hence, modulo a normalization by powers of $\LL$, we recover the
motivic zeta function and the local motivic zeta function at $x$
as the Weil generating series of $\widehat{X}$, resp.
$\mathrm{Spf}\,\widehat{\mathcal{O}}_{X,x}$.
\subsection{Comparison to the motivic nearby cycles}
\begin{theorem}\label{cycles}
We have
$$S(\widehat{X};\widehat{K^s})=\LL^{-(m-1)}\mathcal{S}_f \in\mathcal{M}_{X_s}$$
For any closed point $x$ on $X_s$, we have
$$S(\mathscr{F}_x;\widehat{K^s})=\LL^{-(m-1)}\mathcal{S}_{f,x} \in\mathcal{M}_{x}$$
\end{theorem}
\begin{proof}
The result is obtained by taking a limit $T\rightarrow \infty$ of
the equality in Theorem \ref{comparzeta}, and applying
Proposition-Definition \ref{tube-volume}. Note that
$S(\mathscr{F}_x;\widehat{K^s})$ is well-defined, since
$\mathscr{F}_x$ admits a universally bounded gauge form by
Corollary \ref{gaugeform} and Lemma \ref{universally}.
\end{proof}

\begin{prop}
Assume that $k$ is algebraically closed, and let $\varphi$ be a
topological generator of the absolute Galois group $G(K^s /K)$.
For any integer $d>0$,
$$\chi_{top}S(\mathscr{F}_x(d))=Tr(\varphi^d\,|\,H(\mathscr{F}_x\widehat{\times}_K
\widehat{K^s},\Q_\ell))$$
\end{prop}
\begin{proof}
This is a special case of the trace formula in Theorem
\ref{trace}.   \end{proof}
\begin{cor}
Suppose $k=\C$, denote by $F_x$ the topological Milnor fiber of
$f$ at $x$, and by $M$ the monodromy transformation on the graded
singular cohomology space $H(F_x,\C)$. For any integer $d>0$,
$$\chi_{top}S(\mathscr{F}_x(d))=Tr(M^d\,|\,H(F_x,\C))$$
\end{cor}
\begin{proof}
This follows from the cohomological comparison in
\cite[9.2]{NiSe}.   \end{proof}

\bibliographystyle{hplain}
\bibliography{wanbib,wanbib2}
\end{document}